\documentclass[11pt]{amsart}

\usepackage{palatino} 
\usepackage{pinlabel}
\usepackage{pdfpages}
\usepackage{hyperref}

\renewcommand{\Im}{\operatorname{Im}}

\usepackage{amsmath,amssymb,amsthm,verbatim, amsfonts,amscd,flafter,epsf, epsfig,graphicx,verbatim,pinlabel,mathrsfs}
\usepackage[all]{xy}
\usepackage{epsf}
\usepackage[abs]{overpic}
\usepackage{epstopdf}
\usepackage{color}
\usepackage{mathtools}

\definecolor{red}{rgb}{1,0,0}
\definecolor{blue}{rgb}{0,0,1}
\definecolor{green}{rgb}{0,1,0}

\DeclareMathOperator{\lcm}{lcm}
\newtheorem{theorem}{Theorem}[section]
\newtheorem{lemma}[theorem]{Lemma}

\newtheorem{corollary}[theorem]{Corollary}
\newtheorem{proposition}[theorem]{Proposition}

\theoremstyle{definition}
\newtheorem{definition}[theorem]{Definition}

\newtheorem{remark}[theorem]{Remark}

\newtheorem{example}[theorem]{Example}

\def\F{\mathcal{F}}

\def\id{\mathit{id}}
\def\d{d}
\def\R{\mathbb{R}}
\def\pit{{\widetilde{\pi}}}

\def\E{\mathcal{E}}
\def\M{\mathcal{M}}
\def\A{\mathcal{A}}

   \title{Foliated Open Books}

\author[Vera V\'ertesi]{Vera V\'ertesi}
\address{IRMA \\ University of Vienna}
\email{vera.vertesi@univie.ac.at}

\author[Joan E. Licata]{Joan E. Licata}
\address{Mathematical Sciences Institute, The Australian National University}
\email{joan.licata@anu.edu.au}

\begin{document}
\begin{abstract}
This paper introduces a new type of open book decomposition for a contact three-manifold with a specified characteristic foliation $\F_\xi$ on its boundary.   These \textit{foliated open books} offer a finer tool for studying contact manifolds with convex boundary than existing models, as the boundary foliation carries more data than the dividing set.  In addition to establishing fundamental results about the  uniqueness and existence of foliated open books, we carefully examine their relationship with the partial open books introduced by Honda-Kazez-Matic. Foliated open books have user-friendly cutting and gluing properties, and they arise naturally as submanifolds of classical open books for closed three-manifolds. We define three versions of foliated open books (embedded, Morse, and abstract), and we prove the equivalence of these models as well as a Giroux Correspondence which characterizes the foliated open books associated to a fixed triple $(M, \xi, \F)$.

\end{abstract}
\maketitle

\maketitle




\section{Introduction}
The simplest way to produce a manifold with boundary is to cut a closed manifold, and one goal of this paper is to render this natural operation an effective one in the setting of contact geometry.  We introduce a new topological decomposition of a three-manifold with boundary: a foliated open book. The simplest construction of a foliated open book is an intuitive one: under mild hypotheses, cutting a closed  contact manifold equipped with an open book along a separating surface yields a pair of foliated open books. 

Open book techniques have been responsible for significant progress in contact geometry over the last two decades, leading to both computational and conceptual advances in the field. 
Open books were first used by Thurston and Winkelnkemper \cite{TW} to prove existence of contact structures, and Giroux \cite{Giob} upgraded them as a major tool in the field with the well-celebrated ``Giroux Correspondence'', claiming that contact structures  up to isotopy are described by open books up to positive stabilization. 
Open books were generalized for contact manifolds with convex boundary by Honda,Kazez, and Matic \cite{HKM}, and Van Horn-Morris defined another version
 of open books for contact manifolds with toroidal boundary.  In the latter, pages intersect the boundary torus in circles which give the characteristic foliation on the boundary \cite{Jeremy}. Foliated open books use a similar idea; we will require the characteristic foliation on the boundary to be ``the  same'' as the foliation induced by the intersection of the boundary with the pages, but we adapt the condition of Honda, Kazez, and Matic \cite{HKM} and require that the boundary is convex.

Our main tool for understanding the open book structure near the boundary comes from open book foliations; these originate in the thesis of Bennequin \cite{Be}, who used the singular foliation induced on a disc by the  
angular open book decomposition of $S^3$ to distinguish contact structures on $S^3$. These methods were later named braid foliations and were extensively used  by Birman and Menasco \cite{BM1,BM2,BM3,BM4,BM5,BM6}. Pavelscu revived the notion of braid foliations in general contact structures \cite{Pav} and they were further studied under the name open book foliations by Ito and Kawamuro \cite{IK1,IK2,IK3,IK4,IK5,IK6,IK7,IK8,IK9, IK10}.

Topologically, an open book decomposition identifies the complement of a link as a fibration, and this partition of the manifold into the binding and pages serves to localise the twisting of the contact planes near the binding.  By definition, a contact plane field is non-integrable, but away from the binding, the planes are nearly tangent to the pages: it becomes forgiveable to pretend the manifold is a collection of solid tori with twisting planes together with a foliated fiber bundle.  A foliated open book extends this fiction, decomposing a manifold with boundary into binding and pages.  In this setting, however, the pages need not have constant topology, but may evolve via saddle resolutions.  The idea of an open book adapted to a contact manifold with boundary is not new, but there are a few key differences between the partial open books of Honda-Kazez-Matic \cite{HKM} and the present foliated open books.  First, we allow arbitrarily many non-homeomorphic page types.  Second, the boundary of the manifold inherits a singular foliation   $\F_{\pit}$ whose leaves are the intersections of the pages with the boundary.  

Before continuing, we offer a first example of a foliated open book. 

\begin{example}\label{ex:torus}
Let $M$ be a solid torus embedded in $\mathbb{R}^3=\{r, \theta, z\}$ as shown in Figure~\ref{fig:efobtorus}.  The foliated open book structure on $M$ comes from decomposing the complement of the $z$-axis into surfaces defined by $M\cap \{\theta=c\}$.

\begin{figure}[h!]
\begin{center}
\includegraphics[scale=0.8]{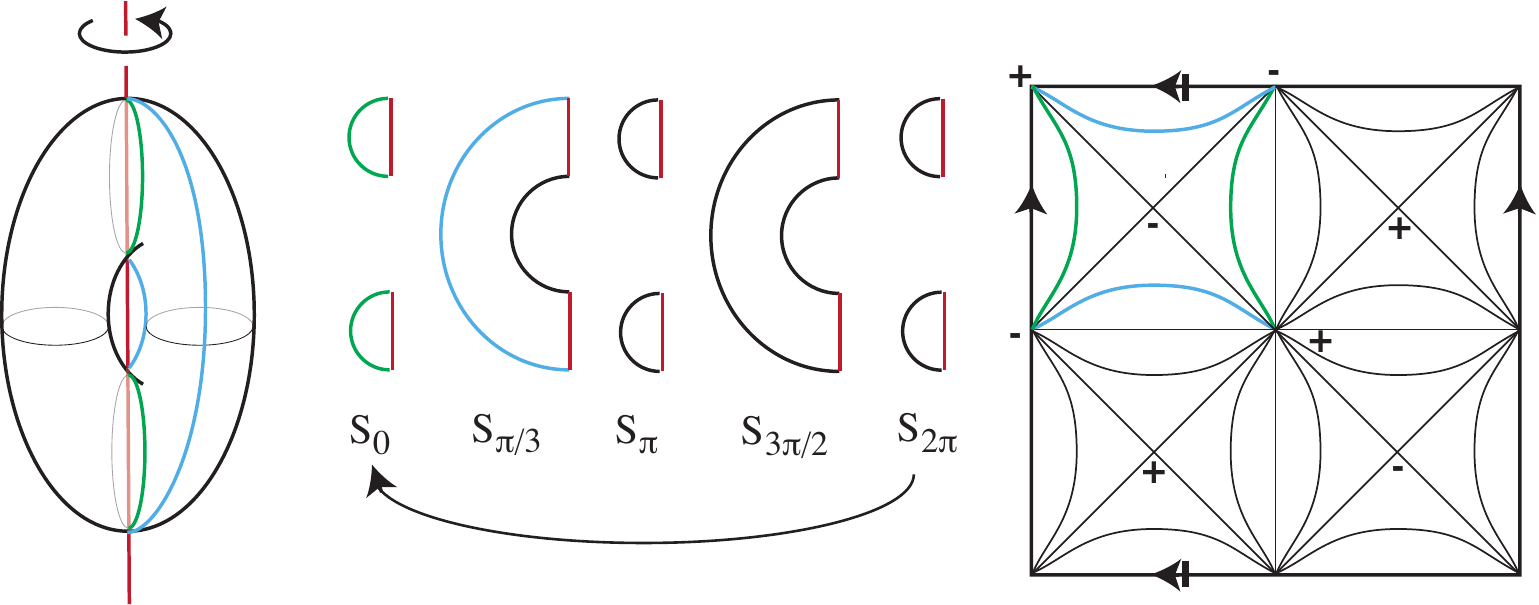}
\caption{ Left: $M=D^2\times S^1$ in $\mathbb{R}^3$.  Center: Selected pages $S_c= M\cap \{\theta=c\}$.  Right: singular foliation on $\partial M$. The  green and blue curves in the three pictures indicate the leaves $\theta=0$ and  $\theta=\frac{\pi}{3}$, respectively. }\label{fig:efobtorus}
\end{center}
\end{figure}

As $c$ changes, the topology of the surfaces changes four times, where each change either joins two components or splits one.  The boundary torus inherits a singular foliation from its intersections with the radial half-planes.  This foliation has elliptic singularities along $\partial M \cap \{r=0\}$ and a hyperbolic singularity corresponding to each saddle resolution.

\end{example}

The singular foliation on $\partial M$ is an intrinsic part of the foliated open book structure.  Singular foliations are ubiquitous in contact geometry, most prominently in the case of the characteristic foliation of a convex surface.  Above, we summarized the compatibility between an open book and a contact structure by requiring the contact planes to be ``nearly tangent'' to the pages, and this translates to foliated open books as a requirement that the characteristic foliation on the boundary ``nearly agrees" with the singular foliation induced by the pages.

A foliated open book $(B, \pi, \F_{\pit})$ is compatible with the contact structrure $\xi=\ker \alpha$ if 
\begin{itemize}
\item[-]$\alpha (TB)>0$;
\item[-] $d\alpha|_{\pi^{-1}(t)}$ is an area form;
\item[-] the singular foliation $\F_{\pit}$ on $\partial M$ whose leaves are the level sets of $\pit=\pi\vert_{\partial M}$ ``agrees with''  
 the characteristic foliation $\F_{\xi}$.
\end{itemize}

Precise definitions appear in Section \ref{sec:fob}.  
Having noted some contrasts between partial and foliated open books before the example, we now make the case that the world has room for yet another notion of an open book; this rests on claims of naturality (not in the categorical sense) and gluing.  As described above, when a contact manifold arises as an embedded submanifold,  the foliated open book structure is immediate from that of the ambient manifold.  This makes it easy to construct examples, whereas constructing partial open books can be difficult in practice even for rather simple manifolds.  Just as cutting is intuitive, so is gluing: keeping track of $\F_{\pit}$ on the boundary allows us to glue contact manifolds with foliated open books and get a new foliated open book.  Although gluing is certainly possible with partial open books, it is less straightforward.

\begin{theorem}[See Proposition \ref{prop:glue}]Suppose that the foliated open books compatible with $(M^L,\xi^L)$ and $(M^R,\xi^R)$ induce the same foliation along their boundary.
Then the contact 3-manifold $(M^L\cup M^R,\xi^L\cup\xi^R)$ formed by gluing them has a compatible honest open book decomposition that restricts to each piece as its original foliated open book.
\end{theorem}

\subsection{Applications and results}
The primary purpose  for defining foliated open books is for applications requiring  cut-and-paste arguments, and an example of this approach appears already in \cite{V}, which establishes the additivity of the support norm for tight contact structures. We anticipate further constructions of open books for manifolds assembled from smaller pieces as well as  applications to the study of support norms. 

Foliated open books are inherently compatible with gluing,  so by design they are the right object to define the contact invariant in bordered Floer homology.  We present this construction, together with its comparison to the gluing result of Honda-Kazez-Matic, in a forthcoming paper with Alishahi, F\"oldv\'ari, Hendricks, and Petkova \cite{AFHLPV}.   

It is also possible to amalgamate partial open books with foliated open books; a limited case of this appear in the proof of Theorem~\ref{thm:fobexistence}], and the general construction  will be explored in \cite{LV2}. This is significant not only as a  bridge between different models, but importantly, as a  tool for further work in Heegaard Floer homology. Specifically, these amalgamations offer a path to characterizing the higher multiplication maps in bordered Heegaard Floer homology using only the data of sutured Heegaard Floer invariants, as first proposed in \cite{Z}.

These advantages would be worth little, however, if foliated open books failed to be either sufficiently precise or sufficiently broad.  We prove the following results:

\begin{theorem}[See Theorems \ref{thm:existencecontact}, \ref{thm:uniquecontact} and \ref{thm:fobexistence}] Every foliated open book supports a unique isotopy class of contact structure, and every contact manifold with a characteristic foliation that ``agrees with'' an open book foliation 
admits a compatible foliated open book.
\end{theorem}

As in the case of honest and partial open books, foliated open books may be stabilized by taking a connected sum with an open book for $S^3$.  

\begin{theorem}[Giroux Correspondence, see Theorem \ref{thm:giroux}] Positively stabilizing a foliated open book preserves its compatibility with a contact structure, and two foliated open books supporting the same contact structure are related by a sequence of positive (de)stabilizations. 
\end{theorem}

Having said something about what foliated open books do, we return to the somewhat neglected question of what foliated open books are.  In fact, we will introduce three distinct objects: embedded foliated open books, Morse foliated open books, and abstract foliated open books. In brief:
\begin{itemize}
\item[-] An \emph{embedded foliated open book} is a pair $(B, \pi)$ which could be the restriction of an honest open book to the submanifold formed by cutting  along a surface with an open book foliation without circle leaves.  See Definition~\ref{def:efob}.

\item[-] A \emph{Morse foliated open book} is a pair $(B,\pi)$ such that $B$ is a properly embedded $1$-manifold and $\pi:M\setminus B\rightarrow S^1$ is a circle-valued Morse function with only boundary critical points. See Definition~\ref{def:mfob}.
\item[-] An \emph{abstract open book} is a tuple $(\{S_i\}, h)$ where $S_{i+1}$ is built from $S_i$ by cutting along a properly embedded arc or by adding a one-handle, and $h:S_{2k}\rightarrow S_0$ is a homeomorphism. See  Definition~\ref{def:afob}.
\end{itemize}

In Section~\ref{sec:fobequiv} we describe how to transform a foliated open book of one flavor into another; having several versions provides technical flexibility in different  settings, as well as more transparency in the relationship to other open books.  

The final section of the paper is devoted to the connections with partial open books.  We characterise when a foliated open book naturally gives rise to a partial open book for the same manifold, and we show these conditions can always be achieved via stabilization.  Furthermore, we show how to turn a partial open book into a foliated open book; as one would suspect from the more detailed data on the boundary of a foliated open book, a given partial open book may give rise to several non-equivalent foliated open books.

\begin{theorem}[See Propositions \ref{prop:fobpobequiv} and \ref{prop:pobfobequiv}]
Any sufficiently positively stabilized foliated open book for $(M, \xi)$ contains a partial open book for $(M, \xi)$ as a submanifold  with page-wise inclusions. Furthermore, any sufficiently positively stabilized  partial open book may be obtained this way. 
\end{theorem}

\subsection{Organization of the Paper} 
The paper assumes familiarity with contact structures, open books and Morse functions, but we found it important to recall the basic terminology and results about open book foliations and characteristic foliations in Section \ref{sec:foliation}. Here we also introduce a technical definition of a  ``preferred gradient-like vector field'' that will be used throughout the paper. In Section \ref{sec:fob} we introduce the three definitions of foliated open books, and after setting up the local models in Section \ref{sec:models}, we prove their equivalence in Section \ref{sec:fobequiv}.  The last part of this section again introduces a technical notion: ``sorted handlebodies''  are the main ingredients in understanding the relationship between partial and foliated open books. 
The main theorems that make foliated open books useful are discussed in Section \ref{sec:operations}, together with a number of examples. We prove the existence and uniqueness of the supported contact structures in Section \ref{sec:existence}. Section \ref{sec:pob}  explains the relationship between  foliated open books and partial open books, and we use  this  to prove the existence of a supporting foliated open book and  a ``Giroux correspondence''. We close the paper with some brief remarks about possible applications of foliated open books in  Section \ref{sec:future}.

\section*{Acknowledgments}

The  authors  would  like  to  thank  the  Isaac  Newton  Institute for Mathematical Sciences, Cambridge, for  support and hospitality  during  the  programme  Homology Theories in Low Dimensional Topology where work on this paper was undertaken. This work was supported by EPSRC grant number  EP$\slash$K032208$\slash$1.

The second author would also like to acknowledge the generosity of the Mathematical Sciences Institute at the Australian National University for hosting her for 3 months through the LIA/CNRS Visitors program. She was also supported by the ANR grant ``Quantact''. She would also like to express her gratitude to the CNRS and the University of Strasbourg, where she was employed while most of the work was carried out.



\section{Foliations and vector fields on surfaces}\label{sec:foliation}

Although we assume familiarity with many of the standard tools of contact geometry, this section carefully introduces two classical singular foliations of surfaces in contact geometry: characteristic foliations and open book foliations.  

\subsection{Signed foliations}\label{sec:sgnfoliation}
\begin{definition}
An \emph{oriented singular foliation} on a surface $\Sigma$  is an equivalence class of smooth vector fields $X$  on $\Sigma$, where two vector fields are equivalent if they differ by multiplication by a smooth positive function. Zeroes of $X$ are the \emph{singular points} and connected components of integral curves of $X$ are called \emph{leaves}. We denote the equivalence class by $\F=[X]$.
  \end{definition}

\begin{figure}[h]
\begin{center}
\includegraphics[scale=0.6]{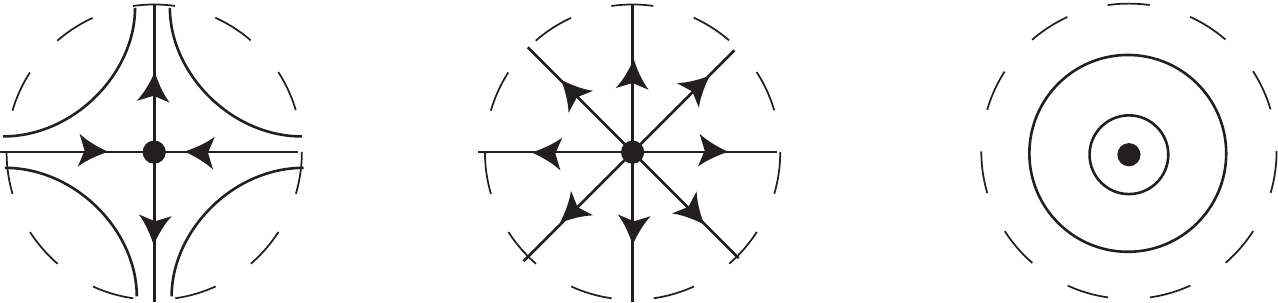}
\caption{Left: hyperbolic point.  Center: elliptic point (source).  Right: center.}\label{fig:sing}
\end{center}
\end{figure}

We will often refer to an oriented singular foliation as simply a foliation.  Both the zeroes and leaves of a foliation are independent of the representative of $\F$.  We will restrict to singular foliations with isolated singularities that are either \emph{centers} ($C$),   
four-prong \emph{saddles} or \emph{hyperbolic points} ($H$),  
or \emph{elliptic points} ($E$). See Figure \ref{fig:sing}. 

\begin{figure}[h]
\begin{center}
\includegraphics[scale=0.6]{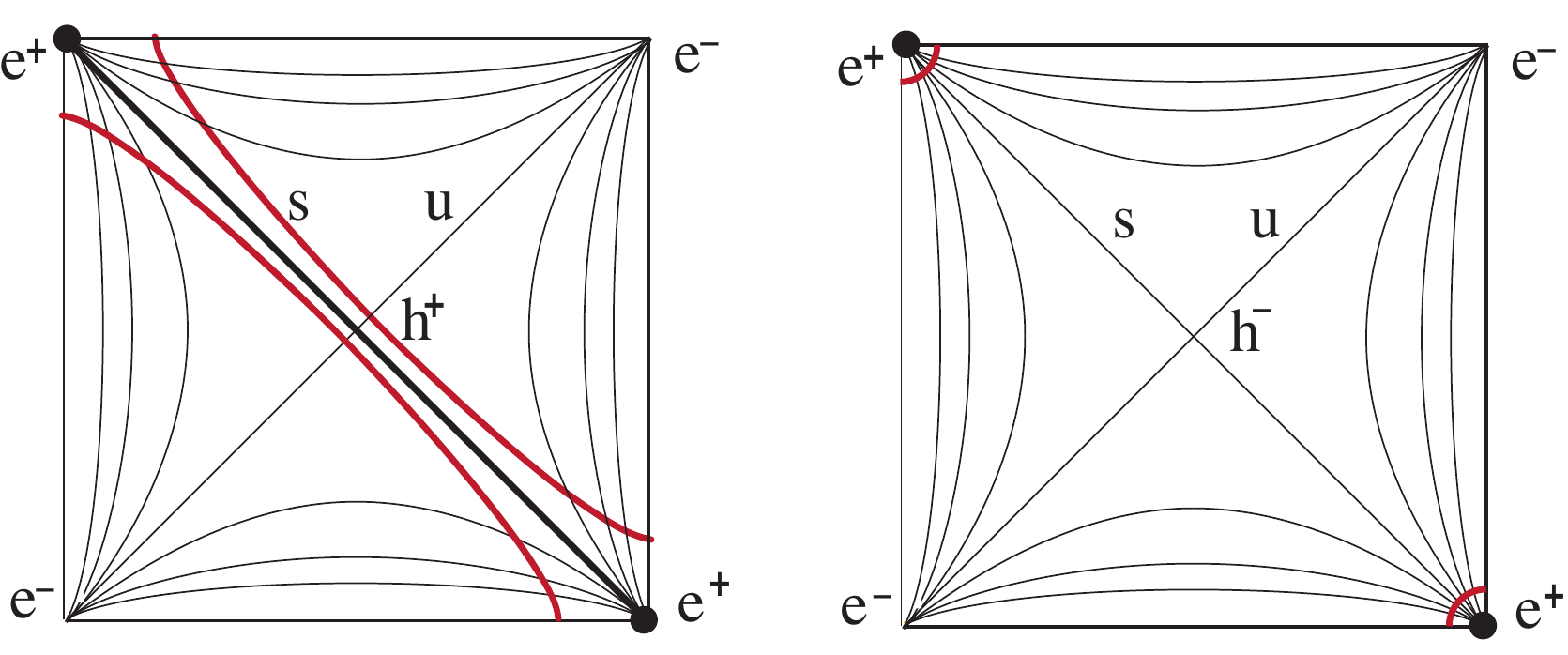}
\caption{The bold curves show $G_{\scriptscriptstyle{++}}$ on a tile defined by a positive (left) and negative (right) hyperbolic point, and the red neighborhood is $\Gamma$.  That the separatrices are labelled with $u$ and $s$ for unstable and stable, respectively.}\label{fig:g++}
\end{center}
\end{figure}

 Oriented singular foliations  have two types of elliptic points: \emph{sources} ($E_+$) and \emph{sinks} ($E_-$).  
The leaves that limit to (from) hyperbolic points are called \emph{(un)stable separatrices}.  Very often, we will further require that our foliation has no circle leaves, and thus, no centers. 
A singular foliation is \emph{Morse-Smale} if there are no leaves which connect a hyperbolic point to a hyperbolic point. In addition to the above data, we will often assign signs (a priori, arbitrarily) to the hyperbolic points.  When $|H_+|=|H_-|$, we call a Morse-Smale foliation with this extra information $(\F,H=H_+\cup H_-)$ a \emph{signed foliation}.  A signed foliation with no circle leaves can be cut along regular leaves into square \textit{tiles}, each of which contains exactly one hyperbolic singularity; see Figure~\ref{fig:g++}.

For signed foliations with no circle leaves we can define the \emph{positive graph} $G_{\scriptscriptstyle{++}}$ embedded into $\Sigma$ as the closure of the union of stable separatrices of positive hyperbolic points. This is a graph whose vertices are 
 the positive elliptic points and whose edges are in one-to-one correspondence with $H_+$; each edge connects the two elliptic points at the ends of the stable separatrices of a single positive hyperbolic point. 
 See Figure~\ref{fig:g++}.  Similarly, one can define the negative graph $G_{\scriptscriptstyle{--}}$ using unstable separatrices and negative hyperbolic points. A \emph{dividing curve} for a signed foliation $(\F,H=H_+\cup H_-)$ is a curve {which is positively transverse to the leaves of $\F$ and bounds} $R_+=N(G_{\scriptscriptstyle{++}})$. Up to isotopy through curves with this property, we could have defined $\Gamma$ as the boundary of $R_-=N(G_{\scriptscriptstyle{--}})$ with the opposite orientation. In fact, up to isotopy through curves  positively transverse to $\F$, $\Gamma$ is the unique curve that is positively transverse to $\F$ and separates $H_+$ from $H_-$. We say that $\Gamma$ \emph{divides} the signed foliation  $(\F,H=H_+\cup H_-)$.

\begin{definition}  Two signed foliations $(\F,H=H_+\cup H_-)$ and $(\F',H'=H'_+\cup H'_-)$ on $\Sigma$ are \emph{ topologically conjugate}  if there is a homeomorphism $\psi$ of $\Sigma$ that takes the leaves of $\F$ to the leaves of $\F'$, while respecting the partition of the hyperbolic points by sign (i.e., $\psi(H_\pm)=H_\pm'$). Two signed foliations are \emph{strongly topologically conjugate} if they are topologically conjugate via a homeomorphism in the identity mapping class of $\Sigma$ and there is a common curve $\Gamma$ on $\Sigma$ that separates both $\F$ and $\F'$.
\end{definition}

Although we have chosen to define a foliation as a vector field, the presence of an area form $\omega$ on $\Sigma$ offers an equivalent definition as a dual $1$-form.  Specifically, 
one may define an oriented singular foliation as the kernel of a 1-form $\beta$, where $\iota_X\omega=\beta$.  In computations, it is often more convenient to use 1-forms, so we will assume a fixed $\omega$ and move freely between these two perspectives.

\subsection{Characteristic foliation}\label{ssec:charfol} 
Characteristic foliations provide natural examples of signed foliations, and here we recall some results relevant to our paper. See \cite{Girbif} for further discussion.  

\begin{definition}
Let $j:\Sigma\hookrightarrow M$ be a smooth embedding of a surface $\Sigma$ into a contact 3--manifold $(M,\xi)$, and let $\alpha$ be some contact form satisfying $\ker \alpha=\xi$. The \emph{characteristic foliation} 
$\F_\xi(\Sigma)=\F_\xi$ is the kernel of the pullback $j^*\alpha\vert_\Sigma$. 
\end{definition}

As noted above, $\F_\xi$ can equivalently be defined by the vector field $X$ satisfying $\iota_X\omega=\alpha\vert_\Sigma$.  Note that with either definition, $\F_\xi$ is well defined only up to multiplication by a positive smooth function on $\Sigma$.  
The leaves of the characteristic foliation are immersed Legendrian curves on $\Sigma$, but the definition presented above is stronger than simply identifying the leaves as sets.

At regular points, the leaves of $\F_\xi$ are tangent to $T_p\Sigma\cap\xi_p$, oriented so that the vectors evaluating positively under $\alpha$ on $T_p\Sigma$ co-orient the leaves. 
The singular points are exactly the points where $T_p\Sigma=\xi_p$. A singular point is positive (respectively, negative) if the orientations of $T_p\Sigma$ and $\xi_p$ agree (disagree). We will see below that this sign agrees with the sign already defined for elliptic points, and that characteristic foliations have no centers. Depending on the sign of a singular point
 $\pm d(\alpha\vert_\Sigma)= d(\alpha\vert_\xi)=d\alpha\vert_{\xi}$; this is non-degenerate, as $\alpha$ is a contact form. This means, again depending on the sign of the singular point, that $\pm d(\alpha\vert_\Sigma)$ is an area form on $\Sigma$. Thus the isolated singular points of $\F_\xi$ can only be elliptic or hyperbolic. The sign of elliptic points is positive for sources and negative for sinks. Moreover, this sign convention for  hyperbolic points automatically makes characteristic foliations signed. 
 
One may recognize the 1-forms that arise as characteristic foliations: 
 
\begin{theorem}[Giroux, \cite{Girbif}] A 1-form $\beta$ on a surface $\Sigma$ is the restriction of some contact form for some embedding of $\Sigma$ in some contact manifold if and only if $\pm d\beta$ is an area form at the singular points.   \qed
\end{theorem}

We distinguish an important class of surfaces in a contact manifold: those whose neighborhoods  have $I$-invariant contact structures. A \emph{contact vector field} is a vector field ${X}$ whose flow preserves the contact structure $\xi$,  and an embedded surface $\Sigma$ is \emph{convex} if there is a contact vector field  transverse to it. Convex surfaces are naturally equipped with \emph{dividing curves} $\Gamma=\{p: {X}_p\in \xi_p\}$ and in fact, the existence of dividing curves for a characteristic foliation detects convexity: 
 
\begin{theorem}[Giroux \cite{Gi}]\label{thm:Gir} An embedded surface $\Sigma$ is convex if and only if there is a curve $\Gamma$ that divides $\beta=\alpha\vert_\Sigma$. \qed \end{theorem}

The usual definition of  the ``dividing curve''  of $\F_\xi$ in contact geometry is at first glance slightly stronger and begins with a fixed vector field $X$ which represents $\F_\xi$.  We require that  $X$ is positively transverse to $\Gamma$; that there exists an area form $\Omega$ such that 
$\Gamma=\{\mathit{div}_\Omega X=0\}$; and that $R_\pm$ is  $\{\pm \mathit{div}_\Omega X>0\}$. We will see in the following that the two notions for dividing are  equivalent, but until then we denote the  stronger, divergence-dependent definition by ``geometrically divide'', and our original, {$X_p\in \xi_p$} definition by ``topologically divide''.

In order to translate the definition of geometric dividing curves to the world of 1-forms, fix $X$ and choose  $\beta'=\iota_X\Omega$, a form which defines the same foliation $\F_\xi$.  Then the characteristic foliation $\F_\xi$ is geometrically divided by $\Gamma$ if and only if  $\Gamma=\{d\beta'=0\}$ and  $\beta'\vert_\Gamma$ orients $\Gamma$. In this case $R_\pm$ is defined by $\{\pm d\beta'>0\}$. The equivalence of the two notions of dividing now follows from the following lemma:

\begin{lemma}\label{lem:beta} If the characteristic foliation $\F_\xi$ is topologically divided by $\Gamma$, then it has a representative 1-form $\beta'$ so that $\Gamma=\{d\beta'=0\}$.
\end{lemma}
\begin{proof}
Though this statement is standard in contact geometry, we recall the proof for later use. Starting from an arbitrary 1-form $\beta$ defining $\F_\xi$, we will look for
$\beta'$ in the form $\beta'=g\beta$, for some positive function $g$. Let $A\cong [-\varepsilon,\varepsilon] \times \Gamma$ be a neighbourhood of $\Gamma$ with coordinates $(u,v)$,  so that $\F_\xi$ is  given by $\partial u$ and $\Gamma=\{u=0\}$.
 
Choose a Morse function $h$ on $\Sigma$ whose gradient flow with respect to some metric  directs $\F_\beta$ and that satisfies the following conditions:

\begin{itemize} 
\item[-] $\Gamma=h^{-1}(0)$;
\item[-] $A=h^{-1}[-\varepsilon,\varepsilon]$, so that $h=u$;
\item[-] $\partial u$ orients the level sets of $h$. 
\end{itemize}
Such a Morse function exists by Theorem B of \cite{Smale}, and we  use it to define
 $g=h^2+1-2\varepsilon$.  Assuming that $\varepsilon<1/2$, the function $g$ is indeed positive on $\Sigma$, $\beta'(T\Gamma)=g\beta(\partial u)>0$, and 
\[\pm d\beta'=\pm dg\wedge\beta\pm gd\beta.\]
The second term  on the right vanishes away the singular points and is a positive multiple of a volume form near the singular points. As for the first term, we have $dg=2hdh$.  Recall that the gradient vector field of $h$ directs $\F_\beta$, so $dh\wedge\beta>0$ and the sign of $h$ is $\pm$ in $R_\pm$. 
This proves that $\pm d\beta'>0$ everywhere on $R_\pm$ as required. 
\end{proof}

\begin{remark}\label{rmk:g}  Notice that in the proof above, $g=h^2+1-2\varepsilon\le 1-\varepsilon<1$ on $A$. We will  use this in Lemma~\ref{lem:approx}.
\end{remark}

From now on we will generally not distinguish the two notions of dividing, although this may require us to change the defining 1-form in order to  assume that some $\Gamma$ is geometrically dividing. For example,  Lemma \ref{lem:approx} will require a specific choice of 1-form.

One can also give a local model for the contact structure in the neighbourhood of $\Sigma$ as follows.  
As above, let $A=N(\Gamma)\cong [-\varepsilon,\varepsilon]\times \Gamma$ with the given coordinates. 
Then the contact structure on $\Sigma\times I(z)$ is given as:
\begin{proposition}[Giroux, \cite{Gi}]\label{prop:Gi} 

Suppose that $\F_\xi$ is (geometrically) divided by $\Gamma$ with respect to the representative $\beta'$.  
Then there is a function $f\colon \Sigma\to [-1,1]$ such that $f^{-1}(\pm 1)=R_\pm'$; on $A$, $f$ depends only on $u$ and  is a monotonically decreasing function of $u$;  $f^{-1}(0)=\Gamma$; and $\alpha=\beta'+fd z$ is a contact form on $\Sigma\times I$. 
\qed
\end{proposition}

\subsection{Open book foliations}\label{sec:obfol} There is yet another foliation related to contact structures,  the \emph{open book foliation} introduced first in the context of the angular open book for the standard contact structure in $S^3$ in \cite{Be} and then later generalized in 
 \cite{Pav} and \cite{IK1}.
Let $\Sigma$ be a surface, this time embedded in a manifold $M$ equipped with an open book $(B,\pi)$.
After a $C^\infty$-isotopy of $\Sigma$,  one may assume  that $\Sigma$ is transverse to $B$.  Define $E=B\cap \Sigma$ to be the set of \emph{elliptic points}.   After possibly applying a further $C^\infty$-isotopy of $\Sigma$,   we may assume  that $\widetilde{\pi}=\pi\vert_\Sigma\colon \Sigma\!\smallsetminus\! E\to S^1$ is a circle-valued Morse function with at most one critical point on each level.  For any such surface, the \emph{open book foliation} is defined in   \cite{IK1} as the level sets of $\pit$. In order to align this  with our use of the word foliation, we will define an equivalence class of vector fields whose integral curves are the level sets of $\pit$.

Consider the 1-form $d\pit$ on $\Sigma\setminus E$ and recall that an elliptic point $e\in E_\pm$ has a neighborhood $D^2_{\varepsilon}(r,\vartheta)$ where $\pit=\pm\vartheta$. Define $\gamma=\pm\Psi^{\varepsilon/2}_{\varepsilon/4}(r^2-1)d {\pit}+d\pit$, where $\Psi^{\varepsilon/2}_{\varepsilon/4}$ is a smooth bump function that equals $1$ for $r<\varepsilon/4$, equals $0$ for $r>\varepsilon/2$, and is strictly monotonically decreasing on $[\varepsilon/4,\varepsilon/2]$. Extend $\gamma$ by $0$ on $e$.

Up to multiplication by a smooth positive function, ${\gamma}$ does not depend on any of the above choices. 
The signs of the elliptic points of $\F_{\pit}$ come from the signs of the intersections of the oriented binding $B$ with $\Sigma$, while each hyperbolic point inherits a sign when we compare the orientation of the level sets of $\pi$ with the orientation of $\Sigma$.  Thus we recover a signed foliation, but with further structure coming from the map  $\pit$.  We summarise the relationship between this data in the next definition.
 
 \begin{definition}\label{def:obfol} An  \emph{open book foliation} on an oriented surface $\Sigma$ is a signed foliation $(\F=[\gamma], H=H_+\cup H_-)$ and a function $\pit$ that satisfy the following  conditions:
 \begin{enumerate}
 \item letting $E$ denote the set of elliptic singularities of $\F$, $\pit: \Sigma\setminus E\rightarrow S^1$ is a circle-valued Morse function with only index 1 critical points exactly at $H$;
 \item on each connected component of $\Sigma$, there is at most one critical point for each critical value of $\pit$;
 \item $\gamma=d\pit$  on $\Sigma\setminus N(E)$, and near each elliptic point, $\gamma$ is as described above;
  \item the partition of $E$ into $E_+\cup E_-$ is induced by  the orientation of $\Sigma$ and $M$;
   \end{enumerate}
 We denote an open book foliation on $\Sigma$ by $(\F_{\pit},\pit,H=H_+\cup H_-)$.
 \end{definition}

Note that  if a surface is embedded in a manifold with an open book, then it inherits an open book foliation, but this definition allows us to define an open book foliation on an abstract surface.  This is justified by the following somewhat vague statement, a more precise version of which is proved as  Propositions~\ref{prop:nbrh} and \ref{prop:obcyl}.
\begin{theorem}
Suppose that  $(\F_{\pit},\pit,H=H_+\cup H_-)$ is an open book foliation on $\Sigma$.  Then there exists a  manifold $M$ with open book decomposition $(B,\pi)$ such that $\Sigma$ embeds into $M$ and the inherited open book foliation on $\Sigma$ agrees with the original one.  Furthermore, $\pi$ is  determined up to diffeomorphism on a neighborhood of $\Sigma$ by the requirement that $\pi|_\Sigma$ induces  $(\F_{\pit}, \pit, H=H_-\cup H_+)$. 
\end{theorem}

Restricting  the open book $(B,\pi)$ to  $N(\Sigma)$ of $\Sigma$ produces an explicit local model for a product neighborhood of $\Sigma$; see Corollary~\ref{cor:lm}.  Whenever a surface is already embedded in a manifold with an open book, however, we assume that $(\F_{\pit},\pit,H=H_+\cup H_-)$ is the induced open book foliation.

A relationship between open book and characteristic foliations is proved in \cite{IK1}:
\begin{proposition}[Theorem 2.21, \cite{IK1}]\label{prop:IK} If $\Sigma$ is a surface in $(B, \pi)$ with  open book foliation $\F_{\pit}$ with no circle leaves, then there is a contact structure $\xi$ supported by $(B,\pi)$ such that $\F_{\pit}$ is strongly topologically conjugate to $\mathcal{F}_\xi$.  \qed
\end{proposition}  
\begin{remark}The original proof of \cite{IK1} claims only topological conjugacy, but in fact the stronger condition follows from their argument. 
\end{remark}
\begin{remark} For an open book foliation, $d(d\pit)=0$ everywhere, while for characteristic foliations {$\F_\beta$, $d\beta\neq 0$ at the hyperbolic points.}  Since we define foliations as 1-forms (or equivalently, vector fields), this is the strongest compatibility of foliations that one might hope for; contrast this to Remark 2.22 in \cite{IK1}, where foliations are viewed only as collections of leaves.
\end{remark}

\begin{example}The existence of the function $\pit$ is  essential  for open book foliations and does not follow from the other conditions.  Figure \ref{fig:notord} shows the singular leaves of a signed foliation which cannot be induced  as the level sets of an $S^1$-valued function. Two separatrices from the shaded hyperbolic point terminate at the shaded elliptic point, while only one separatrix from the unshaded hyperbolic point terminates at the shaded elliptic point;  this cannot occur in an open book foliation. To construct such an example on a closed surface, double the   annulus to yield a signed foliation on a torus.)

\begin{figure}[h]
\begin{center}
\includegraphics[scale=0.7]{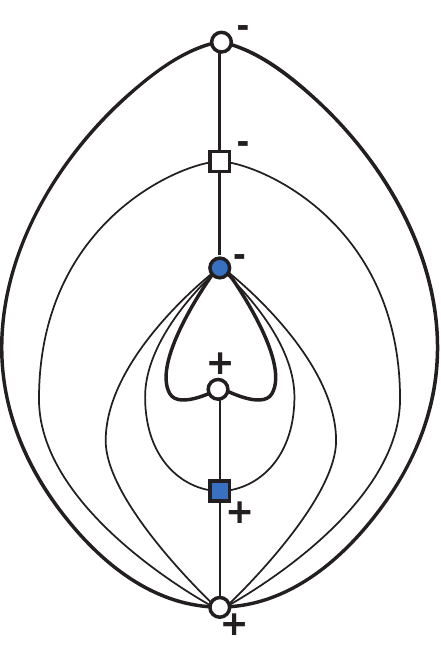}
\caption{In this signed foliation on the annulus, boxes represent hyperbolic points and circles represent elliptic points. }\label{fig:notord}
\end{center}
\end{figure}
\end{example}

\subsection{Gradient-like vector fields}\label{sec:grad}
When the leaves of a singular foliation arise as the level sets of a function, the gradient vector field with respect to some metric has the same singular points but is otherwise transverse.  In the absence of a designated metric, one may consider \textit{gradient-like} vector fields, which share many of the same properties.  
More formally, a vector field $X$ on a cobordism $M$ 
is gradient-like for $\pi$ \ if $X$ is positively transverse to the level sets of $\pi$ away from the critical points; $X$ is tangent to the vertical boundary; and $X$ can be described via the standard Morse model near the critical points.  See \cite{BNR} for details. The Morse function $\pi$ has a such a gradient-like vector field if and only if the critical points of $\pi$ agree with the critical points of $\pi\vert_{\partial M}$.

For a fixed gradient-like vector field defined on a manifold with boundary, let $W^s(h)$ (respectively, $W^u(h)$) denote the (un)stable submanifold of a critical
point $h$.  If the critical point $h$ lies on the boundary, let $w^s(h)$  and $w^u(h)$ denotes the critical submanifolds of the restriction of the gradient-like vector field to $\partial M$. Equivalently, $w^s(h)=W^s(h) \cap \partial M$ and  $w^u(h)=W^u(h) \cap \partial M$.

Given an oriented surface $\Sigma$ with a circle-free open book foliation $(\F_{\pit}, \pit, H=H_+\cup H_-)$, we define a special class of gradient-like vector fields for $\pit$ which are characterized by the relationship between their flowlines and the leaves of $\F_{\pit}$.  These will play an important role in Section~\ref{sec:disj}.  We begin by setting some notation and constructing an example, before stating Definition~\ref{def:pref} at the end of the section. 

\begin{figure}[h]
\centering
  \labellist
         \pinlabel  $e_+$ at 200 180
           \pinlabel  $e_-'$ at -20 180
            \pinlabel  $e_-$ at 200 -10
             \pinlabel  $e_+'$ at -20 -10
            \pinlabel {\textcolor{magenta}{$w^u$}} at 200 40
             \pinlabel {\textcolor{magenta}{$w^u$}} at -20 140
              \pinlabel {\textcolor{magenta}{$w^s$}} at  40 -15
               \pinlabel {\textcolor{magenta}{$w^s$}} at 160 195
                \pinlabel  $h$ at 100 100
                  \endlabellist
\includegraphics[scale=0.6]{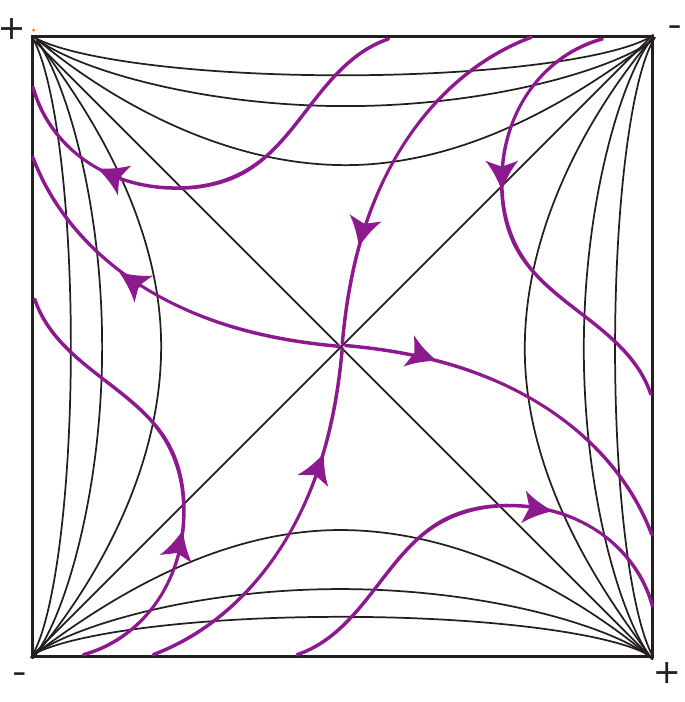}
\caption{The flow of $\nabla\pit$ on a tile.}\label{fig:grad}

\end{figure}

As noted in Section~\ref{sec:sgnfoliation}, the surface $\Sigma$ naturally decomposes into quadrilateral tiles, each of which contains precisely one hyperbolic point and whose boundary consists of regular leaves connecting elliptic point corners.  Let $\nabla\pit$ be a vector field whose flowlines on such a tile  are shown  in Figure \ref{fig:grad}.  (Since we are only interested in the qualitative behavior of the flowlines near the critical points, we are free to assume the required coordinate model.)

More precisely, $\nabla\pit$ vanishes exactly at the singular points of $\F_{\pit}$, is positively transverse to the level sets of $\pit$, and rotates infinitely many times positively (respectively, negatively) near the positive (negative) elliptic points.  By construction, $w^s(h)$ and $w^u(h)$ connect the same pairs of elliptic points as did the separatrices of $h$ in  $\F_{\pit}$. Away from the elliptic points, the flowlines of $\nabla \pit$ are ``close'' to the level sets of $\pit$ in the sense described next.

For ``closeness'', we further require that there is a small $\varepsilon>0$ and a small neighbourhood $N$ of the elliptic points so that for any hyperbolic point $h$ with $\pit(h)=t_h$, 
\begin{itemize}
\item[-] $w^s(h)\setminus N$ intersects  $\pit^{-1}(t)$ only for $t\in [t_h-\varepsilon,t_h]$;
\item[-] $w^u(h)\setminus N$ intersects  $\pit^{-1}(t)$ only for $t\in [t_h, t_h+\varepsilon]$. 
\end{itemize}

\begin{figure}[h]
\begin{center}
\includegraphics[scale=.6]{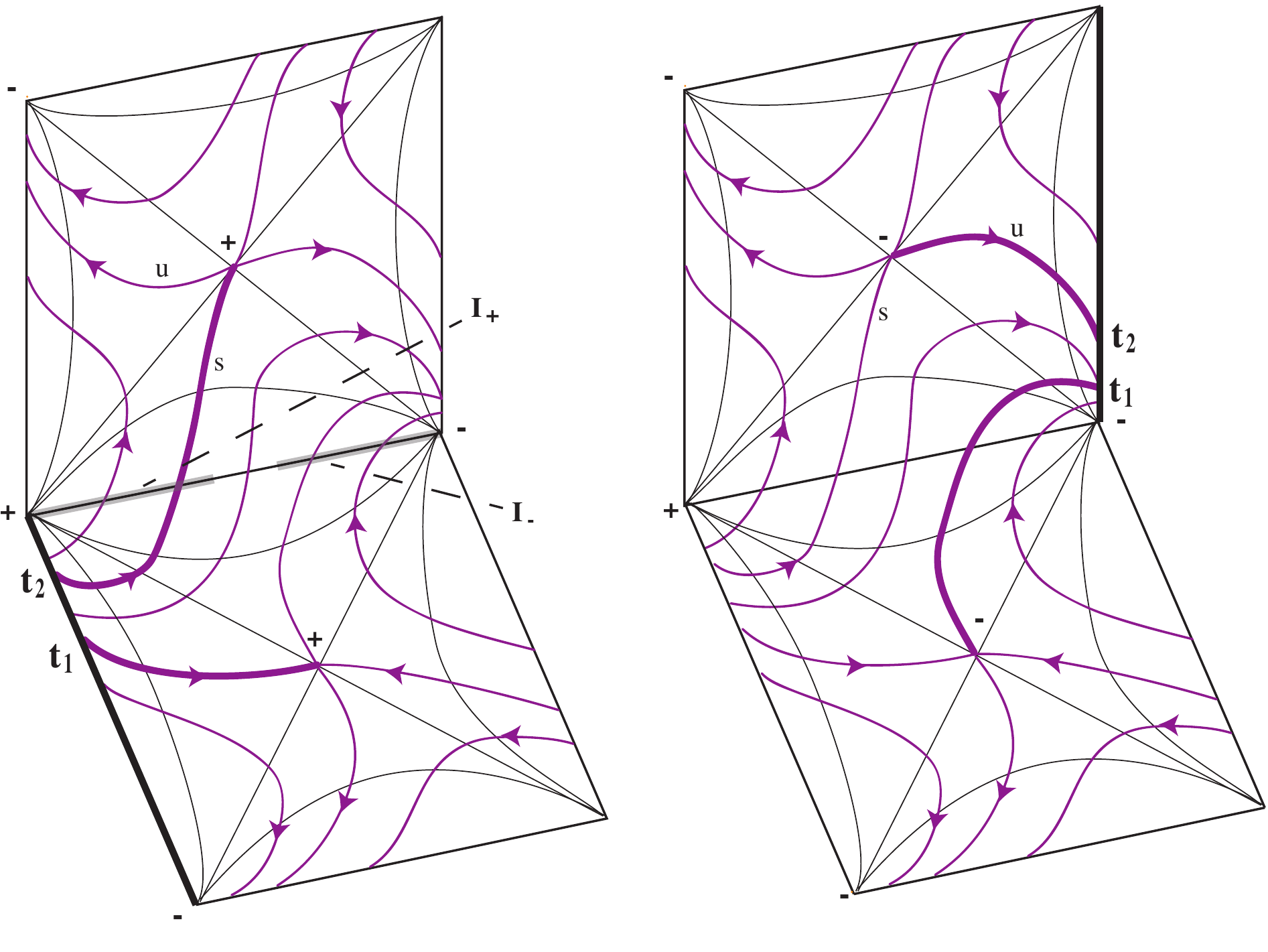}  
\caption{ The order of the positive (respectively, negative) hyperbolic points connected to a fixed elliptic point is encoded in the order that their  stable (unstable) separatrices first intersect $\pi^{-1}(0)$.  Each figure shows two tiles which share a regular leaf along the boundary.  On the left, the intervals $I_{+}$ and $I_{-}$ on a regular leaf are thickened.}\label{fig:orderofseps}
\end{center}
\end{figure}

The hyperbolic points of $\F_{\pit}$ inherit a cyclic order from the codomain of $\pit$, and we can partially recover this order from the vector field constructed above.  Namely, for the set of hyperbolic points of a fixed sign with separatrices to any fixed elliptic point, the order in which the seperatrices hit any regular leaf for the first time agrees with the restriction of the original cyclic order to this subset of hyperbolic points.  In Section~\ref{sec:ptofob}, we will require exactly this property from the gradient-like vector field for $\pit$.  See Figure~\ref{fig:orderofseps}.

\begin{definition}\label{def:pref} A gradient-like vector field $\nabla\pit$ for $\F_{\pit}$ is \emph{preferred} if  the following properties hold: 
\begin{itemize}
  \item[-]  for each regular time $t$ and for each component $I$ of $\pit^{-1}(t)$ with $\partial I=\{e_+, e_-\}$, there exist disjoint subintervals  $I_+$ and $I_-$  such that $I_+$ contains $e_+\cup \bigcup_h \big(w^s(h)\cap I\big)$ and  $I_{-}$ contains $e_-\cup \bigcup_h \big(w^u(h)\cap I\big)$.
 \item[-] for positive hyperbolic singularities $h_1, h_2$ with  critical values $0<\pit(h_1)<\pit(h_2) <1$,  the first intersection of $w^s(h_2)$  with $I_+$ is closer to $e_+$ than the first intersection of $w^s(h_1)$ with $I_+$;
 \item[-] for negative hyperbolic singularities $h_1, h_2$  with  $0<\pit(h_1)<\pit(h_2) <1$,  the first intersection of  $w^u(h_2)$ with  $I_-$ lies closer to $e_+$ than the first intersection of $w^u(h_1)$ with $I_-$. 
 \end{itemize}
\end{definition}

Definition~\ref{def:pref} lists  properties satisfied by the the vector field $\nabla\pit$ constructed above, proving the next lemma.

\begin{lemma} Any open book foliation admits a preferred gradient-like vector field. \qed
\end{lemma}

In Section~\ref{sec:pob} we will consider the question of when a preferred $\nabla \pit$ defined on $\partial M$ may be extended ``nicely'' to the interior of $M$.



\section{Foliated open books } \label{sec:fob}

In Section~\ref{ssec:efob} we introduce the central object of the paper,  embedded foliated open books.  We provide some examples of embedded foliated open books and define the compatibility between these and contact structures.  At the end of the section, we introduce two variations: abstract open books and Morse foliated open books. These offer increased flexibility and perspective, and the rather technical proofs relating these notions are deferred until Section~\ref{sec:fobequiv}.  Suppose throughout that  $M$ is a smooth oriented 3--manifold with boundary.

\subsection{Embedded foliated open books}\label{ssec:efob}
 The key object of this paper ---a foliated open book--- is motivated in part by trying to understand the result of cutting an open book $(B,\pi)$ along an embedded surface that admits an open book foliation with no circle leaves.

\begin{definition}\label{def:efob} An \emph{embedded foliated open book}  is a pair $(B_e,\pi_e)$, where 
$B_e$ is an oriented properly embedded 1-manifold in $M$   
and the function $\pi_e\colon M\setminus B_e\to S^1$ is a regular function such that the the following hold:
\begin{enumerate}
\item\label{efob:1}the restriction $\pit_e=\pi_e\vert_{\partial M}$ is an $S^1$-valued Morse function;  
\item\label{efob:0} the closure $S_t$ of each level set $\pi_e^{-1}(t)$ is a cornered surface with boundary $B\cup \pit^{-1}_e(t)$ and corners   $E=B\cap \pit_e^{-1}(t)$; 
\item\label{efob:2}  the restriction\footnote{We restrict to individual components of $\partial M$ so that cutting an open book along a non-separating surface will still produce a foliated open book.} of $\pit_e$  to each component of $\partial M$  has  a unique critical point for each critical value; 
\item\label{efob:4} the level sets of $\pit_e$ have no circle components.
\end{enumerate}
\end{definition}
\begin{remark}
From item (\ref{efob:4}) it follows that the critical points of $\pit_e$ can only be of index 1. \end{remark}

\begin{example}\label{prop:cut} Suppose $\Sigma$ is a surface  embedded  in a closed manifold $M$ with an open book $(B, \pi)$, and suppose that the open book foliation $\F_{\pit}$ has no circle leaves.  Then the restriction of $\pi$ to each component of the closure of $M\setminus\Sigma$ is an embedded foliated open book. 
\end{example}
\begin{example}
More generally, one may consider a surface surface $\Sigma$ with boundary $\partial {\Sigma}=T$ transverse to the pages of $(B,\pi)$. In this case the open book foliation $\F_{\pit}$ on $\Sigma$  points transversally out of or into $T$, depending on whether component of $T$ in question is positively or negatively transverse to the pages. After cutting along such a $\Sigma$, the  boundary of   $M\setminus\Sigma$ is the double of $\Sigma$ and the foliation is the union along $T$ of $\F_{\pit}$ with its negatively oriented copy.
\end{example}

After cutting along  an embedded $\Sigma$, the new components of the boundary retain open book foliations in the sense of Definition~\ref{def:obfol}.  In this case, the signs of the hyperbolic points are induced by their signs as singularities of the embedded $\Sigma$, but the same reasoning assigns signs to hyperbolic singularities on the boundary of any embedded foliated open book.  

Choose a gradient-like vector field $\nabla\pi_e$ for $\pi_e$; here, only, we consider gradient-like vector fields which can be transverse to the boundary.
Let $H_+$ denote the set of hyperbolic singularities of $\pit_e$ where $\nabla \pi_e$ points out of $M$, let $H_-$ be the set of hyperbolic singularities  of $\pit_e$ where $\nabla \pi_e$ points into $M$, and let $E_+\cup -E_-=B_e\cap\partial M$.  We assume  these signs in the next statement; the proof follows immediately from the definitions just introduced.

\begin{proposition} Suppose that $B_e$ is a properly embedded 1--manifold in $M$ and $\pi_e\colon M\setminus B_e\to S^1$ is a regular function. Then $(B_e, \pi_e)$ is an embedded foliated open book if and only if $(\F_{\pit_e},\pit_e, H=H_+\cup H_- )$ is an open book foliation on $\partial M$ with no circle leaves. \qed
\end{proposition}

Interestingly, even Item~(\ref{efob:0}) in Definition~\ref{def:efob} is forced on $\pi_e$ by the conditions on $\F_{\pit_e}$ near the elliptic points.

In some settings, we will not want to distinguish circle-valued functions that are just reparametrizations of each other, so we introduce an equivalence relation that makes sense in the context of surfaces with open book foliations as well as open books of various sorts.

\begin{definition}\label{def:reparam} Two (embedded foliated) open books $(B, \pi)$ and $(B', \pi')$ are \emph{reparamterizations} if there exist neighborhoods $N(B)$ and $N(B')$ and a diffeomorphism $p: S^1\rightarrow S^1$ satisfying the following:
\begin{enumerate}
\item for each component $B_i$ of $B$ or $B'$, the designated neighborhood intersects each page in a subsurface homeomorphic to $B_i\times I$; and 
\item on the complement of the designated neighborhoods,  $\pi'=p\circ \pi$.  
\end{enumerate}
Two open book foliations are reparameterizations if there exist neighborhoods $N(E)$ and $N(E')$ and a diffeomorphism $p: S^1\rightarrow S^1$ such that on the complement of the neighborhoods of the elliptic points,  $\pi'=p\circ \pi$.
\end{definition}

This notion of equivalence is a reasonable one to consider for several reasons.  Readers familiar with abstract open books will recognise it as an unavoidable indeterminacy in a manifold constructed from the abstract data $(S, h)$, although it is often not stated explicitly. It will be useful to have the flexibility to reparameterize open books when we state gluing theorems in Section~\ref{ssec:glue}, and in some cases we are concerned only with the underlying signed oriented foliation on a surface, rather than the full data of an open book foliation.  

From now on we will refer to embedded foliated open books by the triple $(B_e,\pi_e,\F_{\pit_e})$ to emphasize the open book foliation  on the boundary.   Unless otherwise specified, we will consider $\pi_e$ and $\F_{\pi_e}$ up to reparameterization.

Given a manifold $M$ whose boundary has an open book foliation $(\F_{\pit}, \pit, H=H_-\cup H_+)$,
 we say that an embedded  foliated open book $(B_e,\pi_e,\F_{\pit_e})$ for $M$ is \emph{compatible} with the pair $(M,\F_{\pit})$ if $\F_{\pit_e}=\F_{\pit}$.    In general we call an oriented $3$-manifold $M$ together with an open book foliation on $\partial M$ that is defined up to reparameterization   a \emph{3-manifold with foliated boundary}.  To simplify the notation, we will denote this by $(M,\F_{\pit})$, rather than the more precise $(M,\F_{[\pi]})$.

Two 3-manifolds with foliated boundary $(M,\F_{\pit})$ and $(M',\F_{\pit'})$ are \emph{diffeomorphic} if there is a diffeomorphism $\psi\colon M\to M'$ that  restricts on the boundary to a diffeomorphism of the open book foliations $\F_{\pit}$ to $\F_{\pit'}$.  If $M=M'$ and $\psi$ is an isotopy, we say that $(M,\F_{\pit})$ and $(M,\F_{\pit'})$ are \emph{isotopic} 3-manifolds with foliated boundaries. 

Two embedded foliated open books $(B_e,\pi_e,\F_{\pit_e})$ for $(M,\F_{\pit_e})$ and $(B'_e,\pi'_e,\F'_{\pit_e})$ for $(M',\F_{\pit_e'})$ are \emph{diffeomorphic}  if there is a diffeomorphism $\psi:(M,\F_{\pit_e})\rightarrow (M',\F_{\pit_e'})$ that takes $B_e$ to $B'_e$ and $\pi'_e\circ\psi=\pi_e$. If $M=M'$ and $\psi$ is isotopic to the identity, then we say that the embedded foliated open books $(B_e,\pi_e,\F_{\pit_e})$ and $(B'_e,\pi'_e,\F'_{\pit_e})$ are \emph{isotopic} to each other.

Examples of foliated open books are given in Section \ref{ssec:cutex} and in Corollary \ref{cor:lm}. The interested reader is advised to review these now for better intuition. 

\subsection{Compatible contact structures} 
 Foliated open books, like their partial and closed cousins, are intended as a tool to study contact manifolds.

\begin{definition}\label{def:supportxi} Given  a contact manifold $(M, \xi)$ with boundary, let $\F_\xi$ denote  the characteristic foliation on $\partial M$. The embedded foliated open book ($B,\pi_e, \mathcal{F}_{\pit_e}$)  \textit{supports} the contact structure $\xi$ if  there is a contact 1-form $\alpha$ for $\xi$ such that
\begin{enumerate} 
\item\label{cond1}  $\d\alpha$ is a positive area form on the interior of each \emph{page} $\mathring{S_t}=\pi^{-1}(t)$;
\item\label{cond2}  $\alpha>0$ on $TB$; 
\item\label{cond3} $\mathcal{F}_\pit$ and $\mathcal{F}_\xi$ are strongly topologically conjugate  on $\partial M$.
\end{enumerate} 
\end{definition}

\begin{remark} By Lemma \ref{lem:beta}, $\Gamma$  (geometrically) divides the characteristic foliation. Thus $\partial M$ is automatically convex.
\end{remark}
In Section~\ref{sec:uniquexi} we will prove the following:
\begin{theorem}\label{thm:support}
Any embedded foliated open book supports a unique contact structure.
\end{theorem}

We prove the converse  in Section~\ref{sec:existfob}:
\begin{theorem}\label{thm:fobexistence} Let $(M,\xi,\F_\xi)$ be a contact manifold and assume that $\F_\xi$ is strongly topologically conjugate to an open book foliation $\F_{\pit}$ with no circle leaves. Then there is a foliated open book $(B,\pi,\F_{\pit})$ supporting $\xi$. In particular,  $\partial B=E_+\cup-E_-$ and $\pi\colon M\setminus B\to S^1$ is an extension of $\pit$.
\end{theorem}

The following example shows that the circle free condition is essential for the claim that a supported contact structure always exists.
\begin{example}
Suppose that $S$ is a component of a page of a foliated open book which is topologically a disc, so that $\partial S$ is a circle component of a leaf of the open book foliation on the boundary. Since $\F_\xi$ and $\F_{\pit}$ on $\partial M$ are strongly topologically conjugate by definition, we may isotope $S$ so that the characteristic foliation  has a circle $C$ bounding a disc. By construction, this circle is  Legendrian.  
 
Since $C$  comes in an interval of Legendrian circles in $\F_\pit$, its linking number with a Legendrian push-off ---and hence, its Thurston-Bennequin number--- equals $0$.  On the other hand, we will show in Proposition~\ref{prop:unique} that the restriction of any supported contact structure to a handle body associated to a sequence of consecutive  pages is tight, which contradicts the existence of the overtwisted disc $S$.
 
\end{example}

\subsection{Abstract foliated open books}\label{ssec:afob} 
The discretised version of an embedded foliated open book is an abstract foliated open book. These are combinatorial objects, and they   specify manifolds with foliated boundary only up to diffeomorphism.  
Here we present an independent definition for abstract foliated open books, one which will be unsurprising for readers familiar with other types of open books.  We defer an explanation of  the  equivalence between abstract and embedded foliated open books to Section~\ref{sec:fobequiv}, after we have introduced an intermediate object, called Morse foliated open books, in Section~\ref{sec:morsefob}.

\begin{definition}\label{def:afob}An \emph{abstract foliated open book} is a tuple $(\{S_i\}_{i=1}^{2k}, h)$ where $S_i$ is a surface with 
boundary $\partial S_i=B\cup \alpha_i$ \footnote{By a slight abuse of notation we denote the ''constant'' part of the boundary of $S_i$ by $B$ for all $i$.} and corners at  $E=B\cap \alpha_i$ such that 
\begin{enumerate}
 \item for all $i$, $\alpha_i$ is a union of intervals;
\item the surface $S_{i}$ is obtained from $S_{i-1}$ by either 
\begin{itemize}
\item[-\textbf{(add):}]attaching a 1-handle along two points $\{p_{i-1}, q_{i-1}\}\in \alpha_{i-1}$; or 
\item[-\textbf{(cut):} ] cutting $S_{i-1}$ along a properly embedded arc $\gamma_{i-1}$ with endpoints in $\alpha_{i-1}$ and then smoothing.
(See Figure \ref{fig:saddles}). 
\end{itemize}
\end{enumerate}

\begin{figure}[h]
\begin{center}
\includegraphics[scale=.8]{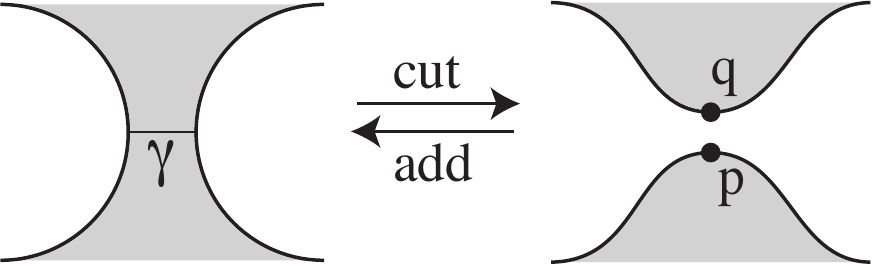}
\caption{The cutting and adding operations on successive pages are inverses of each other.}\label{fig:saddles}
\end{center}
\end{figure}

Furthermore,  $h: S_{2k}\rightarrow S_0$ is a diffeomorphism between cornered surfaces that
preserves $B$ pointwise. 
\end{definition}
Note that the operations \textbf{(add)} and \textbf{(cut)} are opposites of each other; if $S\xrightarrow{\text{\textbf{add}}} S'$ along $p$ and $q$, then $S'\xrightarrow{\text{\textbf{cut}}} S$ along the cocore of the new 1-handle, and vice versa. When we would like to specifically identify the attaching sphere or arc, we will record it under the arrow, as below.

\begin{definition}Two abstract foliated open books $(\{S_i\}_{i=0}^{2k},h)$ and $(\{S_i'\}_{i=0}^{2k'},h')$ are 
 \emph{diffeomorphic} to each other if $k=k'$ and there is a sequence of cornered diffeomorphisms $\{\psi_i\colon S_i\to S_i'\}_{i=1}^k$, all of which agree on $B$;  satisfy $\psi_i(B)=B'$;  and are compatible with the handle attachments and cutting: 
 \begin{itemize}
 \item[-] if  $S_{i-1}\xrightarrow[\gamma_{i-1}]{\textbf{cut}}S_i$, then $S_{i-1}'\xrightarrow[\gamma'_{i-1}]{\text{\textbf{cut}}} S_i'$, where
  $\gamma_{i-1}'=\psi_{i-1}(\gamma_{i-1})$ and $\psi_i$ is the restriction of $\psi_{i-1}$ to $S_i\subset S_{i-1}$;
  \item[-]  if $S_{i-1}\xrightarrow[p_{i-1},q_{i-1}]{\text{\textbf{add}}} S_i$, then $S_{i-1}'\xrightarrow[\psi_{i-1}(p_{i-1}),\psi_{i-1}(q_{i-1})]{\text{\textbf{add}}} S_i'$, where $\psi_{i}$ restricts as $\psi_{i-1}$ on $S_{i-1}$ and maps the attached  handles to each other; and
  \item[-]  
$h'= \psi_0\circ h\circ\psi_{2k}^{-1}$. 
 \end{itemize}

\end{definition}

Next, we define an operation that allows us to freely choose which page is indexed as $S_0$.  This is more than merely notation, as our definition requires the monodromy to be a homeomorphism of this page. 
\begin{definition} 
The \emph{shift} of an abstract open book $(\{S_i\}_{i=0}^{2k},h)$ is  \[(\{S_i[1]\}_{i=0}^{2k},h[1])=(\{S_1,S_2,\dots,S_{2k},S_1'\},h'),\] where
\begin{itemize}
\item[-] if $S_{0}\xrightarrow[\gamma_0]{\text{\textbf{cut}}} S_1$, then $S_1'$ is defined by the relation $S_{2k}\xrightarrow[{h^{-1}(\gamma_0)}]{\text{\textbf{cut}}} S_1'$ and $h'$ is the restriction of $h$  to $S_1'$;

\item[-] if $S_{0}\xrightarrow[{p_0, q_0}]{\text{\textbf{add}}} S_1$, then $S_1'$ is defined by the relation $S_{2k}\xrightarrow[{ h^{-1}(p_0), h^{-1}(q_0) }]{\text{\textbf{add}}} S_1'$ and $h'$ is $h$ extended by the identity on the added 1-handle. 
\end{itemize}
An $r$-fold iteration of the shift operation is called an \textit{$r$-shift}, and denoted by $(\{S_i[r]\}_{i=0}^{2k},h[r])$. One can analogously define $r$-shifts for $r<0$.
\end{definition}
As we will see, the above two moves preserve the diffeomorphism type of the defined 3--manifolds with foliated boundary. Two abstract foliated open books are \emph{conjugates} of each other if they are related by a finite sequence of diffeomorphisms and shifts.

\subsection{Morse foliated open books}\label{sec:morsefob} To conclude this section, we introduce a final type of foliated open book. 
Morse foliated open books form a bridge between embedded foliated open books and abstract foliated open books. Their utility arises from the fact that --- in contrast to embedded foliated open books--- the $S^1$-valued function defining the pages of $M$ has the same critical points as its restriction to $\partial M$.  Thus the boundary critical points have stable and unstable manifolds embedded as submanifolds of $M$.  A Morse foliated open book may be obtained from an embedded foliated open book by slightly modifying $\pi_e$ near the boundary, a procedure described in detail in Section~\ref{sec:mem}.  The pages of embedded foliated open books and Morse foliated open books are essentially the same, and we will not distinguish between them once we have clarified their relationship.

\begin{definition}\label{def:mfob} A \emph{Morse foliated open book} for $M$ is a pair $(B_m,\pi_m)$, where 
$B_m$ is an oriented properly embedded 1-manifold in $M$  
and the function $\pi_m\colon M\setminus B_m\to S^1$ is an $S^1$-valued Morse function such that the following hold: 
\begin{enumerate}
 \item\label{mfob:5} the level sets of $\pit_m=\pi_m|_{\partial M}$ have no circle components;

\item\label{mfob:0} the closure $S_t$ of each level set $\pi_m^{-1}(t)$ is a cornered surface with boundary $B\cup \pit_m^{-1}(t)$ and corners   $E=B\cap \pit_m^{-1}(t)$; 
\item\label{mfob:1} all the critical points of $\pi_m$ are on $\partial M$; 
\item\label{mfob:2} $\pit_m$ is Morse function with the same critical points as $\pi_m$; and  
\item\label{mfob:4} the restriction of $\pit_m$  to each component of $\partial M$ has a unique critical point for each critical value.
\end{enumerate} 
\end{definition}
\begin{remark}It follows from condition (\ref{mfob:5}) that the critical points of $\pit_m$ have index 1. Condition (\ref{mfob:2}) implies that the critical points of $\pi_m$ have index 1 or 2 and that $\pi_m$ has a gradient-like vector field $\nabla \pi_m$ that is tangent to $\partial M$. All gradient like vector fields for $\pi_m$ will automatically be assumed to have this property.
\end{remark}
 As in the case of an embedded foliated open book, the boundary of a Morse foliated open book naturally inherits an open book foliation $\F_{\pit_m}$.  In this case, $\partial B_m=E_+\cup-E_-$, and $H_+$ is the set of index 2  critical points of $\pi_m$, while $H_-$ is the set of index 1  critical points of $\pi_m$.

Just as for embedded open books we can define \emph{diffeomorphism} and \emph{isotopy} of Morse foliated open books. The definitions for \emph{supported} contact structure and \emph{stabilisation}  depend only on the binding and the pages of an embedded foliated open book, rather than the specific function $\pi_e$, so these definitions can be  extended verbatim to Morse foliated open books.


\section{Local Models } \label{sec:models}

This is a  technical section that may be skipped at first reading. As a start, we provide a neighborhood theorem for open book foliations on surfaces in Section~\ref{sec:nhood}, and we construct an explicit model for this neighborhood in Section~\ref{sec:egcircle}. Taken as a whole, this section lays the groundwork for moving between embedded, Morse, and abstract foliated open books.   Such  freedom will be extensively employed in the rest of the paper, and the reader may wish to survey the results  in Section~\ref{sec:fobequiv} which describe the precise notions of equivalence between the various types of foliated open books.

\subsection{Neighborhood Theorem}\label{sec:nhood}
The key theorem extends to open book foliations of the usual neighborhood theorems for characteristic foliations. 
\begin{proposition}[Neighborhood Theorem]\label{prop:nbrh} 
Let $(B,\pi)$ and $(B',\pi')$ be open books for the 3--manifolds $M$ and $M'$ respectively. Let $\Sigma\hookrightarrow M$ and $\Sigma'\hookrightarrow M'$ be embedded surfaces with induced open book foliations $(\F_{\pit}, \pit,H=H_+\cup H_-)$ and $(\F_{\pit}',\pit', H'=H'_+\cup H'_-)$. Suppose that there is a diffeomorphism $\psi\colon \Sigma\to \Sigma'$ that takes the two open book foliations to each other; i.e., $\psi(E)=E'$, $\pit'\circ\psi=\widetilde{\pi}$, and $\psi(H_\pm)=H_\pm'$. Then there are neighborhoods $N=N(\Sigma)$ and $N'=N(\Sigma')$ and an extension of $\psi$ to $\Psi\colon N\to N'$ so that 
\begin{enumerate}
\item $\Psi(B)=B'$;
 \item $\pi'\circ\Psi=\pi$;
 \item $\Psi(B\cap N)=B'\cap N'$ and $\Psi$ maps the trivialisation of $N(B)\cap N\cong (B\cap N)\times D^2$ in which $\pi=\vartheta$ to the trivialisation of $N(B')\cap N'\cong (B'\cap N')\times D^2$ in which $\pi'=\vartheta$. 
 \end{enumerate}
Moreover, one can choose supported contact structures $(M,\xi)$ and $(M',\xi')$ for $(B,\pi)$ and $(B',\pi')$,  respectively, so that $\Psi$ is a contactomorphism between $(N,\xi\vert_{N})$ and $(N',\xi'\vert_{N'})$. 
\end{proposition}

\begin{proof}
Let us first describe local models around the points of $\Sigma\subset M$. For each point $p\in {\Sigma}$ we will give a neighborhood $D_p^2\times I$ of $p$ in $M$ with coordinates adapted to $(B,\pi)$. These coordinates are chosen by a repeated use of the Implicit Function Theorem, and both the neighborhood $D_p^2$ and the interval $I$ might shrink as we make our additional choices. 

In the neighborhood of a \emph{regular point} $p$ of the foliation ${\F}_{\pit}$ we can choose local coordinates $(u,v,z)$ on $M$ so that $v=\pi-\pit(p)$, ${\Sigma}=\{z=0\}$ with orientation $(\partial u,\partial v)$ and coorientation $\partial z$, and the foliation ${\F}_{\pit}$ is directed by $\partial u$. We call such a coordinate system in neighborhood that only contains regular points of ${\F}_{\pit}$ an \emph{adapted coordinate system}. The set of adapted coordinates in a given neighborhood $U$ is convex.

Around an \emph{elliptic point} $e\in E_\pm$ of ${\F}_{\pit}$, we know that $B\pitchfork {\Sigma}$ at $e\in E_\pm$, and $N(B)\cong B(\varphi)\times D^2(r,\vartheta)$. In these coordinates, $e=(\varphi_0,0,0)$, and in a sufficiently small neighborhood we can write $\Sigma\cap N(B)$ as the graph of a function $f\colon \{\varphi_0\}\times D^2\to \R$ that fixes $\{\varphi_0\}\times\{0\}$. Since $\Sigma$ is transverse to $B$ we have $df\neq d\varphi$ at the origin of $\{\varphi_0\}\times D^2$, so  the Implicit Function Theorem allows us to choose coordinates $(z=f-\varphi$, $r,\vartheta)$. In these coordinates $\Sigma=\{z=0\}$, $B$ is oriented by $\partial z$, $\Sigma$ is cooriented by $\pm\partial z$, and $\pi=\vartheta$.

In a neighborhood of a \emph{hyperbolic point} $h\in H_\pm$,  the function $\pit=\pi\vert_{{\Sigma}}$ is Morse  with  an index 1 critical point $h$. By the Morse Lemma and since the differential of $\pi$ in the $I$-direction is nonzero, we can apply the Implicit Function Theorem and choose coordinates $(x,y,z)$ so that ${\Sigma}=\{z=0\}$ and $\pi-\pi(h)=z-y^2+x^2$. The surface ${\Sigma}$ is cooriented by $\pm\partial z$ and ${\F}_{\pit}$ is directed by $y\partial x + x\partial y$.

As a next step we consider the intersections of these coordinate systems.  First,  choosing sufficiently small neighborhoods  ensures that the neighborhoods of elliptic and hyperbolic points described above are all disjoint.

If $q$ is a regular point  in the neighborhood of another regular point $p$ with a coordinate system $(u,v,z)$ as above, then we can write $q=(u_0,v_0,0)$ and construct a new adapted coordinate system: \[(u',v',z')=(u-u_0,v-v_0,z).\] 
 
Let $p=(0, r_0,\vartheta_0)\neq (0,0,0)$ be a regular point on $\Sigma$ in the neighborhood of $e\in E_\pm$. Then $\pi(p)=\vartheta_0$ and the coordinates: \[(u',v',z')=(\pm (r-r_0),\vartheta-\vartheta_0, \pm z)\] give an adapted coordinate system around $p$.

Similarly, let $p=(x_0,y_0,0)\neq (0,0,0)$ be a point on ${\Sigma}$ in the neighborhood of $h\in H_\pm$.

Then
\[(u',v',z')=\big(xy-x_0y_0,(z-y^2+x^2)-(-y_0^2+x_0^2), \pm z\big)\] 

gives an adapted coordinate system in a neighborhood of $p$.

Finally, using these local models for $\F_{\pit}$ and $\F_{\pit}'$ we can define local maps that take the corresponding coordinate systems to each other. Then we can use a partition of unity to construct a global map $\Psi\colon N \to N'$ from the local ones. Near the elliptic points, the map brings the coordinate systems for the elliptic points to each other, so $\psi$ satisfies conclusion (3). Conclusion (1) is automatically satisfied by the construction. As for conclusion (2), we need to check that $\pi'\circ \Psi=\pi$. This is certainly true for the local maps, and the construction of $\Psi$ used their convex combinations. In the above change of coordinate systems, the value of $\pi$ was implicit in the system, so we can assume that we only need to take the convex combination of adapted coordinate systems in the neighborhood of regular points. These form a convex set, thus the value of $\pi'\circ \Psi$ is unchanged while taking their convex combination.

As for the second part, concerning compatible contact structures, let $(S,h)$ be an abstract open book corresponding to $(B,\pi)$. From now on we will assume that $M$ is identified with $M(S,h)$. We briefly recall the Thurston--Winkelnkemper construction \cite{TW} of a contact structure $\widetilde{\xi}=\ker\widetilde{\alpha}$ on $M$ supported by $(B,\pi)$. The construction depends on 
\begin{itemize}
\item[-] a 1-parameter family of 1-forms $(\beta_t)_{t\in [0,1]}$ on the page $S$, such that $d\beta_t$ is an area form on $S$ with total area $2\pi$ and $h_*\beta_0=\beta_1$; and 
\item[-] a sufficiently large constant $\widetilde{C}$. 
\end{itemize}
Then the contact form $\widetilde{\alpha}$ is defined away from $N(B)$ by $\beta_t+\widetilde{C}dt$. In a smaller neighborhood of the binding identified as $B(\varphi) \times D^2(r,\vartheta)$, the contact form $\widetilde{\alpha}$ is given by $2d\varphi+r^2d\vartheta$. In the complement of the two neighborhoods the construction gives an explicit mutual extension of the two 1-forms. 

Let $(S',h')$ be the abstract open book corresponding to $(B',\pi')$, and as before, assume that $M'$ is identified with $M(S',h')$. Take the 1-parameter family of 1-forms  $\widetilde{\beta}'_t=(\Psi\vert_{(S\times\{t\})\cap U})_*\beta_t$ on $(S'\times \{t\})\cap N'$ and choose an area form on $S'$ so that $d\widetilde{\beta}'_t$ has total area less than $2\pi-\varepsilon$ on any closed subset of $(S'\times \{t\})\cap N'$.
Since $\Psi$ respects the monodromies $h$ and $h'$, we have that 
$\widetilde{\beta}_1=h'_*\widetilde{\beta}'_0$ wherever both are defined. Thus we can extend $\widetilde{\beta}'_t$ to an area form $\beta'_t$ on $S'_t$ satisfying the conditions of the Thurston-Winkelnkemper construction, and now a compatible contact form may be defined by a choice of a sufficiently large constant $\widetilde{C}'$. 

Let $C=\max\{\widetilde{C},\widetilde{C}'\}$ and define the contact forms $\alpha$ on $M$ and $\alpha'$ on $M'$ by the Thurston--Winkelnkemper construction using the parameters $\beta_t$ and $C$, and $\beta_t'$ and $C$, respectively. For the corresponding contact structures $\xi=\ker\alpha$ and $\xi'=\ker\alpha'$ we get that $\Psi\vert_N$ is a contactomorphism between $(N,\xi\vert_N)$ and $(N',\xi'\vert_{N'})$, as needed. 
\end{proof}
In order to give an explicit model for the open book in $N(\Sigma)$,  the next subsection  constructs a simple open book for a $\Sigma$-bundle over a circle with a prescribed open book foliation on $\Sigma\times \{0\}$.

\subsection{
Prescribing open book foliations on surfaces}\label{sec:egcircle}
In this section we will show that every circle-free open book foliation $(\F_{\pit}, \pit,H=H_+\cup H_-)$ on a  surface $\Sigma$ is in fact induced by an embedding of $\Sigma$ into an open book.

As $(\F_{\pit}, \pit, H=H_+\cup H_-)$ is  circle-free, we can construct a dividing curve $\Gamma$ and   see 
that $\pit\vert_\Gamma\colon \Gamma\to S^1$ is a covering of degree $n$, where $n=|E_+|=|E_-|$. Let $\{\Gamma_i\}_{i=1}^k$ be the set of connected components of $\Gamma$ and let $n_i$ be the degree of $\pit\vert_{\Gamma_i}$. Note that $n=\sum_{i=1}^k n_i$. Then, just as in Section~\ref{ssec:charfol} we can choose local coordinates $(u,v)$ on $A_i=N(\Gamma_i)$ so that $\pit=n_iv$, $\Gamma=\{u=0\}$, and $\partial u$ directs the level sets $\pit^{-1}(t)$. Define $A=\cup A_i=N(\Gamma)$ and $R_\pm'=R_\pm\setminus  A$. 
 Next we construct a characteristic foliation on $\Sigma$ that is also divided by $\Gamma$. 
 
\begin{lemma}[Approximating $\F_{\pit}$ by $\F_{\xi}$]\label{lem:approx}
Let $\Gamma(\F_{\pit})$ be the dividing curve and $R_\pm(\F_{\pit})$ the corresponding positive and negative regions of the open book foliation $(\F_{\pit}, \pit,H=H_+\cup H_-)$. 
Then there is a 1-form $\beta$ that agrees with $d\pit$ away from a small neighborhood of the singular points and $\pm d\beta>0$ on $E_\pm\cup H_\pm$.  The foliation $\F_\xi$ defined by $\beta$ is a characteristic foliation, has the same singular points and is strongly topologically conjugate to $\F_{\pit}$.

Moreover, we can choose a representative $\beta'=g\beta$ so that $\Gamma=\{d\beta'=0\}$ and $R_\pm=\{\pm d\beta'>0\}$, where $g$ is a positive function and $g<1$ on $A$.
\end{lemma} 

\begin{proof}
For the first statement we start with the 1-form $\gamma$ as in Section \ref{sec:obfol}~ and modify it in the neighborhood of the hyperbolic points of $\F_{\pit}$.
Recall that an elliptic point  $e \in E_\pm$ has a neighborhood $D^2_{\varepsilon/4}(r,\vartheta)$ where $\gamma=\pm r^2 d \vartheta$.  Thus at $e$  we have $\pm d\beta=2r dr\wedge d\vartheta$, which is a positive multiple of a volume form.

In a neighborhood of a  hyperbolic point $h\in H_\pm$, we can chose local coordinates  so that $\pit=x^2-y^2$. Thus $\gamma=d\pit=2xdx-2ydy$ and $d^2\pit=0$.  Let $\beta=d\pit\pm \Psi^\varepsilon_{\varepsilon/2} xdy$. Then, as $d\Psi^\varepsilon_{\varepsilon/2}=0$ at $r=0$, at $h$ we have that $\pm d\beta= dx\wedge dy$, a positive multiple of the volume form. 

Extend the 1-form on $\Sigma\setminus N(H)$ as $d\pit$ to  obtain a 1-form $\beta$ that satisfies the conditions of the first claim. 

Without the assumption $g<1$ on $A$, the second statement is standard in contact geometry; the stronger statement is Remark~\ref{rmk:g} after  the proof of Lemma \ref{lem:beta}.
\end{proof}

We use the function $f$ from Proposition \ref{prop:Gi} and the 1-form $\beta'=g\beta$ from Lemma \ref{lem:approx} to construct both a contact structure and an open book for a $\Sigma$-bundle over  $S^1$. Consider the product $\Sigma\times \R$ with the contact structure defined as the kernel of the 1-form $\alpha'=\beta'+f dz$. Let $B'=E\times \R$, where the orientation of $\{e\}\times \R$ is given by $\pm\partial z$ for $e\in E_\pm$. Consider the function \[\pi'=\pit+fz\colon \Sigma\times \R\to S^1,\]  
where adding $f(p)z(p)\in \mathbb{R}$ to $\pit(p)\in S^1$ (for $p\in \Sigma$) indicates translation by the image of $f(p)z(p)$ in the quotient $S^1=\mathbb{R}\slash \mathbb{Z}$. 

We would like to glue $\Sigma\times \{0\}$ to $\Sigma\times\{l\}$ for some $l\in\mathbb{N}$ to get an open book and a contact structure for a $\Sigma$-bundle over  $S^1$ . As $f=\pm 1$ on $R_\pm'$, we have $\pi'(x,0)=\pi'(x,l)$ on $R_\pm'$. 
On $A_i=N(\Gamma_i)\cong\Gamma_i\times [-\varepsilon,\varepsilon]$, with coordinates $(u,v)$ chosen as in Lemma \ref{lem:approx}, we have $\pi'(u, v,l)-\pi'(u,v,0)=lf(u,v)$, where $f(u,v)$ is independent of $v$  and monotonically decreases from $1$ to $-1$ in $u$.  This means that each level set of $\pi'\vert_{\Sigma\times\{l\}}$ restricted to $A_i$ consists of $n_i$ parallel curves connecting the points $(-\varepsilon,v)$ with $(\varepsilon,v+2l/n_i)$.  When $2l$ is divisible by $n_i$, the endpoints of each component therefore have the same $v$ coordinate, as $v$ parmeterizes $S^1\cong \R/\mathbb{Z}$. And in fact, the curve intersects an arbitrary  $\{v=c\}$ segment positively in $\frac{2l}{n_i}\in \mathbb{Z}$ times. 

Set  $l=\lcm\{n_i\}_{i=1}^k$ and define a diffeomorphism
 \[\psi=\prod_{i=1}^k D_{\Gamma_i}^{\frac{2l}{n_i}},\] where $D_{\Gamma_i}$ is the right-handed Dehn twist along  $\Gamma_i$.  Let  $M_\psi(\Sigma)=\Sigma\times [0,l]/(\psi(x),0)\sim (x,l)$ be the mapping torus of $\psi$.  Then $\pi'(x,0)=\pi'(\psi(x),l)$ on $S^1$ and $\pi'$ descends as $\pi$ to  $M_\psi(\Sigma)\setminus  B$, where $B=B'\slash\sim$.

Notice, too, that since $f$ depends only on $u$ and the Dehn twists were along the $v$-direction in the $A_i$, the 1-form $\alpha'$ descends  to $M_\psi$ as a 1--form $\alpha$, giving a contact structure $\xi$ on $M_\psi$.

\begin{proposition}[Open book for $(M_\psi,\xi)$]\label{prop:obcyl} Via the construction above, any $\F_{\pit}$ determines an open book decomposition $(B,\pi)$ for $M_\psi$. Moreover, $(B,\pi)$ supports the contact structure $\xi=\ker\alpha$.  
\end{proposition} 
Before proving the above proposition, we reformulate what it means for an open book decomposition to support a contact form.

\begin{lemma}\label{lem:ob} Let $B$ be a (positively) transverse  knot in $(M,\xi=\ker\alpha)$ and suppose that the fibers of the map $\pi\colon M\setminus  B\to S^1$  are Seifert surfaces for $B$. Then the pair $(B,\pi)$ is an open book  supporting the contact form $\alpha$ if and only is $d\pi\wedge d\alpha>0$. \qed
\end{lemma}

The proof is straightforward and is left to the reader.

\begin{proof}[Proof of Proposition \ref{prop:obcyl}.]

In a neighborhood of $B$ the function $\pi$ restricts as $\pm(\vartheta+z)$,
 so we see that $\pi^{-1}(t)$ indeed enters $N(B)$ as a Seifert surface. If we evaluate $\alpha$ on $TB=\pm\partial z$ we get $\alpha(\pm\partial z)=1>0$, as needed.

For the second step,  by Lemma \ref{lem:ob}~ it is sufficient to show that $d\pi\wedge d\alpha >0$ on $M_\psi\setminus  B$. Since $\alpha$ is a contact form,  $\alpha\wedge d\alpha>0$, yielding the following:
\[
(\beta'+fdz)\wedge (d\beta'+df\wedge dz)=dz\wedge(\beta'\wedge df+fd\beta')>0.
\]
 This in turn is equivalent to the condition that 
 $
\omega_0=\beta'\wedge df+fd\beta'$
 is an area form on $\Sigma$. 

We compute $d\pi \wedge d\alpha$:
 
\[
(d\pit+zdf+fdz)\wedge(d\beta'+df\wedge dz)=dz\wedge(d\pit\wedge df+fd\beta')
\]
Thus we require that  
 $\omega'=d\pit\wedge df+fd\beta'$ is an area form on $\Sigma$.  To check this, we compare the given form to the area form $\omega_0$ given by the contact condition and substitute $\beta'=g\beta$:
\[\omega'=\omega_0+(d\pit-g\beta)\wedge df.\]
In the following we will prove that $(d\pit-g\beta)\wedge df\ge 0$, which implies that  $\omega'>0$.

On $R_\pm'$ the differential $df$ is 0, thus $\omega'=\omega_0$, so it is indeed an area form.

On $A$,  $\beta=d\pit$, so we have   
\[(d\pit-g\beta)\wedge df=(1-g)d\pit \wedge df.\]
Recall that on $A$, $f$ depends only on $u$ and is decreasing, while  
$\pit$ depends only on $v$ and is increasing.    
Thus $d\pit \wedge df$ is an area form, and as $g<1$ on $A$, so is $(1-g)d\pit \wedge df$, as desired. 
\end{proof}

\begin{figure}[h]
\begin{center}
\includegraphics[scale=.6]{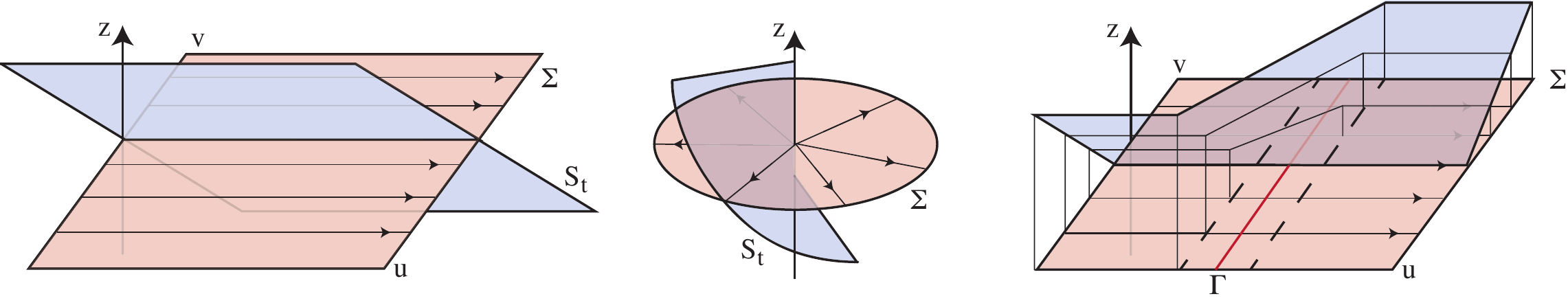}
\caption{   Local models for the open book of Proposition~\ref{lem:ob}.  The blue surfaces are pages of $(B,\pi)$, while the red surfaces are fibers $\Sigma \times \{t\}$.  The model near a hyperbolic point, which is not shown, is a blue hyperboloid intersecting a horizontal red plane.}\label{fig:locmod1}
\end{center}
\end{figure}

Note that the pages of $(B,\pi)$ are $l$-fold covers of $\Sigma\setminus  E$ via the map from $\pi^{-1}(t) \to \Sigma\setminus  E$ defined by  $(x,z)\mapsto x$. 

\begin{example}\label{ex:obprod}  Cutting the above construction at $\Sigma\times\{0\}$ gives a contact 3-manifold $(\Sigma\times [0,l],\xi)$, with a compatible embedded foliated open book that is the restriction of $(B,\pi)$. We illustrate this in specific case where $\F_\pit$ on the torus is as given in Figure \ref{fig:efobtorus}. In the resulting foliated open book, each page is a copy of the original surface, but cut open along a different leaf of $\F_{\pit}$; respresentative pages are shown in  Figure \ref{fig:fobforfoliation}. 
 
\begin{figure}[h]
\begin{center}
\includegraphics[scale=.6]{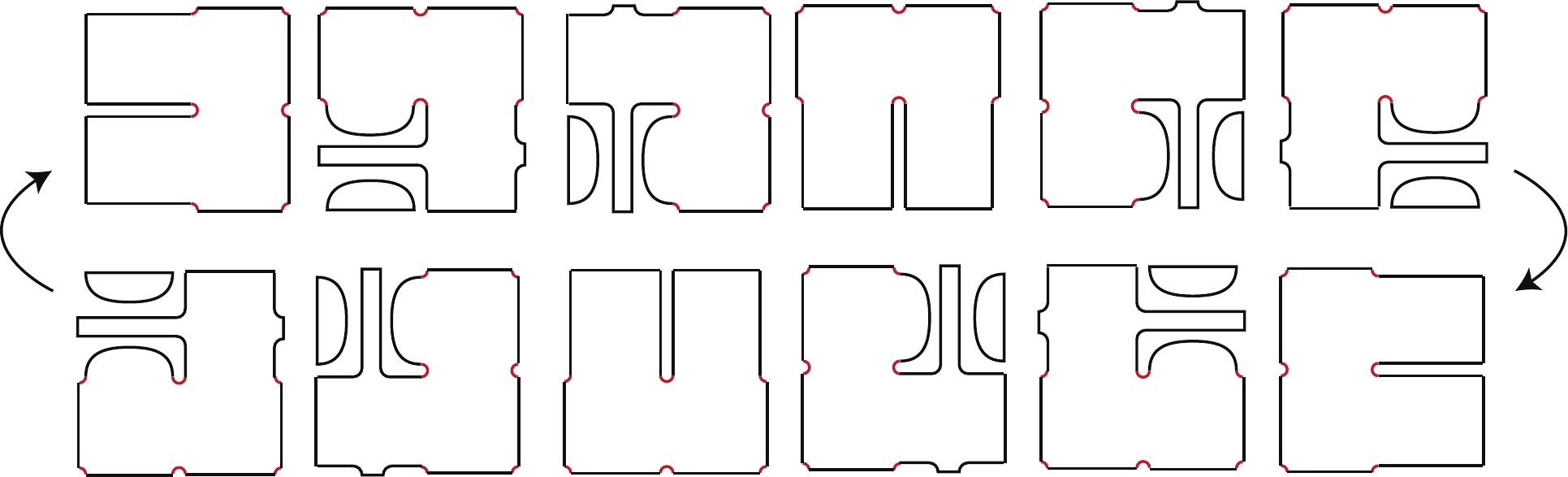}
\caption{Representative pages in the foliated open book constructed by applying Example~\ref{ex:obprod} to the foliation of the torus given by Figure~\ref{fig:efobtorus}  }\label{fig:fobforfoliation}
\end{center}
\end{figure}
In order to see why Figure~\ref{fig:fobforfoliation} shows the correct pages, Figure~\ref{fig:fobforfoliation2} indicates how a given leaf embeds in $\Sigma\times[0,1]$.  Note that away from $A=N(\Gamma)$, each leaf is covered exactly once. This builds up surfaces homeomorphic to  $R_\pm$ as $t$ varies, and these are joined across $A$ by vertically twisting bands.  

 \begin{figure}[h]
\begin{center}
\includegraphics[scale=.2]{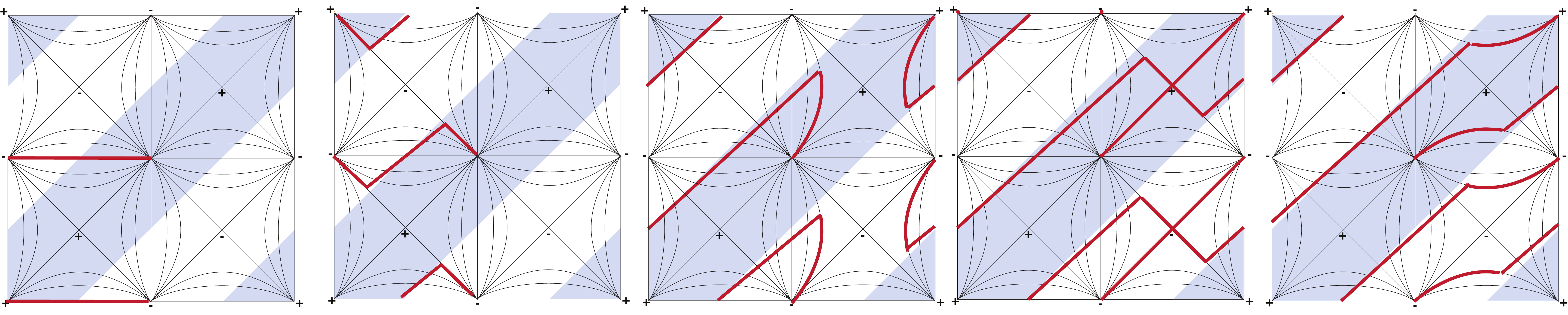}
\caption{ The red curves indicate how the page  $(\pi')^{-1}(0)$ intersects selected surfaces  $\Sigma \times \{t\}$ when Proposition~\ref{prop:obcyl}  is applied to the foliated torus from Figure~\ref{fig:efobtorus}. On each surface $\Sigma \times \{t\}$, the blue region is $R_+$ and the white region is $R_-$.  The $t$-values shown lie in $[0, \frac{1}{2}]$.}\label{fig:fobforfoliation2}
\end{center}
\end{figure}
\end{example}

\begin{corollary}[Local model]\label{cor:lm}
Suppose that $\Sigma$ is an embedded surface in a 3-manifold $M$ with an open book $(B,\pi)$ that induces an open book foliation $(\F_{\pit}, \pit,H=H_-\cup H_+)$ on $\Sigma$ with no circle leaves. Then $\Sigma$ has a neighborhood $N(\Sigma)\cong \Sigma\times [-\eta,\eta]$ so that $\pi=\pit+fz$, where $f$ is the  function from Proposition \ref{prop:Gi} that is $\pm 1$ on $R_\pm'$ and monotonically increasing in $u$ on $A=N(\Gamma)$.

Moreover, for any 1-form $\beta'$ that is (geometrically) divided by $\Gamma$, and restricts to $A$ as $gd\pit$ with $g$  as in Lemma \ref{lem:approx}, there is a contact structure $\xi$ on $M$ compatible with $(B,\pi)$ so that $\xi\vert_{N(\Sigma)}=\beta'+fdz$. \qed
\end{corollary}

\begin{remark}This local model for $\Sigma\times [-\eta, \eta]$ is itself a foliated open book which appears as a submanifold of an open book containing the surface $\Sigma$ with the induced open book foliation.
\end{remark}

Notice that  the induced open book foliations on the nearby surfaces $\Sigma\times\{z\}$ are strongly topologically conjugate  to $\F_{\pit}$. One can also understand the restrictions of each page $\pi^{-1}(t)$ to $N(\Sigma)$; the pages are homeomorphic to the surfaces $\pit^{-1}[t-\eta,t+\eta]\cap R_+'$ and $-\left(\pit^{-1}[t-\eta,t+\eta]\cap R_-'\right)$ glued together with a twisted band. If the interval $ [t-\eta,t+\eta]$ contains only regular values of $\pit$, then $\pi^{-1}(t)\cap N(\Sigma)$ is just a union of $|E_+|=|E_-|$ strips; each strip is a rectangle in $N(\Sigma)$ with a pair of its opposite sides  embedded into $\Sigma\times\{-\eta\}$ and $\Sigma\times \{\eta\}$. If $[t-\eta,t+\eta]$ contains a singular value, then one component of  $\pi^{-1}(t)\cap N(\Sigma)$ is a ``cross'', embedded as a saddle surface.  
The intersection $\pi^{-1}(t)\cap N(\Sigma)$ is singular when $t\pm \eta$ is a singular value.  
 
\subsection{Expansion}\label{sec:exp} In order to understand the relationship between Morse and embedded open books, we will modify our local model by introducing cancelling pairs of critical points.  Let  $\Sigma$ be an orientable surface embedded in a $3$-manifold $M$ whose open book structure induces an open book foliation $\F_{\pit}$ on $\Sigma$.  The above corollary provides a local model for a neighborhood $\Sigma\times[-\eta, \eta]\subset M$, which we take as the starting point for the next result.

\begin{proposition}[Expansion]\label{prop:exp} 
Let $(B_e, \pi_e)$ be the embedded foliated open book for  $\Sigma\times [-\eta, \eta]\subset M$ described in Corollary~\ref{cor:lm}. Then there is an $S^1$-valued Morse function $\pi'\colon \Sigma\times [-\eta,\eta]\setminus B \to S^1$ so that the following hold:
\begin{enumerate}
\item $\pi'$ is $C^0$-close to $\pi$ and agrees with $\pi$ outside $\Sigma\times[-\varepsilon,\varepsilon]$, for some sufficiently small $\varepsilon<\eta$;
\item for any hyperbolic point $h\in H_\pm$ of $\F_{\pit}$, the function $\pi'$ has a pair of canceling critical points $p_{h}^\pm=h\times \{\pm\varepsilon\}\in\Sigma_\pm=\Sigma\times \{\pm \varepsilon\}$;
\item the index of $p_{h}^\pm\in \Sigma_\pm$ is $2$ if the signs of $h$ and $\Sigma_\pm$ agree, and $1$ otherwise; 
\item $\pi'$ has no other critical points;
\item\label{it:new} each function $\pi'\vert_{\Sigma_\pm}$ is  Morse, and the union of their critical points is the  set of critical points of  $\pi'$; 
\item $\pi'$ defines an open book foliation $(\F_{\pi'}^\pm, \pi'\vert_{\Sigma_\pm},H\times\{\pm\varepsilon\}=H_+\times\{\pm\varepsilon\}\cup H_-\times\{\pm\varepsilon\})$ on $\Sigma_\pm$ that is strongly topologically conjugate  to $\F_{\pit}$.   
\end{enumerate} \qed
\end{proposition}
Note that condition (\ref{it:new}) implies that there is a gradient-like vector field for $\pi'$ that is tangent to both $\Sigma_+$ and $\Sigma_-$.

\begin{proof}[Proof of Proposition \ref{prop:exp}] 
Let $z'$ be the interval coordinate and $u$ and $v$ the surface coordinates in the local model for $\Sigma\times [-\eta,\eta]$ as in Corollary ~\ref{cor:lm}. 
For sufficiently small $z'$, the surfaces $\Sigma\times \{z'\}$ have induced open book foliations that are strongly topologically conjugate  to $\F_{\pit}$. We will modify $\pi$  near the singular points of  $\F_{\pit}$. For simplicity, we give the explicit computation only for $h\in H_+$, but the modification near negative hyperbolic points is similar. 

First, reparametrize the $I$ coordinate by $z$, where $z'=z^3+\rho z$ for some $\rho>0$. Note that $z\in [-\zeta,\zeta]$, where $\zeta^3+\rho\zeta=\eta$. Now define a path of functions $\pi_s$ that agrees with $\pi$ outside of a three-ball neighborhood $B_{\varepsilon}^3(h)$ of each hyperbolic point $h\in H_+$ and that is defined as follows on $B_{\varepsilon}(h)$:
\[\pi_s=(1-\Psi_{\varepsilon/2}^{\varepsilon}(r))\pi+\Psi_{\varepsilon/2}^{\varepsilon}(r)(\pit+z^3-sz)=\pit+z^3+z(\rho-(s+\rho)\Psi_{\varepsilon/2}^{\varepsilon}(r)).\]

As before,  $\Psi_{\varepsilon/2}^{\varepsilon}$ is a smooth bump function that equals $1$ for $r< \varepsilon/2$ and $0$ for $r>{\varepsilon}$, and where  $r^2=u^2+v^2+z^2$ is the 3-dimensional distance from $h$ in $B_{\varepsilon}^3$.
Clearly $\pi_{-\rho}=\pi$, and we will show that 
\begin{itemize}
\item[-] $\pi_s$ is a regular function for $s<0$; 
\item[-] $\pi_s$ is non-Morse only for $s=0$,
 with an embryonic critical point at each hyperbolic point of $\F_\pit$; and 
 \item[-] for $s>0$, $\pi_s$ is Morse with two cancelling critical points of index 1 and 2 associated to each hyperbolic point $h$ of $\F_\pit$. 
 \end{itemize}
 In order to prove  the preceding claim, in the following we determine the critical points of $\pi_s$ for a fixed $s$. Away from $D_{\varepsilon}^3(h)$, the function $\pi_s=\pi$ and thus it is regular, while its differential on $D_{\varepsilon}^3(h)$ is
\[d\pi_s=d\pit+\big(3z^2+\rho-(s+\rho)\Psi_{\varepsilon/2}^{\varepsilon}\big)dz-z(s+\rho)d\Psi_{\varepsilon/2}^{\varepsilon}dr\]

For $r>\varepsilon$ the $\Sigma$-component of  $d\pi_s$ equals $d\pit\neq 0$, and thus $d\pi_s\neq 0$.  In the thickened sphere defined by $F=\{\varepsilon/2\le r\le\varepsilon\}$, observe that the coefficients of  $du$ and  $dv$ in $d\pit$ are $2u$ and $-2v$,  respectively; thus their proportion is $u:-v$. However, in $z(s+\rho)d\Psi_{\varepsilon/2}^{\varepsilon}dr$ the proportion of the coefficients of $du$ and $dv$ is $u:v$. If in their sum both coefficients are 0, we must have $u=v=0$. Thus the possible critical points of $\pi_s$ in $F$ are all of the form $(0,0,z)$ with $\varepsilon/2\le z \le\varepsilon$. At these points, the $dz$ coordinate of the equation $d\pi_s=0$ simplifies to
\[3z^2+\rho-(s+\rho)\Psi_{\varepsilon/2}^{\varepsilon}-z(s+\rho)d\Psi_{\varepsilon/2}^{\varepsilon}=0\]
On $F$ we have $3z^2\ge \varepsilon^2/4$, while $\vert\rho-(s+\rho)\Psi_{\varepsilon/2}^{\varepsilon}\vert$ can be  bounded above by $\varepsilon^2/8$ if we choose both $s$ and $\rho$ small. 
As for the third term,
\[-z(s+\rho)d\Psi_{\varepsilon/2}^{\varepsilon}\le\varepsilon/2(s+\rho)K,\] 
where $K=\max_{\varepsilon/4\le r\le\varepsilon/2}\{-d\Psi_{\varepsilon/2}^{\varepsilon}\}>0$. Thus if we choose $|s+\rho|$ sufficiently small we can ensure that this third term is less than $\varepsilon^2/8$, and it follows that $d\pi_s$ is nonzero on $F$.

For $r<\varepsilon/2$ we have 
$d\pi_s=\d\pit+(3z^2-s)dz$. Thus for $\sqrt{\frac{s}3}<\varepsilon/2$,  there are critical points at $z=\pm \sqrt{\frac{s}3}$.

In summary, for sufficiently small $\rho, |s|>0$, the function $\pi_s$ is a Morse function for $s\neq 0$. It is regular whenever $s<0$, it has an embryonic critical point at the hyperbolic points $h$, and it has a pair of cancelling index 1 and 2 critical points for each hyperbolic point $h$ on the surfaces $\Sigma\times\{\pm \sqrt{\frac{s}3}\}$.

In the following we will examine the singularities and the leaves of the foliation induced by $\pi_s$ on $\Sigma\times\{z\}$.
Notice that if we restrict $\pi_s$ to $\Sigma\times\{z\}$, then we get a function $\pit+C$ (for $C=z^3-zs$) in $B_{\varepsilon/2}^3$ and $\pit+C'$ (for $C'=z^3+\rho z$) outside $B_{\varepsilon/2}^3$; thus the level sets of $\pit$ agree with the level sets of $\pi_s\vert_{\Sigma\times\{z\}}$ on $B_{\varepsilon/2}^3\cup\Sigma\setminus  B_{\varepsilon}^3 $. From the previous computation we  see that $\pi_s\vert_{\Sigma\times \{z\}}$ has no critical points in $F=\{\varepsilon/2\le r\le\varepsilon\}$. The fibers of $\pi_s(t)$ on $F\cap \Sigma$ connect the level sets $\pit^{-1}(a)$ on $D^2_{\varepsilon/4}$ to the fibers of $\pit(a-(\rho+s)z)$ of $\Sigma\setminus \cup_h D^2_{\varepsilon/2}$. 
This means that if we choose both $|z|$ and $|\rho+s|$ sufficiently small, the foliation induced by
 $\pi_s\vert_{\Sigma\times\{z\}}$ will be strongly topologically conjugate to $\F_{\pit}$.

Now fix $s>0$  and $\rho>0$ sufficiently small so that the results above hold and let $\pi'=\pi_s$, $\Sigma_\pm=\Sigma\times\{\pm \sqrt{\frac{s}3}\}$. 
The computation above shows that the restriction of $\pi'$ to $\Sigma_\pm$ is Morse, since the original $\F_{\pit}$ was an open book foliation.   Since  the levels sets of $\pi'\vert_{\Sigma_\pm}$ are strongly  topologically conjugate to $\F_{\pit}$, there are no other critical points of $\pi'\vert_{\Sigma_\pm}$, as required for Item~\ref{it:new}.

\end{proof}

\begin{corollary} The function $\pi'$ from Proposition~\ref{prop:exp} induces a Morse foliated open book on each  component of $M$ cut along $\Sigma_\pm$.  \qed
\end{corollary}

\subsection{Local Models for Morse Foliated Open Books}
The techniques above may be adapted to the case of Morse foliated open books.  We carefully state the result we will rely on in the next section, but leave the  analogous proof to the reader.

Given 3--manifolds $M$ and $M'$ with properly embedded 1--manifolds $B$ and $B'$, respectively, suppose  that  $\pi:M\setminus B\rightarrow S^1$ and $\pi':M'\setminus B'\rightarrow S^1$  are circle-valued Morse functions.  Let $\Sigma\hookrightarrow M$ and $\Sigma'\hookrightarrow M'$ be embedded surfaces with induced open book foliations $(\F_{\pit}, \pit,H=H_+\cup H_-)$, $(\F_{\pit}',\pit', H'=H'_+\cup H'_-)$ with no circle leaves; we require further that $\pi$ and $\pi'$ have  gradient-like vector fields tangent to the  surfaces $\Sigma$ and $\Sigma'$ and thus any hyperbolic points of these open book foliations coincide with critical points of the corresponding Morse functions on the ambient manifold. Note that we allow $\Sigma$ to embed into $M$ or $M'$  as the boundary.

\begin{proposition}[Morse Neighborhood Theorem]\label{prop:mnbrh} 
If there is a diffeomorphism $\psi\colon \Sigma\to \Sigma'$ that takes the two open book foliations to each other, then  there are one-sided neighborhoods $N_\pm=N_{\pm}(\Sigma)$ and $N_{\pm}'=N_{\pm}(\Sigma')$ and an extension of $\psi$ to $\Psi_\pm\colon N_{\pm}\to N_{\pm}'$ so that 
\begin{enumerate}
\item $\Psi_{\pm}(B)=B'$;
 \item $\pi=\pi'\circ\Psi$;
 \item $\Psi_{\pm}(B\cap N_{\pm})=B'\cap N_{\pm}'$ and the trivialisation of $N(B)\cap N_{\pm}\cong (B\cap N_{\pm})\times D^2$ in which $\pi=\vartheta$ maps to the trivialisation of $N(B')\cap N_{\pm}'\cong (B'\cap N_{\pm}')\times D^2$ in which $\pi'=\vartheta$. \qed
 \end{enumerate}

\end{proposition}

We consider one-sided neighborhoods for two reasons.  First, if the surfaces containing the critical points lie in the interior of $M$, we will in practice  always cut along them require the local models only after cutting.  Second, if we restrict to one-sided neighborhoods, then proof of Theorem~\ref{prop:nbrh} may be modified to prove Proposition~\ref{prop:mnbrh} simply by selectively replacing $z$ by $\pm z^2$ when constructing adapted coordinate systems near hyperbolic points.


\section{Equivalence of foliated open books}\label{sec:fobequiv} 

 In this section, we use the local models constructed in Section~\ref{sec:models}  to relate embedded, abstract, and Morse foliated open books.   After this section, we will move freely between these different notions of foliated open books and use  whichever is most convenient.
 
 \subsection{{From embedded to Morse and back}}\label{sec:mem}
 Away from a neighborhood of the boundary, embedded foliated open books and Morse foliated open books are indistinguishable.  In this section, we apply results from the previous section to show how to transform one type to the other.  
 
Starting first with an embedded foliated open book $\E=(B_e, \pi_e, \F_{\pit_e})$, identify $\partial M$ with $\Sigma\times\{0\}$ in $\Sigma\times [-\eta, 0]$, half of  the standard neighborhood constructed in Proposition~\ref{prop:nbrh}.  Define the Morse foliated open book $\M(\E)=(B, \pi', \F_{\pit'})$ on $M\setminus \Sigma\times [-\varepsilon,0]$ to be the result of first modifying $\pi_e$ to $\pi'$ as in Proposition \ref{prop:exp} and then removing  $\Sigma\times [-\varepsilon,0]$.

Similarly, suppose that $\M=(B_m, \pi_m, \F_{\pit_m})$ is a Morse foliated open book, and identify $\partial M$ with $\Sigma\times\{0\}$ in the standard half-neighborhood $\Sigma\times [-\eta, 0]$ from Proposition~\ref{prop:mnbrh}.  Define the embedded foliated open book $\E(\M)=(B, \pi, \F_{\pit})$ on $M\setminus \Sigma\times [\frac{-\eta}{2},0]$ to be the result of removing  $\Sigma\times [\frac{-\eta}{2}, 0]$. 
 
\begin{proposition}\label{prop:eandm2} With the notation above 
\begin{enumerate}

\item\label{it:egy} $\M(\E)$ is a Morse foliated open book.  Furthermore, there is a diffeomorphism $\psi^{\E\to \M}\colon M\to M\setminus \Sigma\times [-\varepsilon,0]$ that takes the bindings and the pages of the two foliated open books to each other. In particular, for any contact structure $\xi$ supported by $\E$ the contact structure $\psi^{\E\to \M}_{*}\xi$ is supported by $\M(\E)$;
\item\label{it:ketto}   $\E(\M)$ is an embedded foliated open book.  Furthermore, there is a diffeomorphism $\psi^{\M\to \E}\colon M\to M\setminus \Sigma\times [-\frac{-\eta}2,0]$ that takes the bindings and the pages of the two foliated open books to each other. In particular, for any contact structure $\xi$ supported by $\M$ the contact structure $\psi^{\M\to \E}_{*}\xi$ is supported by $\E(\M)$;
\item\label{it:harom} The embedded foliated open books $\E$ and $\E(\M(\E))$ are diffeomorphic;
\item\label{it:negy} The Morse foliated open books $\M$ and $\M(\E(\M))$ are diffeomorphic.
\end{enumerate}

\end{proposition}

\begin{proof}

For the proof of Item (\ref{it:egy}),  recall from Proposition~\ref{prop:exp} that $\pi'$ has exactly two critical points for each hyperbolic singularity of $\F_{\pit_e}$, and only one of each pair survives as a boundary critical point of $\M(\E)$.  Thus $\pi'$ satisfies the conditions of Definition~\ref{def:mfob}.  The required diffeomorphism will be constructed as composition $\phi_\E\circ \phi_\M^{-1}$, where $\phi_\E$ and $\phi_\M$ are each diffeomorphisms from the respective foliated open books to the common manifold $\E\setminus (-\eta,0]$.  Let $\varepsilon<\eta'<\eta$ be a real number such that $\pi=\pi'$ on $\Sigma\times[-\eta, -\eta']$, again viewed inside the neighborhood of the boundary in $\E$ and $\M(\E)$.  Choose diffeomorphisms $\phi_*$ that map $\Sigma \times [-\eta,0]$ (respectively, $\Sigma\times[-\eta, -\epsilon]$) to   $\Sigma \times [-\eta, -\eta']$.  Away from the critical points, these may be chosen to restrict to each page as an isotopy to a proper subpage; since the open book foliation on $\Sigma\times \{z\}$ is strongly topologically conjugate to the open book foliation on the boundary of each of $\E$ and $\M(\E)$,  the induced isotopy of plane fields takes a supported contact structure to a supported contact structure. 

For Item (\ref{it:ketto}), we observe first that the restriction of $\pi_m$ to the complement of the boundary of $\M$ is regular, and as the induced open book foliation $\F_{\pit}$ is strongly topologically conjugate to $\F_{\pit_m}$, the conditions of Definition~\ref{def:efob} are satisfied.  The proof of the second half of the claim proceeds as in the proof of Item (1), via a composition of isotopies from the original manifolds to a common submanifold.  

The proofs of Items (\ref{it:harom}) and (\ref{it:negy})  proceed as the second half of Item (\ref{it:egy}).
\end{proof}

\subsection{From Morse to abstract and back}\label{sec:morsetoabs} 
Suppose that $(B_m,\pi_m, \F_{\pit_m})$ is a Morse foliated open book for $(M,\F_{\pit})$ and fix a gradient-like vector  field $\nabla\pi_m$ for $\pi_m$.  Let $M'$ be the  3-manifold  formed by  blowing up $B$ to $S^1\times B$. Extend $\pi_m$ to $M'$ continuously, now calling the extension $\pi$. A local model near the binding allows us to assume that $\nabla \pi=\partial t$ on the blown-up $S^1\times B$ boundary component; here, $t$ parameterizes $S^1$.
Let $\{h_i\}_{i=1}^{2k}$ be the set of critical points of $\pi$  with corresponding critical values $t_1^*=\pi(h_1)<\dots<t_{2k}^*=\pi(h_{2k})$ with respect to the cyclic order on $S^1$. Choose regular points $\{t_i\}_{i=0}^{2k}$  so that $t_0<t_1^*<t_1<\dots<t_{2k}^*<t_{2k}<t_0$. Let $S_i={\pi}^{-1}(t_i)$. This means that $S_i$ is a surface with cornered boundary $(B\times \{t_i\})\cup \alpha_i$, 
where $\alpha_i={\pit^{-1}(t_i)}$ and the corners are $E\times\{t_i\}=(B\times\{t_i\})\cap\alpha_i$.  Note that as we are working in the blow-up of $M$, $S_i$ is closed.

If $h_i$ is an index 2 critical point of $\pi$, then its stable manifold $W(h_i)^s$ intersects $S_{i-1}$ in a properly embedded arc  
$\gamma_{i-1}^+$ with endpoints on $\alpha_{i-1}$, and $S_{i}$ is obtained from $S_{i-1}$ by  cutting along $\gamma_{i-1}^+$ and smoothing. 
In the notation of Definition \ref{def:afob}, we have $S_{i-1}\xrightarrow[{\gamma_{i-1}^+}]{\textbf{cut}}S_i$.

If $h_i$ is a critical point of index 1, then its stable 
manifold $W(h_i)^s=w(h_i)^s$ intersects $S_{i-1}$ in two points $\{p_{i-1},q_{i-1}\}\in \alpha_{i-1}$, and $S_{i}$ is obtained from $S_{i-1}$ by  
gluing a 1-handle along this attaching sphere. In the notation of Definition \ref{def:afob}, we have $S_{i-1}\xrightarrow[p_{i-1}, q_{i-1}]{\textbf{add}}S_i$.

The flow of the gradient-like vector field gives a diffeomorphism $h\colon S_{2k}\to S_0$ that is constant on $B$. 

\begin{proposition}[from Morse to abstract]\label{prop:mtoa} With the notation from above:
\begin{enumerate}
\item Any Morse foliated open book $\mathcal{M}=(B_m,\pi_m, \F_{\pi_m})$ and  choice of gradient-like vector field  $\nabla\pi_m$ determines an abstract foliated open book $\mathcal{A}=(\{S_i\},h)$ as above. 

\item If $\nabla\pi_m$ and $\nabla'\pi_m$ are both  gradient-like vector fields for $\pi_m$, then the corresponding abstract foliated open books $\mathcal{A}=(\{S_i\},h)$ and  $\mathcal{A}'=(\{S_i\},h)$ are diffeomorphic. 

\item If some reparameterization of the Morse foliated open book $\mathcal{M}'=(B_m',\pi_m', \F'_{\pi_m})$ is diffeomorphic to  $\mathcal{M}=(B_m,\pi_m, \F_{\pi_m})$ then the corresponding abstract foliated open books $\mathcal{A}$ and $\mathcal{A}'$ (for any choices of gradient-like vector field) are conjugate to each other. 
\end{enumerate}
\end{proposition}

\begin{proof} The proposition is a consequence of the construction above. For the first claim, we see that the critical submanifolds of a fixed critical point determine a pair of abstract pages related by a single handle addition or deletion.  In order to assemble these pairs into an abstract open book, it remains to identify the two copies of $S_i$ associated to the critical points $t_i^*$ and $t_{i+1}^*$, but this is accomplished by the flow of the gradient-like vector field.  Although the choice of gradient-like vector field could change the  arcs $\gamma_{i-1}$ or the attaching points $p_{i-1},q_{i-1}$, the diffeomorphism types of the pages $S_i$ remain the same; this establishes the second claim. The last claim is simply the translation of reparameterization (specifically, by rotation) to the abstract setting.  
\end{proof}

In the rest of this subsection we describe how to use abstract foliated open books to build up a 3-manifold with foliated boundary together with an embedded foliated open book. 

\begin{figure}[h]
\begin{center}
\includegraphics[scale=0.8]{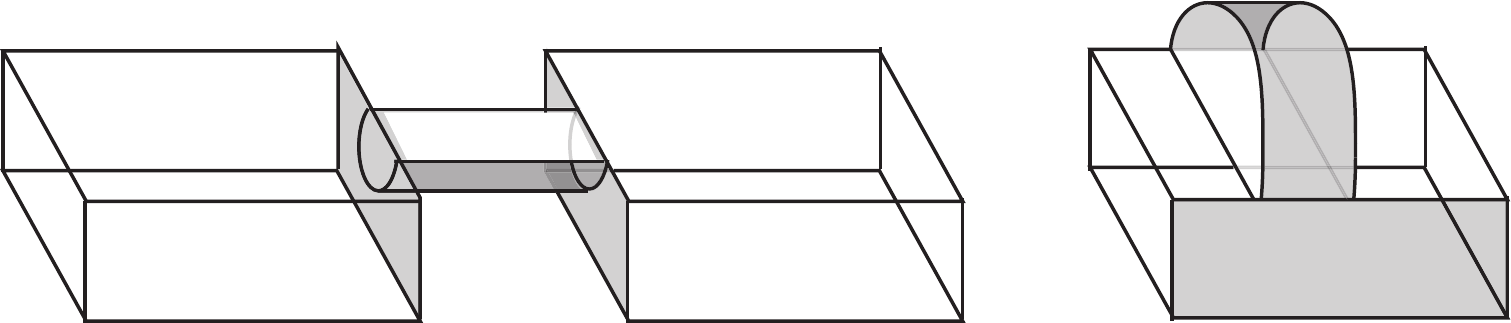}
\caption{The cobordism between $S_{i-1}$ and $S_{i}$, shown before smoothing the attached half handles.  The shaded regions are in the foliated boundary.  Left: an index one boundary critical point corresponding to half a one-handle or the ``add'' operation.  Right: an index two boundary critical point corresponding to half a two-handle or the ``cut'' operation.}\label{fig:cobordism}
\end{center}
\end{figure}

To each handle addition and deletion, we associate a cornered handlebody $H_i$   We first consider the case of $S_{i-1}\xrightarrow[p_{i-1}, q_{i-1}]{\textbf{add}}S_i$. Define $H_i$ to be the cornered handlebody   obtained by attaching a 3-dimensional boundary 1-handle to  $S_{i-1}\times I$. See Figure \ref{fig:cobordism}. 

  By construction,  $H_i$ is a cornered cobordism between $S_{i-1}$ and $S_{i}$; thus $H_i$ is a handlebody with boundary $\partial H_i=-S_{i-1}\cup V_i\cup S_{i}$ and with codimension-1 corners at $\partial (-S_{i-1})\cup \partial S_{i}\cup (E\times I)$ and codimension-2 corners at $E \times \partial I$.   The vertical boundary $V_i$ is the union of two types of components: 
\begin{itemize}
\item[-] products of cornered circles with the interval $I$;
\item[-] a (cornered) pair of pants $P_i$ . 
\end{itemize}
Components of the first type  correspond to components of $\partial S_{i-1}$ which are preserved in $\partial S_{i}$, while  the second type is associated to the handle addition.

In the case of $S_{i-1}\xrightarrow[\gamma_{i-1}]{\textbf{cut}}S_i$, then the cobordism $H_i$ may be built up analogously.  In this case, we  attach a 3-dimensional boundary 2-handle to the thickened $S_i$ along $\gamma_{i-1}$; this turns the previous construction upside down and the boundary and corners may be identified analogously. 
Each handlebody $H_i$  is naturally equipped with a Morse function $\pi_i\colon H_i\to I$ with the following properties: 
\begin{itemize}
\item[-] $\pi_{i}\vert_{-S_{i-1}}\equiv 0$;
\item[-] $\pi_{i}\vert_{S_{i}}\equiv 1$;
\item[-] $\pi_i$ has a  gradient-like vector field $\nabla\pi_i$ that is tangent to $V_i$ and to the $I$-component at $B\times I$;
\item[-]  $\pi_i$ has a unique critical point with critical value $1/2$ and this is located on the pair of pants component of the boundary. The index is $1$ in the case of  handle addition and $2$ in the case of  saddle resolution.
\end{itemize}

We may glue the pairs $(H_i, \pi_i)$  along the $S_i$ via the identity and  glue $S_{2k}$ to $S_{0}$ via $h$ to get a 3--manifold $M'$ with boundary. After rescaling, the maps $\pi_i$ glue to a map $\pi\colon M\to S^1=\R/\mathbb{Z}$. Finally, collapse $B\times S^1$ to $B$ to obtain the 3-manifold $M=M(\{S_i\},h)$. Then $\pi$ restricts as a map $\pi\colon M\!\smallsetminus\! B\to S^1$ and level sets of $\pit=\pi|_{\partial M}$  induce a  foliation $\F_{\pit}$ on $M$. Note that different  choices of $\pi_i$ and  scaling yield different parametrizations of $\F_{\pit}$; thus abstract foliated open books describe diffeomorphism classes of 3-manifolds with foliated boundaries $(M,\F)$. The discussion above can be summarized by the following proposition.

\begin{proposition}\label{prop:atom} With the notation and definitions from above:
\begin{enumerate}
\item Any abstract foliated open book $\A=(\{S_i\},h)$ defines up to diffeomorphism a  3-manifold $(M(\A),\F(\A))$ with foliated boundary and a Morse foliated open book $\M(\A)=(B_m(\A),\pi_m(\A), \F_m (\A))$ compatible with it; 
\item If the abstract foliated open books $\A$ and $\A'$ are conjugate to each other, then the above diffeomorphism-types of  manifold with Morse foliated open book $\big(M(\A),\F(\A), \F_m(A)\big)$, and $\big(M(\A'),\F(\A'),\F_m(\A')\big)$ are the same;
\item For any abstract foliated open book, $\A(\M(\A))$ is equivalent to $\A$;
\item For any Morse foliated open book,  $\M(\A(\M))$ is diffeomorphic to $\M$. 

\end{enumerate}
\end{proposition}

The foliation $\F_{\pit}$ can be described directly from the data of an abstract open book. The regular leaves of $\F_{\pit}$ are $\alpha_i\times \{t\}$ for $t\in (t_i^*,t_{i+1}^*)$, and each leaf is an  interval that connects elliptic points.  Each singular leaf corresponds to a critical value $t_i^*$.

The discussion above allows us to state the obvious definition of a contact structure supported by an abstract open book.

\begin{definition} Suppose that $\A=(\{S_i\}, h)$ is an abstract foliated open book defining a manifold $(M(\mathcal{A}), \mathcal{F}(\mathcal{A}))$.  The contact structure $\xi$ is \emph{supported} by $(\{S_i\}, h)$ if $\xi$ is supported by a Morse foliated open book $(B_m(\mathcal{A}), \pi_m(\mathcal{A}), \F_{ \pit_m(\mathcal{A})})$ compatible with the abstract data.  
\end{definition}

Together with the statements of Proposition~\ref{prop:eandm2} we have the following equivalence for a fixed diffeomorphism class of 3-manifolds $(M,\F)$ with foliated boundary:
\[\frac{\left\{\text{abstract FOBs}\right\}}{\text{conjugacy}} \leftrightarrow \frac{\left\{\text{embedded FOBs}\right\}}{\text{diffeomorphism}}  \leftrightarrow \frac{\left\{\text{Morse FOBs}\right\}}{\text{diffeomorphism}}
\]

From now on we will move freely between these notions, using whichever is most convenient in our discussion. 

\subsection{Sorted handlebodies}\label{sec:sort}

In this final subsection, we further examine the relationship between Morse and abstract foliated open books.  The key definition of a \emph{sorted} foliated open book will be important in Section~\ref{sec:ptofob}.  Throughout this section, we will consider a Morse foliated open book equipped with a fixed  preferred gradient-like vector field $\nabla \pit$ on the boundary, and we indicate this choice by the following tuple: $(B, \pi, \F_{\pit}, \nabla \pit)$.   (Recall that ``preferred'' was defined in Definition~\ref{def:pref}.)

As in Section~\ref{sec:morsetoabs}, let $M'$ denote the manifold formed by blowing up $B$ into $S^1\times B$.
  We denote the the fibers $\pi^{-1}(t)$ in $M'$ by $S_t$, and we note that these are closed and disjoint.
Assume that $t=0$ is a regular point of $\pi$, and  let  $\overline{M}$ denote the formal closure of $M'\setminus S_0$, which is a cobordism between $S_0$ and $S_1$.  
We call a gradient-like vector field $\nabla \pi$ on $M\setminus B$ \textit{preferred} if it extends $\nabla \pit$, and a  preferred vector field induces a vector field ---still denoted $\nabla \pi$--- on  $\overline{M}$.
In the following we will  work exclusively in $M'$ or $\overline{M}$. Informally these are the submanifolds one obtains when gluing together some or all consecutive handlebodies $H_i$ of the previous section.

In the definition below, assume that $t<t'$ are regular values for $\pi:M'\rightarrow S^1$.
\begin{definition}\label{def:sortedcobord} 
The submanifold $_tH_{t'}:=\pi^{-1}(t, t')$ is \emph{sorted}  if the stable and unstable manifolds of critical points of $\nabla\pi\vert_{{}_t H_{t'}}$  
are disjoint in $_tH_{t'}$. The Morse foliated open book $(B, \pi, \F_{\pit}, \nabla \pi)$ is \emph{sorted} if $\overline{M}$ is a sorted submanifold.
\end{definition}
We identify $\overline{M}$ with $_0H_{1}$, and the statements below hold for $t=0, t'=1$, as well. 
\begin{remark}The fact that $\nabla\pi$ is preferred forces an order on the stable and unstable submanifolds near the boundary, which may obstruct this disjointness. See Example~\ref{ex:notsep}. 
\end{remark}

\begin{definition}\label{def:r} 

Define $_tR_{t'}\subset {}_tH_{t'}$ to be the minimal $\nabla \pi\vert_{{}_tH_{t'}}$-invariant subset of ${}_tH_{t'}$  containing a cornered neighbourhood of $\partial M\cap {}_tH_{t'}$ that is disjoint from $[t,t']\times \mathit{int}(B)$. Define $R={}_0R_{1}\subset \overline{M}$ and  let $R_{t}=R\cap S_t$.
\end{definition} 

Since the critical points of $\pi\vert_{{}_tH_{t'}}$  all lie on $\partial M$, the invariant subset $_tR_{t'}$ must contain all the critical submanifolds $W^u(h)$ and $W^s(h)$ for critical points $h$ with critical value between $t$ and $t'$. This implies that the complement of $_tR_{t'}$ is a product:

\begin{lemma}\label{lem:ppersists}  Suppose that $(B, \pi, \F_{\pit}, \nabla \pi)$ is a sorted Morse foliated open book.  Then ${}_tH_{t'}\setminus {}_tR_{t'}$ is diffeomorphic to the product $(S_t\setminus R_t)\times[t,t']$.

The subsets $S_{t}\setminus R_{t}$ are isotopic for all values of $t$.

\end{lemma} 

\begin{proof} 
The restriction of $\nabla \pi$ to ${}_tH_{t'}\setminus {}_tR_{t'}$ has no critical points, so the flow of $\nabla \pi$ defines a diffeomorphism between level sets. 
\end{proof}

\begin{figure}[h!]
\begin{center}
\includegraphics[scale=0.9]{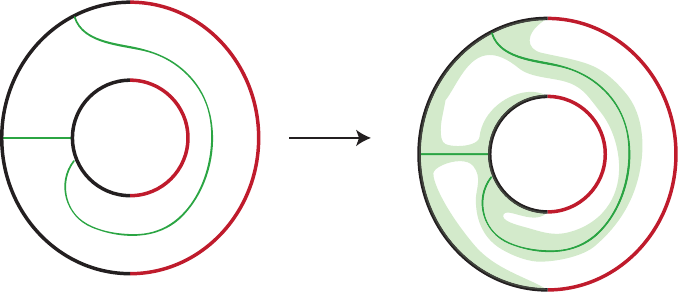}
\caption{ Left: $S_0$ from the sorted abstract foliated open book from Example~\ref{ex:torusstab}, shown with its two $\gamma_i^+$ curves. Right: In the corresponding Morse foliated open book, the shaded region shows $R_0$ and the white region shows $S_0\setminus R_0$. }\label{fig:ttopob}
\end{center}
\end{figure}

The discussion  in the previous subsection allows us to characterize sorted open books in the abstract setting, as well.

When a preferred $\nabla \pi$ is fixed, the  critical submanifolds  define a collection of curves on each regular page in $_tH_{t'}$.  We call these cutting and co-core arcs \textit{sorting arcs}.  Slightly abusing notation, we will let $\gamma_k^{{\pm}}$ denote a sorting arc associated to $h_k$ on any page. With no further conditions on $\nabla \pi$ imposed, observe that these arcs may intersect each other,  and after passing a critical point whose corresponding submanifold  intersects another critical submanifold, some arcs may be cut.  If $h_k$ is a positive hyperbolic point in a foliated open book, then we write $W^s(h_k)\cap S_t=\gamma^+_k$  on all pages $S_t$ for $t<\pi(h_k)$; if $h_k$ is a negative hyperbolic point, then $W^u(h_k)\cap S_t=\gamma^-_k$ on all pages $S_t$ for $t>\pi(h_k)$.  

When $\mathcal{A}=(\{S_i\},h)$ is the abstract open book associated to a Morse foliated open book as in Proposition~\ref{prop:mtoa}, we may  record the sorting arcs $\gamma_k^\pm$ on the abstract pages; we  denote the associated abstract  foliated open book by $(\{S_i\}, h, \{\gamma_k^\pm\} )$ when we want to keep track of this extra data, in parallel with the quadruple used for a sorted Morse foliated open book.   Similarly, we may write $_iH_j$ for the handlebody constructed from the abstract pages $S_k$ with $k\in \{i, i+1, \dots, j\}.$  This notation will be used in Proposition~\ref{prop:sorthandlebody} and again in Proposition~\ref{prop:pobfobequiv}.

When $(B, \pi, \F_{\pit}, \nabla \pit)$ is a sorted Morse foliated open book, the sorting arcs on the associated abstract open book $(\{S_i\}, h, \{\gamma_k^\pm\} )$ are disjointly embedded, and vice versa.  As the next example shows, requiring a Morse foliated open book to be sorted may bound the Euler characteristic of the pages from above.

\begin{example}\label{ex:notsep}
The foliated open book for the solid torus introduced in Example~\ref{ex:torus} is not sorted, as shown in Figure~\ref{fig:notsep}.  Definition~\ref{def:pref} dictates the order of the attaching spheres and endpoints of the cutting arcs along each component of $\alpha_0$, and we see that as $t$ increases, it becomes impossible to connect pairs of blue and green dots by disjoint arcs in the interior of the page.   Example~\ref{ex:torusstab} shows a sorted foliated open book for the same $\F_{\pit}$.
\begin{figure}[h!]
\begin{center}
\includegraphics[scale=0.8]{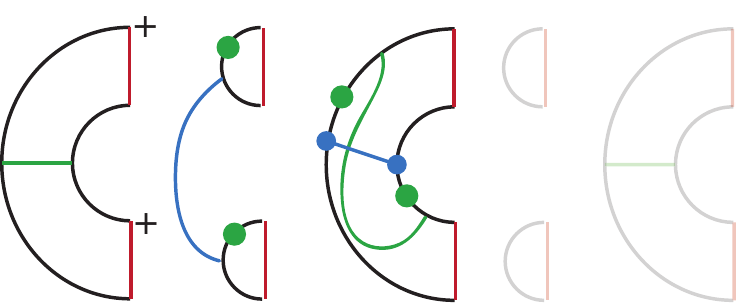}
\caption{ In $S_2$, the specified endpoints cannot be connected by disjoint sorting arcs, so the foliated open book is not sorted.}\label{fig:notsep}
\end{center}
\end{figure}
\end{example}



\section{Operations on open books: cutting, gluing, and stabilization}\label{sec:operations}
In this section, we describe natural operations on open books of various sorts. An initial  cutting result was presented in Example~\ref{prop:cut}, and here we extend this to consider the contact structure supported by an open book cut along a surface with an open book foliation.  We also explain how our focus on boundary foliations simplifies gluing, and we use this to define  stabilization for embedded foliated open books and Morse foliated open books.  We also discuss stabilization for abstract open books, which  is analogous to the operation familiar from closed manifolds.

\subsection{Cutting}\label{ssec:cutting}
One of the key advantages of  foliated open books is that they are natural structures with respect to cutting and gluing. 

\begin{theorem}[Cutting]\label{thm:cut} 
Given a 3-manifold $M$ with an open book $(B,\pi)$, suppose that an embedded surface $\Sigma$ admits an open book foliation $\F_{\pit}$ with no circle leaves. Then there is a contact structure $\xi$ on $M$ supported by $(B,\pi)$ such that the restriction of $\xi$ to the closure of each component of $M\setminus \Sigma$ is supported by the embedded foliated open book obtained as the restriction of $(B,\pi)$ onto the closure of this component. 
\end{theorem}

In fact, this result and others in this section apply equally well when the ambient manifold is a  partial or foliated open book rather than an (honest) open book; our choice to state them narrowly avoids excessive and unpleasant notation.  

\begin{proof}
Fix a contact structure $\xi$ as in Proposition \ref{prop:IK} so that the characteristic foliation $\F_\xi$ on $\Sigma$ is strongly topologically conjugate to $\F_{\pit}$ on $\Sigma$. Then the boundary conditions of Definition \ref{def:supportxi} are automatically satisfied by the components of the closures of the cut-open manifolds, so the restriction of $\xi$ is indeed supported by the restriction of the open book on this component. 
\end{proof}

This cutting can also be understood in the abstract setting. For ease of notation, we further restrict our attention to the case when $\Sigma$ is connected and separates $M$ into $M_L\cup M_R$.

As in Proposition \ref{prop:exp}, we first modify $\pi$ to $\pi'$ in a model neighbourhood $\Sigma \times [\eta, \eta]$  so that canceling critical points are created on $\Sigma \times \{\pm \epsilon\}$.  Remove   $\Sigma \times (\epsilon, \epsilon)$  to get  Morse foliated open books on $M_L'=M_L\setminus \Sigma\times (\epsilon, -\eta]$ and $M_R'=M_R\setminus \Sigma\times [\eta, \epsilon)$. Then we construct the abstract foliated open books $(\{S_i^L\}_{i=1}^{2k},h^L)$ and $(\{S_i^R\}_{i=1}^{2k},h^R)$
corresponding to $(B\cap M_L',\pi'\vert_{M_L'})$ and  $(B\cap M_R',\pi'\vert_{M_R'})$. 
Note that   the two Morse foliated open books have pairs of critical points with the same critical values, but always of opposite type (i.e., one each of index $1$ and $2$).   Then the number $k$  of sequentially distinct pages in the two abstract foliated books is the same.

Moreover, through the two deformation retractions that recover the embedded foliated open book pages from  the Morse foliated open book pages, we can identify the boundaries $\alpha_i^R$ and $\alpha_i^L$ with each other (through $\pi^{-1}(t_i)\cap \Sigma$), and we will call this simply $\alpha_i$.   With this identification,  we can view each $\alpha_i$ as a properly embedded separating arc in $S$,  rather than as  two parallel arcs  $\alpha_i^R$ and  $\alpha_i^L$.  Similarly, we may view $\gamma_{i-1}$ as  embedded on $S^R_{i-1}\subset S$ or $S^L_{i-1}\subset S$. If $S_{i-1}^L\xrightarrow[\gamma_{i-1}]{\text{\textbf{cut}}} S_i^L$, then $S_{i-1}^R\xrightarrow[\partial\gamma_{i-1}]{\text{\textbf{add}}} S_i^R$, and vice versa.  In this case, $h^L$ and $h^R$ are just the restrictions of $h$ to $S_{2k}^L$ and $S_{2k}^R$, respectively. This perspective allows us to  formulate cutting solely in terms of the cutting arcs $\alpha_i$ on $S$, as follows.

\begin{figure}[h!]
\begin{center}
\includegraphics[scale=0.6]{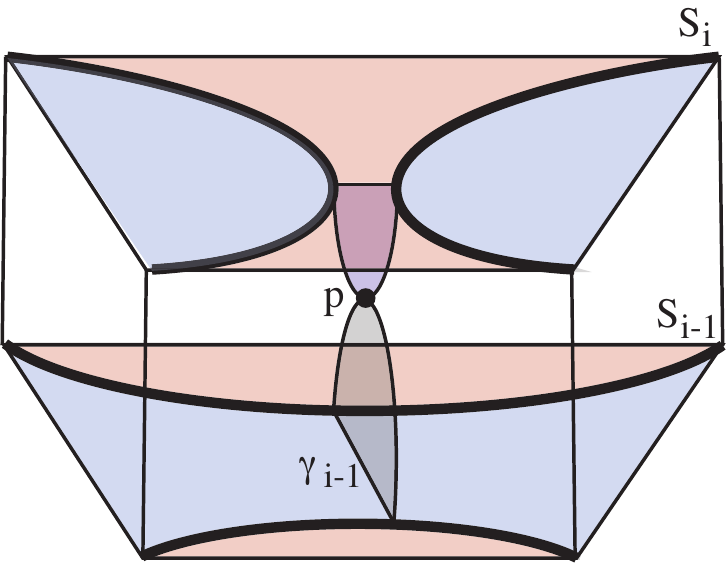}
\caption{ The red $S_{*}^R$ and the blue $S_{*}^L$ are separated by the bold $\alpha_*$. Note that $p$ is an embryonic critical point.  }\label{fig:cut}
\end{center}
\end{figure}

Given  an abstract open book $(S,h)$ for $M$, define a \emph{sequence of cutting arcs} as a set of properly embedded separating arcs $\{\alpha_i\}_{i=0}^{2k}$ for some $k$, so that
\begin{itemize}
\item[-] $h(\alpha_{2k})=\alpha_0$; and
\item[-] $\alpha_{i+1}$ is obtained from $\alpha_i$ by surgery along an arc $\gamma_i$ that  intersects $\alpha_i$ only at its boundary.
(This means that $\alpha_{i+1}$ is the smoothing of $\alpha_i \!\smallsetminus\! \partial \gamma_i$ and two parallel copies of $\gamma_i$.) See Figure \ref{fig:cut}.
\end{itemize}  

Then $S\setminus \alpha_i=\mathring{S}_i^R\cup \mathring{S}_i^L$, where $S_i^R$ and $S_i^L$ are cornered surfaces and the boundary orientation of $S_i^L$ ($S_i^R$) recovers the orientation of $\alpha_i$ ($-\alpha_i$). The diffeomorphism $h$ restricts to $S_{2k}^L$ and $S_{2k}^R$ as $h^L_{2k}$ and $h^R_{2k}$.

\subsection{Gluing}\label{ssec:glue}
The cutting operation for abstract open books can be reversed to glue abstract foliated open books with compatible boundaries. As above, for notational convenience we restrict  to the case of  gluing two abstract foliated open books along their respective connected boundaries to obtain an (honest) open book.   

 Let $(\{S_i^L\}_{i=0}^{2k},h^L)$ and $(\{S_i^R\}_{i=0}^{2k},h^R)$ be abstract foliated open books, and  assume that there is an orientation reversing pairing $\psi$ of the corners $\partial B^L$ and $\partial B^R$ so that:
\begin{itemize}
\item[-] $S_{i-1}^L\xrightarrow{\text{\textbf{cut}}} S_i^L$  if and only if $S_{i-1}^R\xrightarrow{\text{\textbf{add}}} S_i^R$;
and vice versa;
\item[-] the components of $\alpha_{i-1}^L$ containing the attaching sphere (for the handle or for the cutting arc) on $S_{i-1}^L$ have endpoints on $\partial B^L$ which are paired with the endpoints of  the components of $\alpha^R_i$  containing the attaching sphere on $S_i^R$.
\end{itemize}

The final condition allows us to glue the pages $S_i^L$ and $S_i^R$ along orientation-reversing maps  $\psi_i\colon \alpha_i^L\to \alpha_i^R$ that extend $\psi$ so that the core of a handle added to yield $S_i^L$ is identified with the cutting arc of $S_{i-1}^R$, and vice versa. This yields a well defined page $S=S_i^L\cup_{\psi_i}S_i^R$, and the monodromy $h=h^L\cup h^R$ gives an honest abstract open book. Translating the above to embedded open books and using Theorem \ref{thm:cut} gives the following:

\begin{theorem}[Gluing]\label{prop:glue} 
Suppose that the embedded foliated open books $(B^L,\pi^L)$ and $(B^R,\pi^R)$ define the 3-manifolds with foliated boundary
 $(M^L,\F_{\pit^L})$ and $(M^R,F_{\pit^R})$, and assume that there is an orientation reversing diffeomorphism $\varphi\colon \partial M_L\to \partial M_R$ that takes the foliation $\F_{\pit^L}$ to (a possible reparameterization of) $\F_{p(\pit^R)}$. 

Then there are contact structures $\xi^L$ and $\xi^R$ compatible with  $(B^L,\pi^L)$ and $(B^R,p(\pi^R))$, respectively, so that $\xi=\xi^L\cup_\varphi\xi^R$ is a contact structure $\xi$ on the glued-up manifold $M=M^L\cup_\varphi M^R$  that is compatible with the glued-up honest open book
\[(B,\pi)=(B^L\cup_{\varphi\vert_{E^L}} B^R,\pi=\pi^L\cup_{\varphi\vert_{\partial M^L\setminus E^L}} p(\pi^R)).\]

\end{theorem}
Obviously, cutting is the inverse operation for gluing. 
\begin{remark}
The Gluing Theorem easily generalises to the case when  the resulting manifold has a partial or foliated open book, rather than an honest open book. However, the statement does not generalize to self-gluing, since the required  reparametrization  cannot be done globally. Self-gluing is therefore only possible if $\varphi$ maps  one foliation to the other without reparametrization.
\end{remark}

\subsection{Examples}\label{ssec:cutex}
\begin{figure}[h!]
\begin{center}
\includegraphics[scale=0.8]{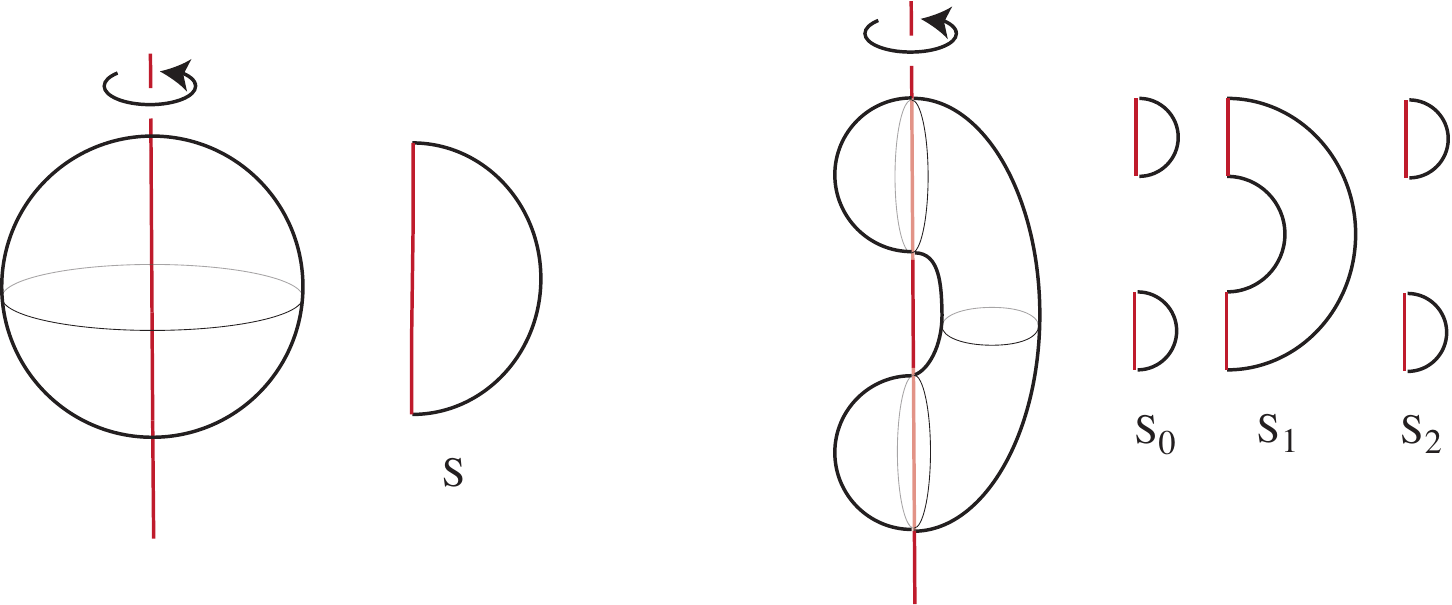}
\caption{ Left: A ball intersecting the binding once.  Right: A ball intersecting the binding twice. }\label{fig:exball}
\end{center}
\end{figure}

The cutting theorem allows us to construct many examples for foliated open books. The next two examples are Darboux balls embedded in the standard contact structure with the angular open book:
\begin{example}\label{ex:ball}
The first picture in Figure~\ref{fig:exball}  shows a foliated open book for the ball with a unique homeomorphism-type of page.  The complement of $B^3$ in the open book with disk-like pages for $S^3$ is a diffeomorphic foliated open book.
\end{example}

\begin{example}\label{ex:tel}
The right-hand picture in Figure~\ref{fig:exball} shows the ball from Example~\ref{ex:ball} after a finger move is performed. The associated foliated open book $(\{S^{\textit{tel}}_{i}\}_{i=0}^2,\id)$ has two distinct homeomorphism-types of pages, while the map is still the identity. The complement of $(\{S^{\textit{tel}}_{i}\}_{i=0}^2,\id)$ in $(D^2, \id)$ is again conjugate to $(\{S^{\textit{tel}}_{i}\}_{i=0}^2,\id)$.
\end{example}

Notice that the pages in the examples above do not define partial open books.  In the first case, there is a single homeomorphism-type of page, and in the second, the ``big" page cannot be built up from the complement of the ``small" page by adding one-handles.  

\begin{example}\label{ex:stab}
In this example, we identify two embedded foliated open books as  submanifolds of  the  open book $(B_0,\pi_0)$ associated to  the Hopf fibration on $S^3$. For $S^3=\{|z|^2+|w|^2=1\}\subset \mathbb{C}^2$, the binding $B_0$ is the set $\{zw=0\}=\{z=0\}\cup\{w=0\}$, oriented  by $\partial w$ and $\partial z$ on the respective components, and  $\pi_0 (w,z)=zw/|zw|$. The pages of this open book are annuli, and $(B_0, \pi_0)$ supports the standard tight contact structure on $S^3$.  

Consider the arc $\gamma_0=\{\Im z=0, zw/|zw|=1\}$ on the page $S_1=\pi_0^{-1}(1)$. Let $N(\gamma_0)$ denote a three-dimensional neighborhood of $\gamma_0$ that intersects pages $\pi^{-1}(t)$  in a 
rectangle  for $t\in(1-\epsilon, 1+\epsilon)$ 
and in a pair of half discs containing $\partial \gamma_0$ for $t\in(1+\epsilon, 1-\epsilon)$.    The open book foliation has a pair of hyperbolic singularities on $\pi^{-1}(1\pm \epsilon)$ as shown in Figure~\ref{fig:exball2}.

 \begin{figure}[h!]
\begin{center}
\includegraphics[scale=0.7]{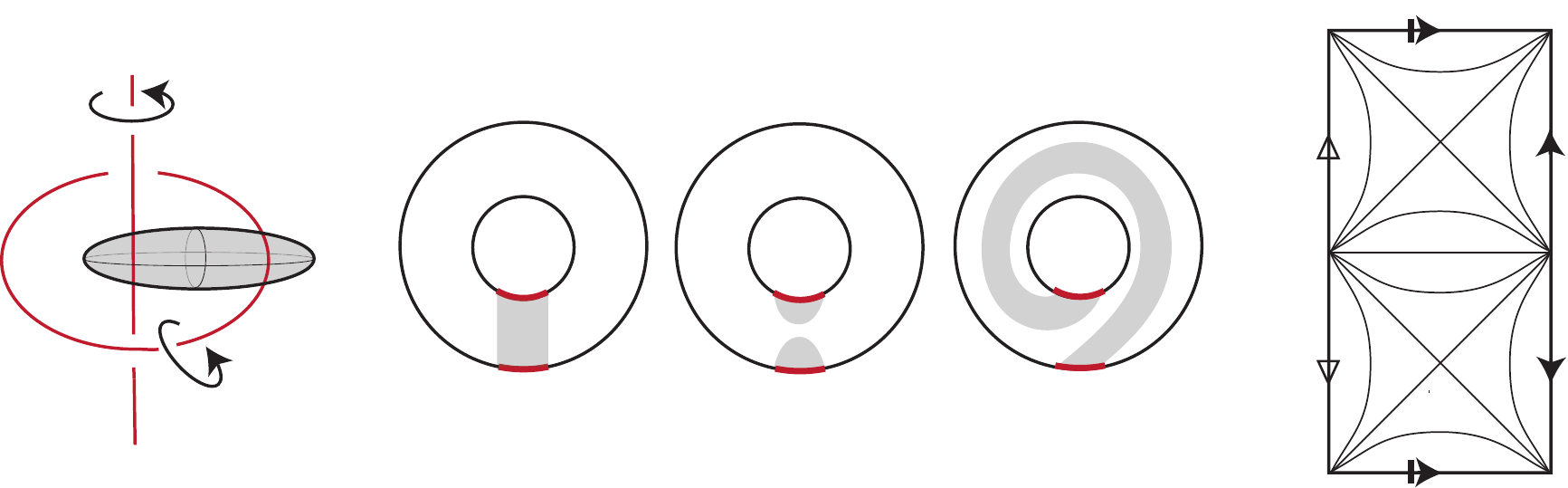}
\caption{ Cutting along the sphere shown on the left separates $S^3$ into a pair of distinct foliated open books.  The center figure shows the pages of the associated abstract foliated open books, and the open book foliation on the cutting sphere is shown on the right. }\label{fig:exball2}
\end{center}
\end{figure}

The foliated open book for $N(\gamma_0)$ may be identified with the foliated open book  $(\{S^{\textit{tel}}_{i}\}_{i=0}^2,\id)$ from Example~\ref{ex:tel}, but the foliated open book for $S^3\setminus N(\gamma_0)$ is distinct. This open book will be the key to defining stabilization in Sections~\ref{ssec:estab} and \ref{ssec:absstab}.

\end{example}

\subsection{Stabilization of embedded foliated open books}\label{ssec:estab}

Naturally, we would  like to understand the relationship  between embedded foliated open books that support the same contact structure.  As is the case with other versions of open books, foliated open books admit an operation called \textit{stabilization}  that preserves the contactomorphism class of the supported contact structure.  Although stabilization in other contexts is often defined in terms of abstract open books, we present it here in the embedded setting as application of the gluing results above; the abstract version is discussed immediately afterward in Section~\ref{ssec:absstab}.

Let $(B_e,\pi_e, \F_{\pi_e})$ be an embedded foliated open book for $M$ and let $\gamma$ be an arc properly embedded in some page $\overline{\pi_e^{-1}(t)}$.  As in Example~\ref{ex:stab}, let $N(\gamma)$ denote a three-dimensional neighborhood of $\gamma$ that intersects nearby pages in a single disc and other pages in a pair of discs near $\partial \gamma$.  This choice implies that up to reparameterization, the open book foliation on $\partial N(\gamma)$ matches the open book foliation on $\partial N(\gamma_0)$.

 \begin{definition}\label{def:Mstab}  The \emph{(positive) stabilization} of  $(B_e,\pi_e)$ along  $\gamma$ is  the manifold formed by gluing the foliated open books $S^3 \setminus N(\gamma_0)$ and $M\setminus N(\gamma)$ along their boundary two-spheres.
 \end{definition}
  
This operation is a refinement of the connect sum that respects the open book structures on the two manifolds; it also defines positive stabilization for (honest) open books.  As every change happens in the interior of the contact 3-manifold,  the following is an immediate consequence:

\begin{proposition}
Any positive stabilization of $(B_e,\pi_e, \F_{\pit_e})$ supports the same contact structure as $(B_e,\pi_e,  \F_{\pit_e})$.
\end{proposition}

Note that as stabilization is defined as a strictly interior operation, the definition above may be applied verbatim to Morse foliated open books. Negative stabilization may be defined the same way by using the open book given induced by the function $z\overline{w}/|z\overline{w}|$ on $S^3$.

Moreover, in Section~\ref{sec:giroux} we will prove the following result.

\begin{theorem}\label{thm:giroux}[Giroux Correspondence for foliated open books]  Any pair of foliated open books supporting $(M, \xi, \F_\xi)$  are isotopic after a sequence of positive stabilizations.  
\end{theorem} 
The proof of this theorem will rely on the proof of the analogous statement for partial open books. \footnote{The argument implicitly uses the proof of Giroux correspondence for (honest) open books, and it does not give an independent proof for Giroux correspondence in the classical case.}

\subsection{Stabilizing abstract foliated open books}\label{ssec:absstab}
It will be convenient to stabilize abstract foliated open books without invoking the equivalence of Proposition~\ref{prop:mtoa}, so we conclude this section with a reformulation of stabilisation adapted to the abstract case.

\begin{definition}\label{def:astab} Given an abstract open book  $(\{S_i\},h)$, let $\gamma$ be a  properly embedded arc in $S_r$ whose endpoints $p,q$ lie on $B\subset \partial S_r$.
The \emph{stabilization along  $\gamma$} is defined as follows: 
\begin{enumerate}
\item first perform an  $r$-shift of $(\{S_i\},h)$ to $(\{S_i[r]\},h[r])$;
\item\label{item:handle}   define a new abstract foliated open book $(S_i'[r]\}_{i=0}^{2k},h'[r])$ by
 \begin{itemize}
 \item[-] $S_i[r]'=S_i[r]\cup H$, where $H$ is a 1--handle with attaching sphere $p\cup q$; and
 \item[-] $h'[r]=D_{\overline{\gamma}}\circ h[r]$, where $D_{\overline{\gamma}}$ is a right-handed Dehn twist along the circle formed by the $r$-shift of $\gamma$ and the core of $H$; here, $h[r]$ also denotes its extension to $H$ by the identity;
 \end{itemize}
 \item perform a $(-r)$-shift  to obtain $(\{S_i'\},h')$, where  $S_i'$ is still obtained from $S_i$ by a handle attachment along $p$ and $q$.
  \end{enumerate}
 \end{definition}
\begin{remark}
The shift is indeed necessary in the definition in order to obtain an object  invariant under conjugation, as  in some cases the arc $\gamma$ cannot be ``found'' on the original $S_0$. 
\end{remark}
\begin{proposition}[Equivalence of stabilizations] 
The various definitions of stabilization are consistent: 
\begin{enumerate}  
\item Any stabilization of a Morse foliated open book can be realized as a stabilisation of a corresponding abstract foliated open book;
\item Any stabilization of an abstract foliated open book induces a stabilization of a corresponding Morse foliated open book. 
\end{enumerate}
\end{proposition}

The result follows from Definitions \ref{def:astab} and \ref{def:Mstab}.

\begin{example}\label{ex:torusstab}
As an example, we will stabilize the solid torus seen in Examples~\ref{ex:torus} and \ref{ex:notsep}.  Figure~\ref{fig:torusstab} shows a sequence of pages which differ from the pages of Example~\ref{ex:torus} by the addition of a handle connecting the two components of the original binding. At each step, the bold curve is either the arc that is cut along to yield the next page or the co-core of the handle that was attached in the step from the previous page.  We record cutting arcs on all previous pages and cocores of added arcs on all subsequent pages, and with the opposite convention, we record the endpoints of these sorting arcs  by bold dots on the boundary. Note that as the arcs are disjoint,  the foliated open book is now sorted.

\begin{figure}[h!]
\begin{center}
\includegraphics[scale=0.9]{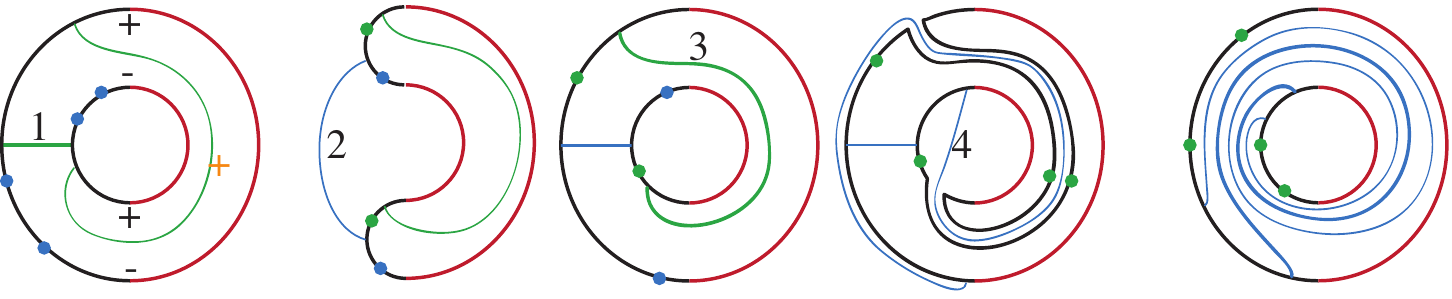}
\caption{ A stabilized abstract foliated open book for the  solid torus from Example~\ref{ex:torus}. }\label{fig:torusstab}
\end{center}
\end{figure}
The  correct identifications between successive pages may be realized by horizontal translation in the figure, and the left hand side is identified with the right by a right handed Dehn twist along the core of the annulus.

\end{example}

\begin{remark}\label{rmk:distribmon}

There is an alternative way to think about stabilizations that occur away from $S_0$ if we distribute the monodromy of an abstract foliated open book throughout time instead of concentrating it at $t=0=1$. In this setting, we have a sequence of pages $\{P_i,X_i\}_{i=0}^{2k-1}$ that are again cornered surfaces with a fixed boundary part $B$. We go from $P_{i-1}$ to $X_i$ by the usual cutting or gluing, but then we ``twist'' the page by a diffeomorphism $h_i$ that fixes $B$ to get the new page $P_i$. For $0\le i \le 2k-1$, cyclicly we have: 
${P}_{i-1} \xrightarrow[p_{i-1},q_{i-1}\text{ or }\gamma_ {i-1}^+]{(\text{\textbf{add} or \textbf{cut}})}{X}_{i}\xrightarrow[\cong]{h_{i}}{P}_{i}$.

When this identification $h_i$ is the identity for all $i\neq 0$,  we recover the old definition by setting $S_i=P_i$ for $0\le  i\le 2k+1$, $S_{2k}=X_0$ and $h=h_{2k}$.

As before, each relation $P_{i-1} \xrightarrow[p_{i-1},q_{i-1}\text{ or }\gamma_ {i-1}^+]{(\text{\textbf{add} or \textbf{cut}})}{X}_{i}$ defines a cornered handlebody $H_i$, but now these handlebodies are glued using the diffeomorphisms $h_i$.

This flexibility provides a more natural description of stabilization along a curve $\gamma\subset P_r$ whose endpoints lie on $B$.  The new pages are 
obtained by attaching a 1-handle to each original page along the endpoints of $\gamma$, so $P_i'=P_i\cup H$ and $X_i'=X_i\cup H$.  As before, denote by $h_i$ the  extension of the diffeomorphism $h_i$ to $H$ by the identity. Now let $h_i'=h_i$ for $i\neq r$ and $h_r'=D_{\overline{\gamma}}\circ h_r$.

\end{remark}

\begin{example}Figure~\ref{fig:torusdistr} shows the abstract open book from Example~\ref{ex:torusstab} with the distributed monodromy described in Remark~\ref{rmk:distribmon}.  As above, the sorting arcs are shown in blue and green, and at most steps, $h_i$ is the natural embedding
 induced by the page.  However,  $h_2:X_2\rightarrow P_2$ is a positive Dehn twist.  One may easily see that this is equivalent to shifting the original abstract open book indices by $2$. 

\begin{figure}[h!]
\begin{center}
\includegraphics[scale=0.9]{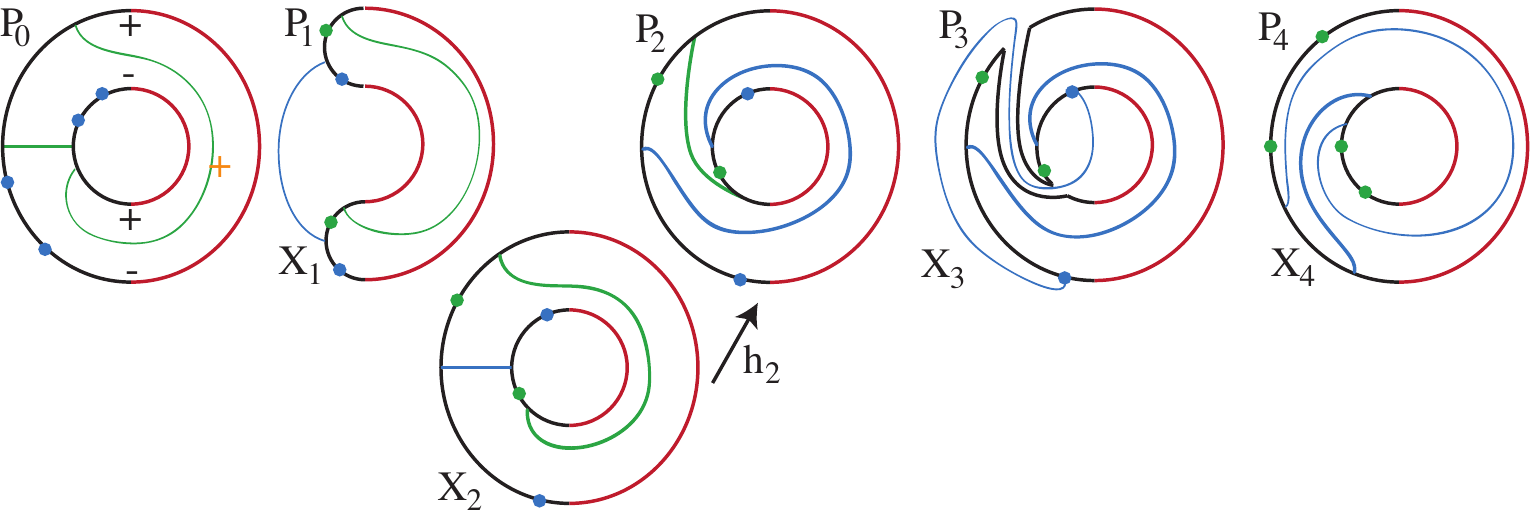}
\caption{ Example~\ref{ex:torus}. }\label{fig:torusdistr}
\end{center}
\end{figure}

We employ this perspective sparingly, as it requires keeping track of more data at each step, but it provides a useful alternative to the formal shifting of Definition~\ref{def:astab}.
\end{example}



\section{From Foliated open books to contact 3-manifolds}\label{sec:existence} 

In this section we prove the fundamental existence and uniqueness results that ensure foliated open books are useful tools for studying contact manifolds.  
\subsection{Existence of supported contact structure}\label{sec:existencexi} 

Having established the connection between abstract foliated open books and embedded foliated open books, we are ready prove the existence of a supported contact structure. This is  easiest in the abstract setting.
\begin{theorem}\label{thm:existencecontact} Any abstract foliated open book supports a contact structure.
\end{theorem}
\begin{proof} 
Let $(\{S_i\},h)$ be an abstract foliated open book.
The basic idea of the proof is to  find an (honest) abstract open book $(S',h')$ for a closed 3--manifold so that $(\{S_i\},h)$ embeds into it pagewise and $h$ is the restriction of $h'$ to $S_{2k}$. Then we construct the contact 3--manifold corresponding to $(S',h')$ and use the cutting result (Theorem \ref{thm:cut}) to get a contact structure supported by $(\{S_i\},h)$.

The condition that $\F_{\pit}$ does not contain circles translates to each $\alpha_i$ being a union of intervals.  
We will find a surface $S'$ and embeddings $\iota_i\colon S_i\hookrightarrow S'$ which satisfy the following properties: 
\begin{itemize}
\item[-] for all $i$,  $\iota_i$ embeds $B$ in a fixed  $B'\subset \partial S'$;
\item[-]   for all $i$, $\iota_i$ properly embeds $(\alpha_i, \partial \alpha_i)$  into $(S', \partial B')$;
\item[-]  if $S_{i-1}\xrightarrow[\gamma_{i-1}]{\textbf{cut}}S_i$, then $\iota_{i-1}(S_{i-1})\setminus  \iota_{i}(S_{i})$ is diffeomorphic to a neighbourhood of $\iota_{i-1}(\gamma_{i-1})$
with a pair of opposite edges in $\iota_{i-1}(\alpha_i)$ and the other pair of edges in $\iota_{i}(\alpha_{i})$; 
\item[-] if $S_{i-1}\xrightarrow[p_{i-1}, q_{i-1}]{\textbf{add}}S_i$, then $\iota_{i}(S_{i})\setminus  \iota_{i-1}(S_{i-1})$ is diffeomorphic to a  
rectangle with a pair of opposite edges in $\iota_{i-1}(\alpha_{i-1})$ and the other pair of edges  in $\iota_{i}(\alpha_{i})$.
\end{itemize}

We construct $S'$ inductively as follows. Start with any $S'_0$ such that $S_0$ embeds into it and satisfies the first two conditions for $i=0$. Now assume that we have defined $\iota_j$ and $S_j'$ for $j\le i$ so that all the conditions are satisfied for $j\le i$.  If $S_{i}$ is related to $S_{i-1}$ by cutting along the arc $\gamma_{i-1}^+$, then we can set $S_{i}'=S_{i-1}'$ and keep all the embeddings; then  define $\iota_i=\iota_{i-1}\vert_{S_i}$. If, on the other hand, $S_{i}$ is related to $S_{i-1}$ by a 1-handle addition, then we would like to realize the 1-handle addition inside $S'_{i}$. This means that we need a path $\gamma_{i-1}'$ in $S_{i-1}'\setminus \iota_{i-1}(S_{i-1})$ connecting the attaching points of the handle. If there is such a path, then we keep $S_{i}'=S_{i-1}'$ and let the image of $\iota_{i}$ be the union of the image of $\iota_{i-1}$ and a tubular neighborhood of $\gamma_{i-1}'$. If there is no such path, then we add a handle to $S_{i-1}'$ to connect the two regions of $S_{i-1}'\setminus \iota_{i-1}(S_{i-1})$ containing the attaching points.  Here we use the fact  that none of the components of $S_{i-1}'\setminus \iota_{i-1}(S_{i-1})$ is disjoint from the boundary, which follows from the fact that $\alpha_{i-1}$ has no circles. Call the new surface $S_{i}'$.   To complete the inductive step, compose all embeddings with the obvious inclusion $\iota\colon S_{i-1}'\hookrightarrow S_{i}'$ and let $\iota_{i}(S_{i})$ extend $\iota_{i-1}(S_{i-1})$ by an embedding which maps the additional handle to a neighborhood of $\gamma_i'$.

 Consider the embeddings $\iota_0\colon S_0\hookrightarrow S_n'$ and $\iota_n\colon S_n\hookrightarrow S_n'$.  The compositions $\iota_0^{-1}\circ h\circ \iota_n\colon \iota_n(S_n)\to \iota_0(S_0)$ might not be extendable to a diffeomorphism of $S_n'$, but by embedding $S_n'$  into a yet larger  surface $S'$, we may assume that the composition extends to the desired diffeomorphism $h':S'\rightarrow S'$. 

The embeddings $\iota_i$ determine a submanifold $M$ corresponding to the abstract foliated open book $(\{S_i\},h)$ in the manifold $M'$ corresponding to the open book $(S',h')$.   By Theorem \ref{thm:cut}  there is a contact structure $\xi$ compatible with $(S',h')$ that restricts to $M$ as a contact structure compatible with 
$(\{S_i\},h)$, as needed.  
\end{proof}

\subsection{Uniqueness of supported contact structure}\label{sec:uniquexi} 

In this subsection we show that an abstract foliated open book supports a unique contact structure.

\begin{theorem}\label{thm:uniquecontact} Any two contact structures supported by a foliated open book are isotopic. 
\end{theorem}
As a first step we will argue  that the handlebody $H_i$ has a unique tight contact structure compatible with a prescribed class of characteristic foliations on its boundary, up to contactomorphism fixing the boundaries.   Second, we show  that any contact structure compatible with the foliated open book restricts on each $H_i$ to a tight structure. The structure of this argument follows \cite{Torisu}, and together, these statements imply the desired uniqueness.

To begin, we describe the decoration induced on $\partial H_i$ by a compatible contact structure $\xi$ on $M$. 

Recall that the manifold obtained by gluing  the handlebodies $H_i$ along $S_i$, denoted $M'$, can be obtained from $M$ by blowing up the transverse curve $B$ to $B\times S^1$. This means that a compatible contact structure $\xi$ on $M$ induces a contact structure $\xi'$ on $M'$ that restricts to each $H_i$ as $\xi'_i$.  

The boundary of $H_i$ is separated into two pages, $S_{i-1}$ and $S_i$, together with the vertical component $V_i= (B\times I)\cup W_i$. Here, $W_i$ consists of rectangles $\alpha_i \times I$  coming from  the unchanged intervals $\alpha_i$, together with a ``saddle'' component corresponding to the two components of $\alpha_i$ that change. 
The characteristic foliation on $B\times S^1\subset \partial M'$ is given by the parallel circles $\{b\}\times S^1$ for $b\in B$, and thus by $\{b\}\times I$ on the $B\times I$ components of $\partial H_i$. 

The open book foliation on $W_i$ is given by the level sets of $\pit_i$. By assumption, the characteristic  foliation on $W_i$ is strongly topologically conjugate  to the foliation given by the level sets of 
$\pit_i$. This means that there are no saddle-saddle connections and no circles in $\F_{\xi'_i}$, so $W_i$ is a convex surface whose dividing curve $\Gamma$ may be  obtained as the boundary of the neighbourhood of the positive graph $G_{\scriptscriptstyle{++}}$ for either $\F_{\xi_i}$ or $\F_{\pit_i}$. 

Note, too, that since $\ker \alpha$ is supported by the open book, the Reeb vector field $R_\alpha$ is transverse to the pages. As $R_\alpha$ is a contact vector field, 
the pages (in particular, the $S_i$'s) are convex with empty dividing sets. 

This discussion establishes the first part of the following proposition:

\begin{proposition}\label{prop:unique} A contact structure $\xi$ supported by the abstract foliated open book $(\{S_i\}, h)$ induces on $H_i$ a tight contact structure $\xi_i$ with the decorations described  above. 

Furthermore, up to isotopy fixing $\partial H_i$, there is a unique tight contact structure on $H_i$ with these boundary decorations.
\end{proposition} 
 
\begin{proof} 
The tightness of $\xi_i$ follows from gluing (Theorem \ref{prop:glue})  and Torisu's Theorem~1.1 in \cite{Torisu} as follows. Suppose that $\xi'$ is a contact structure on $M$ compatible with $(\{S_i\},h)$. Then  as in the proof of Theorem \ref{thm:existencecontact}, we can construct an abstract foliated open book $(\{S_i''\},h'')$ compatible with $\xi''$ so that the two foliated open books glue together into a closed contact manifold $(M',\xi')$ compatible with the abstract open book $(S',h')$. Then by Torisu's statement we know that $\xi$ restricted to any submanifold $S'\times[t,t']\slash\sim\subsetneq M'$ is tight. Since $H_i$ is a subset of the blowup at $B$ of $S'\times[t_{i-1},t_i]$, it must also be tight. 

For uniqueness, we will smooth the boundary of  $ H_i$ (in the interior of $H_i$ and arbitrarily close to $\partial H_i$) in order to obtain a handlebody $H_i'$ with convex boundary and dividing curve $\Gamma_i$ which is determined by the decoration on $\partial H_i$. Since the contact structure on a collar neighbourhood of $\partial H_i$ is determined by the characteristic foliation on $\partial H_i$,  it suffices to show that there is a unique tight contact structure compatible with the dividing curve $\Gamma_i$ on $\partial H_i'$. 

In the above description, all components of $\partial H_i$ were convex except $B\times I$. 
In a neighborhood of each component of $B\times I$, we connect  $-S_{i-1}$ to $S_i$ via a nearby convex surface $F_i$ that has a 1-component dividing curve which is parallel to $B$.  See the third picture of Figure \ref{fig:boundary1}. 
  
 \begin{figure}[h!]
\begin{center}
\includegraphics[width=\textwidth]{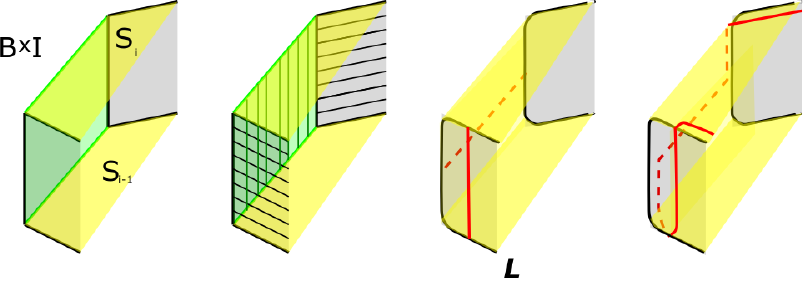}
\caption{Smoothing corners in $H_i$. The first figure shows a part of $\partial H_i$ near $B\times I$. The second figure shows the characteristic foliations, the third one is the picture after merging $S_{i-1}$ with $S_i$ and smoothing the corners, with the dividing curves indicated with red. The fourth picture depicts the dividing curve after completing the smoothing.}\label{fig:boundary1}
\end{center}
\end{figure}

The above process allows us to smooth out the codimension 2 corners as well,  leaving a 3-manifold $H_i''$ whose boundary consists of two convex surfaces, $F_i$ and a smoothed $W_i$, that intersect in a Legendrian curve $L$. See the third picture of  Figure~\ref{fig:boundary1}.
   
Recall that $W_i$ is convex and that the dividing curve $\Gamma\subset W_i$ is  the boundary of a neighborhood of the stable separatrices of the positive hyperbolic points. 
 In the blown-up picture of the boundary, the positive elliptic points correspond to the intervals $E_+\times I$; as usual, $\partial B=-E_-\cup E_+$.  Each product component of $W_i$ contains a vertical component of $\Gamma$ of the form $\{p\}\times I$, where $p$ is a point on $\alpha_{i-1}$ close to $E_+$. The saddle component contains the same vertical dividing curve if $S_{i}$ is obtained from $S_{i-1}$ by a 1-handle addition (as in this case the corresponding hyperbolic point is negative), but the dividing curve is as shown in Figure \ref{fig:unique2} if $S_i$ is obtained from $S_{i-1}$ by cutting along an arc $\gamma_{i-1}$.

\begin{figure}[h!]
\begin{center}
\includegraphics[scale=0.6]{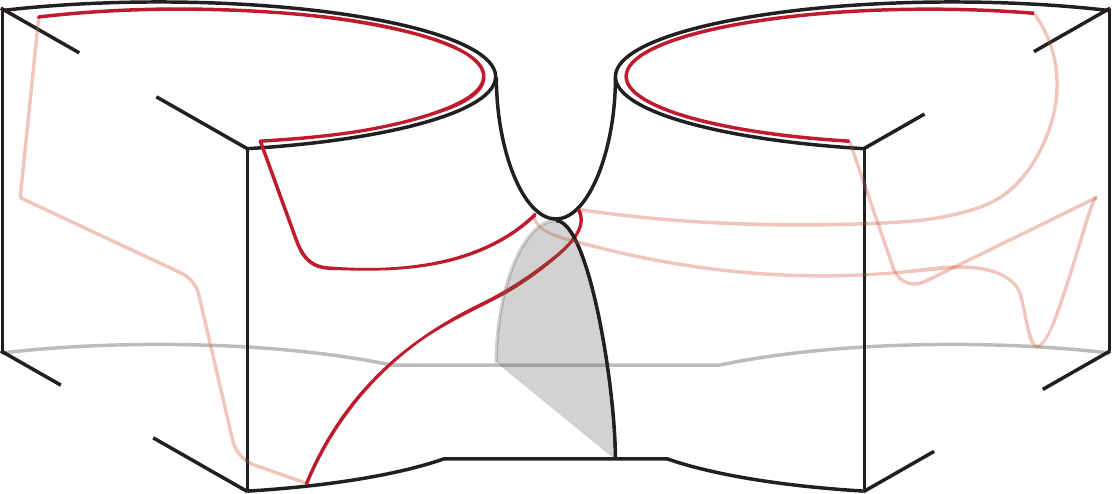}
\caption{ The smoothed $\Gamma$ intersects a decomposing disc twice on a cornered handle body associated to a positive hyperbolic point }\label{fig:unique2}
\end{center}
\end{figure}

In summary, the Legendrian curve $L$ intersects $\Gamma$ on $F_i$ once at every component corresponding to $E\times I$, while on $S_i$ it intersects $\Gamma$ near the points $E_+\times \{0,1\}$ on $\alpha_{i-1}$ and $\alpha_i$, respectively. See, again, the third picture in Figure  \ref{fig:boundary1}.

Next we can apply standard smoothing along the Legendrian curve $L$ to obtain the handlebody $H_i'$ with dividing curve $\Gamma_i$ as seen in Figure~\ref{fig:unique2}.

Finally, we will prove that $H_i'$ has a unique contact structure by describing a disc decomposition of $H_i'$ such that the boundary of each disc $D$ intersects $\Gamma_i$ in two points. As a first case, assume that $S_i$ is obtained from $S_{i-1}$ by a 1-handle addition. Then choose arcs $\{a_j\}$ properly embedded into $S_{i-1}$ with endpoints on $\alpha_{i-1}$ that cut up each component of $S_{i-1}$ into discs, and let $D_j$ be the (smoothed) $a_j\times I$. After the smoothing on the fourth picture of Figure~\ref{fig:boundary1}, it is clear that the boundary of each $D_j$ intersects $\Gamma_i$ in  two points. 

This suffices to cut the tight $H_i'$ into tight balls, proving the proposition in the case of handle addition.  If $S_i$ is obtained from $S_{i-1}$ by a cutting along a curve $\gamma_{i-1}$, then we can first choose arcs $a_j$ on $S_i$ that cut up $S_i$ into discs and again take $D_j=a_j\times I$ to complete the proof.

\end{proof}
This also finishes the proof of Theorem \ref{thm:uniquecontact}.

In the proof of Proposition \ref{prop:pobfobequiv}, we will require the extension of this approach to a concatenation of multiple $H_i$ handlebodies. Recall Definition~\ref{def:sortedcobord}.

\begin{proposition}\label{prop:sorthandlebody} Any sorted submanifold $_tH_{t'}$ supports a unique tight contact structure compatible with the boundary decorations described above. 
\end{proposition}

The argument is a minor extension of the analysis above.

\begin{proof} 
Recall from Lemma~\ref{lem:ppersists} that any pair of pages in the sorted $_tH_{t'}$, cut along their respective sets of sorting arcs, are diffeomorphic after smoothing; we view this surface as a minimal page that persists for all $t$.  

As seen in Figures \ref{fig:boundary1} and \ref{fig:unique2}, each  critical submanifold $W^s(h_+)$ and $W^u(h_-)$ intersects the dividing curve $\Gamma$ associated to a smoothed $_tH_{t'}$ twice.  After decomposing along these discs, the result is a product cobordism which may be further decomposed along discs guided by arcs $\{a_i\}$ which cut the minimal page into discs, as above.
  \end{proof}


\section{Relationship to Partial Open Books}\label{sec:pob}\label{sec:ptofob}
Partial open books offer a well established tool for studying contact manifold with convex boundary, and in this section, we explore the connections between partial open books and foliated open books.  
Section~\ref{sec:buildup} describes how to build an abstract foliated open book from the data defining an abstract partial open book, while Section~\ref{sec:pobinfob}  reverses the process and constructs an abstract partial open book from an abstract foliated open book.  In each  case, some conditions on the initial open book are imposed, and we show in Section~\ref{sec:disj} that these can be achieved via positive stabilization.  These constructions are not inverses, as the construction of a foliated open book from a partial one allows some choice,  but Section~\ref{sec:fobpobequiv} shows that these transformations preserve the supported contact structure.    Section~\ref{sec:existfob} applies the relationship between partial and foliated open books to prove another existence result: if $(M, \xi)$ is a contact manifold with appropriate characteristic foliation, then there is a compatible foliated open book (Theorem~\ref{thm:fobexistence}).  Finally, we prove  a Giroux Correspondence for foliated open books in Section~\ref{sec:giroux}. 

\subsection{Foliated open books from partial open books}{\label{sec:buildup}
The intuition behind transforming a partial open book into a foliated open book is straightforward; in a partial open book, a large change in Euler characteristic is concentrated at two pages where the $S$ and $P$ handlebodies meet.  For appropriate $S$ and $P$, this may be distributed as 
 a sequence of small changes across many pages  to yield a foliated open book.

We begin with a topological condition describing when a surface may be built up from a subsurface by successive attachment of one-handles.  Note that this relationship appears in the definition of an abstract partial open book \cite{EO}.  

\begin{lemma}\label{lem:builtup} Suppose that $S$ is a surface with non-empty boundary and $X$ is a nonempty subsurface with cornered boundary such that each component of $\partial X$ is either contained in $\partial S$ or is polygonal with alternating edges in $ \partial X\cap \partial S$ and $\partial X\setminus \partial S$. Then $S$ can be built up from $X$ by successive attachments of 1-handles if and only if the boundary of any component of $\overline{S\setminus X}$ intersects $\partial X$ in at least two intervals.
\end{lemma}

\begin{proof}

For the ``if" direction, note first that each handle that creates a new component of $S\setminus X$ must attach directly to $X$, so the number of interval components of $\partial S \cap \partial X$ on the new component (i.e., the added handle) is two.  No subsequent one-handle attachments can decrease the number of components of $\partial X \cap \partial S$ on the boundary of any connected component of $S\setminus X$.

For the ``only if" direction, we describe the co-cores of the handles which build $S$ up from $X$. Fix a connected component $C$ of $\overline{S\setminus X}$.  We can simplify the topology of $C$ by cutting along arcs with boundary on $\partial S\setminus \partial X$.  Each cut turning $C$ to $C'$ is along an arc which is the co-core of a one-handle which attaches to $C'$ to yield $C$.  

Suppose that successive cutting yields a region $C''$ that is polygonal with $2n$ edges, $n>2$. Then we can add an additional $n-1>1$ cutting arcs with endpoints on $\partial C''\setminus \partial X$, parallel to all but one of the components of $\partial S\cap\partial X$ on $\partial C''$. This yields a collection of bigons, each with one edge in $\partial S\cap \partial X$ and the other  disjoint from $X$.  Taken collectively in the original surface, the cutting arcs are the co-cores of a set of handles that build  $C$ up from $X$.  
\end{proof}

If $S$ and $X$ are as in the above lemma, we say that $S$ \emph{can be built up from $X$}. With this phrasing, we recall that the definition of an abstract partial open book $(S,P,h)$ requires that $S$ can be built up from $S\setminus P$. To create foliated open books, however, we need that $S$ can be built up both from $P$ and $h(P)$. This property can be achieved by positive stabilization:

\begin{lemma}\label{lem:pobtofob} Any partial open book $(S, P, h)$ may be positively stabilized to some $(S'', P'', h'')$ with the property that $S''$ can be built up from $P''$ and from $h(P'')$.
\end{lemma}

\begin{proof}     
In a partial open book stabilization $(S,P,h)\rightarrow (S',P',h')$, the added one-handle  becomes part of $P'$, so choosing the attaching sphere to have at least one component on $\partial S\setminus \partial P$ implies that  $|\partial P' \cap (\partial S'\setminus \partial P')|>|\partial P \cap (\partial S\setminus \partial P)|$. Iterating, we  eventually stabilize to a partial open book $(S'',P'', h'')$ with the property that each component of $\partial S''\setminus \partial P''$ meets $\partial P''$ in at least two intervals.  Thus $S''$ can be built up from $P''$, as desired. One may similarly stabilize to ensure that $S''$ is also built up from $h(P'')$. 
\end{proof}

Suppose  now that $\mathcal{P}=(S,P,h)$ is a partial open book such that $S$ may be built up from  $P$ and $h(P)$.  This is equivalent to the statement that $S$ may be cut either along a set of properly embedded arcs $\{\delta_i^+\}_{i=1}^k$ to yield $P$ or along a set of arcs $\{\delta_j^-\}_{j=1}^k$ to yield $h(P)$. These arcs are the co-cores of the  handles added in the building-up construction. 

Set $S=S^{\mathcal{P}}_k$ and define $S^{\mathcal{P}}_{k+j}=S^{\mathcal{P}}_{k+j-1}\setminus \delta_j^+$ for $j\in \{1, \dots, k\}$.  Similarly, define $S^{\mathcal{P}}_{k-j}=S^{\mathcal{P}}_{k-j+1}\setminus \delta_j^-$ for $j\in \{1, \dots k\}$.  This yields a collection of abstract pages $\{S_i\}_{i=0}^{2k}$ that are related to each other by gluing  (for $i\le k$) and cutting (for $k\le i$). By construction, $S_{2k}^\mathcal{P}=P$ and $S_0=h(P)$, thus we can define $h^\mathcal{P}$ to be $h\vert_P\colon S_{2k}^{\mathcal{P}}=P\to h(P)=S_0\subset S$.  It is clear that this data defines a foliated open book, but the indexing of the cutting arcs is arbitrary, and different choices will yield distinct foliated open books.  We would like to ensure that this process yields a sorted foliated open book.  As a first step, we apply a $k$-shift to so that the ``big" page is now $S_0$.

Interpreting an indexed set of cutting arcs $\{\delta_i^\pm\}$ as the intersections of stable critical submanifolds with $S_k$, we must verify the conditions of Definition~\ref{def:pref}.   If the intersections are not properly ordered, then reindexing the arcs or performing handle slides may correct this.  Alternatively, one may positively stabilize the original partial open book so that no interval of $\alpha_k=\partial S_k\setminus B$ contains an endpoint of more than one $\delta$ arc.  The ordering conditions are then trivially satisfied, so the resulting foliated open book is sorted with respect to a preferred gradient-like vector field.

\begin{definition}\label{def:inducedabstract} A partial open book $(S,P,h)$ is \emph{sufficiently stabilized} if $S$ may be built up from $P$ and from $h(P)$, and furthermore, if there is a choice of associated cutting arcs that yields an abstract foliated book $\A(\mathcal{P})=(\{S^{\mathcal{P}}_i\}, h^{\mathcal{P}})$ which is sorted with respect to a preferred gradient-like vector field.  For any such choice of cutting arcs $\{\delta_i^\pm\}$, we say that $\A(\mathcal{P})=(\{S^{\mathcal{P}}_i\}, h^{\mathcal{P}})$ is an abstract foliatied open book \emph{induced by} the sufficiently stabilized partial open book $\mathcal{P}=(S,P,h)$.
\end{definition}

We will show in Section~\ref{sec:fobpobequiv} that this induced foliated open book  is associated to the ``same'' contact 3-manifold with boundary.

\subsection{Partial open books from foliated open books}\label{sec:pobinfob}
Transforming a foliated open book into a partial open book requires more effort  than the reverse process, and it relies on the Morse model of a foliated open book to define a ``minimal" page that can play the role of $P$ and a ``maximal" page that can play the role of $S$. As in the previous construction, some carefully chosen preliminary stabilizations may be required before this is possible. This section relies heavily on the discussion of sorted foliated open books from Section~\ref{sec:sort}.
 
We will prove the following result in Section~\ref{sec:disj}:
\begin{proposition}\label{prop:disj}
Any foliated open book $\A=(\{S_i\},h)$ may be positively stabilized to a sorted foliated open book. 
\end{proposition}

As an illustration, consider Examples~\ref{ex:notsep} and \ref{ex:torusstab}, which show how the (unsorted) foliated open book introduced in Example~\ref{ex:torus} may be stabilized to a sorted foliated open book.
 
Suppose, assuming Proposition~\ref{prop:disj}, that $\A=(\{S_i\}, h, \{\gamma_k^\pm \})$ is a sorted  foliated open book for $(M, \xi, \F_{\pit})$. 
Set $S^\A=S_0$ and let $P^\A=\overline{S_0\setminus R_0}$, where $R_0$ is the subsurface of $S_0$ from Definition~\ref{def:r}. See Figure \ref{fig:ttopob}.
Lemma~\ref{lem:ppersists} implies that $P^\A$ embeds in each $S_t$ as subsurface. Define $\iota$ to be the map which embeds $P^{\A}$ into $S_{2k}$, and set ${h^\A}=h\circ \iota$. 
 
 \begin{definition} Given a sorted foliated open book $\A=(\{S_i\}, h, \{\gamma_k^\pm\})$, the \textit{associated triple}  is $\mathcal{P}(\A)=({S^\A}, {P^\A}, {h^\A})$.
  \end{definition}

  \begin{proposition}\label{prop:fobtopob}   Given a sorted abstract foliated open book $\A=(\{S_i\}_{i=1}^{2k}, h, \{\gamma_k^\pm\})$ with at least one of $k>0$ or $|\partial M|>1$, there exists a positive stabilization $\A'=(\{S_i'\}_{i=1}^{2k}, h')$ with the property that the associated triple  $\mathcal{P}(\A')=({S_{\A'}}, {P_{\A'}}, {h_{\A'}})$ is a partial open book.    \end{proposition}
 We call a foliated open book as above \emph{sufficiently stabilized}, and the obtained partial open book $\mathcal{P}(\A)$ an \emph{associated partial open book}.
  
 Proposition~\ref{prop:fobtopob} is a partner to Proposition~\ref{lem:pobtofob}, in the sense that these assert the existence of procedures for turning partial open books into foliated open books and vice versa.  In fact, Proposition~\ref{prop:pobfobequiv} will show that these are the correct procedures, in that they preserve the supported contact structure.  The proof of Proposition~\ref{prop:fobtopob} rests on a pair of lemmas (\ref{lem:rplus}, \ref{lem:nocircle}) which will be proven in the next subsection.

 \begin{remark} 
  We note that the hypotheses $k>0$ or $|\partial M|>1$ in Proposition~\ref{prop:fobtopob} eliminate only the case of a foliated open book whose unique boundary component is a sphere with the ``trivial" foliation whose only singularities are a pair of elliptic points.  Gluing a trivial foliated open book for a  Darboux ball (Example~\ref{ex:ball}) to $M$ yields a closed manifold with an open book decomposition. Much of the machinery in this section is building towards a proof of a Giroux Correspondence for foliated open books, but we note that the classical Giroux Correspondence applies in this case.  We thus ignore the  $k=0$,$ |\partial M|=1$ case in what follows.
  \end{remark}

\subsection{Sorting via stabilization}\label{sec:disj}
 
Propositions~\ref{prop:disj} and \ref{prop:fobtopob} each assert that an abstract foliated open book may be stabilized to achieve some given property.  In order to prove these claims, it will be useful to identify a subsurface of the foliated boundary which is isotopic to $R_+=N(G_{\scriptscriptstyle{++}})$.

We assume a preferred gradient-like vector field $\nabla\pit$ on $\partial M$  and its extension $\nabla\pi$ on $M$. Recall that $w^s(h)$ and $w^u(h)$ are the stable and unstable flows of $\nabla\pit$ corresponding to the critical point $h$. Fix a regular time $t=0$, and set $w^s(h)'$ to be the subinterval of $w^s(h)$ that contains $h$ and intersects $\pit^{-1}(0)$ exactly at its two endpoints. Set \[R_+'=\pit^{-1}[0,\eta]\cup \bigcup_{h_+}N(w^s(h_+)'),\] where $\eta$ is sufficiently small to ensure  the leaves of $\pit^{-1}[0,\eta]$ are all regular. Then we claim the following:
\begin{lemma}\label{lem:rplus}
$R_+\subset \partial M$ is isotopic to $R_+'\subset \partial M$ in $\partial M$.  If the foliated open book is sorted, then $R_+'$ is isotopic to ${S^\A} \setminus {P^\A}$ in $M$.
\end{lemma}

\begin{proof}
Each of these regions is a neighborhood of a graph whose vertices are the positive elliptic points, and each edge connects the same pair of positive elliptic points. 

\begin{figure}[h]
\begin{center}
\includegraphics[scale=1.5]{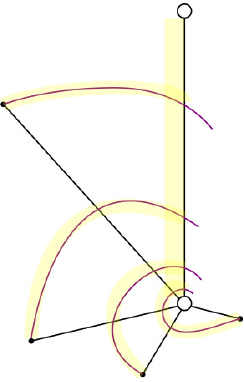}
\caption{$R_+'$ in the star of a positive elliptic point.}\label{fig:circord}
\end{center}\end{figure}

Because $\nabla \pit$ is preferred,  
 the stable separatrices $w^s(h_+)$ intersect each interval $I_{+}\subset \pit^{-1}(0)$ in the same the cyclic order as the separatrices of $\F_{\pit}$ which terminate at the positive elliptic point at the end of $I_+$ 
See Figure \ref{fig:circord}. Taking the union of the neighborhood of $I_+$ together with the separatrices recovers the original the cyclic order; this may be seen by retracting the $0$-leaf to the positive elliptic point and bringing the endpoints of the separatrices along.

Now suppose that the foliated open book is sorted, and recall the construction of the associated triple $({S^\A}, {P^\A}, {h^\A})$. For each positive hyperbolic point, the arc $w^s(h_i)'$ can be pushed down to $\gamma_i^+\subset S_0$ along the half-disc $W^s(h_i)$.  Extending these ``pushes'' to an isotopy with support ${}_0R_1$ gives  an isotopy of $R_+'$ onto $S^\A\setminus P^\A$, as desired.

\end{proof}

Having defined the triple associated to a sorted foliated open book, Lemma~\ref{lem:nocircle} is a first step towards showing that $({S^\A}, {P^\A}, {h^\A})$ meets the topological conditions to define an abstract partial open book.

 \begin{lemma}\label{lem:nocircle} Suppose that $({S^\A}, {P^\A}, {h^\A})$ is the  triple associated to some sorted foliated open book.  Then no component of $\partial {P^\A}\setminus \partial {S^\A}$ is a circle.
 \end{lemma}
In fact, Lemma~\ref{lem:nocircle} is a corollary of Lemma~\ref{lem:rplus}:
\begin{proof}[Proof of Lemma~\ref{lem:nocircle}]
Components of $\partial {P^\A} \setminus \partial {S^\A}$ are homeomorphic to components of $\partial R_+'\setminus \pit^{-1}(0)$.  When $R_+'$ deforms to $R_+$ on $\partial M$, the boundary of $R_+'$ maps to $\Gamma$ and $\pit^{-1}(0)$ maps to $\Gamma\cap \pit^{-1}(0)$.  Recall from Section~\ref{sec:egcircle} that $\pit|_\Gamma$ is a covering of $S^1$.  Thus every component of $\Gamma$ intersects $\pit^{-1}(0)$, and it follows that $\Gamma\setminus \pit^{-1}(0)$ can have no circle components.  
\end{proof}
Now we are ready to prove Proposition~\ref{prop:fobtopob}:
\begin{proof}[Proof of  Proposition~\ref{prop:fobtopob}] 

Suppose that $\A=(\{S_i\}, h)$ is a sorted foliated open book with $k>0$.  We will find a  positive stabilization $\A'=(\{S'_i\}, h')$ with the property that in the associated triple $({S_{\A'}}, P_{\A'}, h_{\A'})$, the surface $S_{\A'}$  can be built up from each  of $P_{\A'}$, $h_{\A'}(P_{\A'})$,  $S_{\A'}'\setminus P_{\A'}$, and ${S^\A} \setminus {h_{\A'}}({P_{\A'}})$ by successive attachment of one-handles.

As a first step, observe that any positive stabilization of the original foliated open book induces a positive stabilization of the associated triple, where the latter will be defined exactly as in the case of an abstract partial open book. Stabilizing the original foliated open book at $S_0$ adds a $1$-handle to $S_0$,  and hence, ${S^\A}$, with the attaching sphere lying completely on $B$.  The stable manifolds which cut out ${P^\A}$ from ${S^\A}$ remain disjoint from this handle, so the entire handle becomes part of ${P^\A}$, and the monodromy changes by the appropriate Dehn twist; this stabilizes the triple $({S^\A}, {P^\A}, {h^\A})$.

Note  that we may achieve only a restricted class of partial open book stabilizations this way.  In the case of an abstract partial open book, there are no restrictions on the location of the attaching sphere, but each time a component of the attaching sphere lies on $\partial S^\A\setminus \partial B$, a pair of new intersections between $B$ and $\partial M$ is created.  Although this preserves the dividing set, and hence, the partial open book up to the relevant notion of equivalence, this represents a fundamental change in the foliated boundary by introducing new pair of elliptic points.  In what follows, therefore, we consider only stabilizations of the triple along attaching spheres on $B=\partial {S^\A} \cap \partial {P^\A}$; any such stabilization may be achieved by a stabilization of the original foliated open book at $S_0$.

Returning to the statement of the lemma, we note that Lemma~\ref{lem:nocircle} implies $\partial {P^\A} \setminus \partial {S^\A}$ consists only of intervals.  
We  will show that there exists a sequence of stabilizations along $\partial {P^\A} \cap \partial {S^\A}$ that ensures the resulting $S_{\A'}$ may be built up from each of the specified subsurfaces, and we then show that once this property is achieved for a pair $(S_{\A'}, X)$, it persists for the images of the pair under any further stabilizations.

 First, we observe that as ${P^\A}$ and $h({P^\A})$ are defined by cutting ${S^\A}$ along arcs, it is immediate that ${S^\A}$ may be built up from ${P^\A}$ and $h({P^\A})$. This also implies that if $k>0$, then $|\partial {P^\A} \cap \partial {S^\A}|\geq 2$.
 
According to Lemma~\ref{lem:builtup}, $S^\A$ can be built up   from $S^\A\setminus {P^\A}$ exactly when the boundary of each component of ${P^\A}$ intersects the boundary of $S^\A\setminus {P^\A}$ in at least two intervals. When we stabilize, each added one-handle becomes part of ${P^\A}$, so we may stabilise until ${P^\A}$ is connected.  In the case, the observation of the previous paragraph ensures that $S^\A$ can be built up from $S^\A\setminus {P^\A}$ whenever $k>0$. 
  
Now consider a map $h^\A: {P^\A} \rightarrow {S^\A}$ that is the identity near  $\partial {P^\A} \cap \partial S^\A$.  After stabilizing, we are free to assume that   $h({P^\A})$ is also  connected, so the stabilized ${S^\A}$ can built up from the stabilized ${S^\A} \setminus h({P^\A})$ as above.

Finally, we consider the case $k=0$, when the entire boundary consists of spheres, each of which has the trivial foliation with two elliptic points.  Suppose first that $|\partial M|>1$.  Since the topology of the page remains constant for all $t$, the page is connected and the unique component of $P^\A$ satisfies $|\partial  P^\A|\cap |\partial (S^\A\setminus  {P^\A} )|=|\partial M|>1$.
\end{proof}

Now we prove Proposition~\ref{prop:disj}, the claim that any foliated open book may be positively stabilized to be sorted:

\begin{proof}[Proof of Proposition ~\ref{prop:disj}]

 Recall from Section~\ref{sec:sort} that the sorting arcs $\gamma_k^\pm$ on $S_i$ can be understood as the intersection between the page $S_i$ and the (un)stable submanifold of the critical point $h^\pm_k$ with respect to a preferred gradient. If these arcs are disjoint, then the foliated open book is sorted by definition. On the other hand, when the stable and unstable submanifolds of two critical points intersect,  then these arcs might become disconnected on the pages. In the following,  we will  consider sorting arcs in a a non-sorted foliated open book  only on the pages where these arcs are still connected.

 Given an abstract foliated open book  $(\{S_i\},h)$  for a 3-manifold $M$ with foliated boundary, we will show that it admits a sorted stabilisation $(\{S'_i\},h',\{\gamma_k^{\pm'}\})$.

Let $j$ be the  least index such  that the sorting arcs are not disjoint on $S_j$.  Then $\gamma^+_j$ intersects one or more  $\gamma^-$ arcs; let $k$ denote the largest index such that $\gamma^-_k\cap \gamma^+_j\neq \emptyset. $ We will define a stabilization at $S_j$ to remove this intersection;  the process will introduce no new intersections, so inductively, it suffices to prove the proposition.

As in the definition of abstract stabilization,  choose an arc $\gamma$ on $S_j$ that connects to $B$ near the endpoints of $\gamma^+_j$, follows $\alpha_t$ to the endpoints of $\gamma^+_j$, and then runs parallel to $\gamma^+_j$ except for one intersection in the interior  of $\gamma^+_j$. The hypotheses on $\nabla\pit$ ensure that $\gamma$ is disjoint from all sorting arcs while it runs  along $\alpha_t$, and it intersects $\gamma^+_j$  exactly once. See the first picture of Figure \ref{fig:gamma}.

\begin{figure}[h]
\begin{center}
\includegraphics[scale=1]{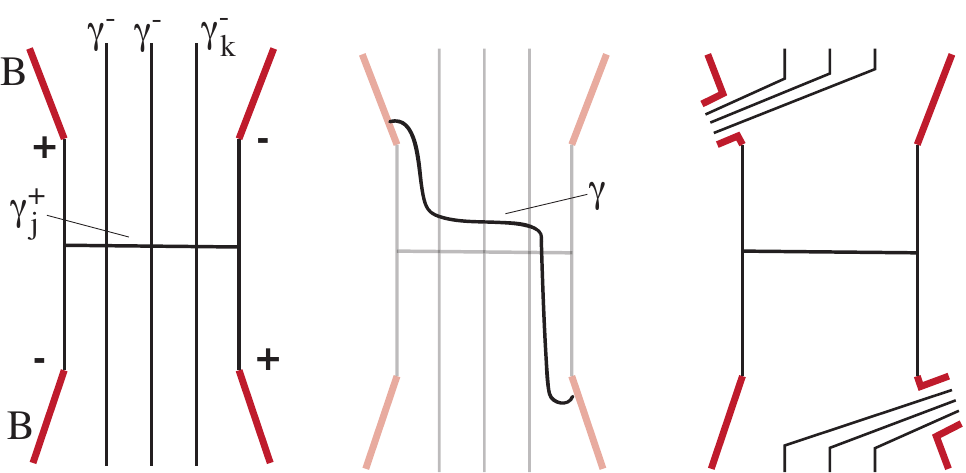}
\caption{Stabilizing the abstract foliated open book along $\gamma$ removes intersections between sorting arcs.}\label{fig:gamma}
\end{center}\end{figure}

Now  perform the stabilization  along $\gamma$.  Consequently,  $\gamma^-_i$  undergoes a right-handed Dehn twist for each $i<j$, removing the targeted intersection point.  Note that this may remove other intersection points, as well, but it cannot introduce new intersections between sorting arcs, so it may be repeated until the foliated open book is sorted. 
\end{proof}

\subsection{Preservation of the contact structure}\label{sec:fobpobequiv}
In this section we prove that the operations turning partial to foliated open books and vice versa preserve the contactomorphism-type of the associated contact structures. 
\begin{proposition}\label{prop:fobpobequiv}
Suppose that the partial open book $\mathcal{P}=(S,P,h)$ is sufficiently stabilized and let $\A(\mathcal{P})=(\{S_i^\mathcal{P}\},h^\mathcal{P}))$ denote an induced foliated open book.
Then $(M(\mathcal{P}),\xi(\mathcal{P}))$ is contactomorphic to  $(M(\mathcal{A}(\mathcal{P})),\xi(\mathcal{A}(\mathcal{P}))$. \end{proposition} 
As partial open books only define contact structures up to contactomorphism and  gluing ``$I$-invariant'' contact structures to the boundary, this statement is as strong as it could be with the given definitions.

The converse statement is as follows:
\begin{proposition}\label{prop:pobfobequiv}
Suppose that the foliated open book ${\A}=(\{S_i\},h)$ is sufficiently stabilized and let $\mathcal{P}({\A})=(S^{\A},P^{{\A}},h^{\A})$ denote the associated partial  open book. Then $(M({\A}),\xi({\A}))$ is contactomorphic to  $(M(\mathcal{P}({\A})),\xi(\mathcal{P}({\A}))$.
\end{proposition} 
We first prove Proposition \ref{prop:pobfobequiv}:

\begin{proof}[Proof of Proposition \ref{prop:pobfobequiv}] 

 In order to prove the proposition, we cut each of the contact manifolds $\big(M(\A),\xi(\A)\big)$, $\big(M(\mathcal{P}(\A)),\xi(\mathcal{P}(\A))\big)$ into a pair of cornered handlebodies which have convex boundary and isotopic dividing sets after smoothing.  The contact structures supported by the foliated and partial open books, respectively, restrict to each handlebody as the unique tight contact structure compatible with this decoration, so it suffices to show that the smoothed handlebodies coming from each structure have matching dividing curves.

 For $({S^\A}, {P^\A}, {h^\A})$, we consider the handlebody  ${S^\A}\times[\frac{-3}{4}, \frac{-1}{4}]$ and its complement in $M(\mathcal{P(\mathcal{A})})$:
 \[\big({S^\A}\times[\frac{-1}{4}, 0]\big) \cup \big({P^\A}\times[0,1]\big))\cup \big(S\times [-1, \frac{-3}{4}] \big)/\sim\]
Just as in the proof of Theorem \ref{thm:existencecontact} we may embed each of these  inside the product ${S'}\times I$, one of the  handlebodies used to construct the open book for an associated closed manifold.  The tight contact structure on this restricts to a tight contact structure on each of our handlebodies, and we examine the induced dividing curves, applying standard disc decomposition arguments to show that there is a unique tight restriction.  
 
 Near the binding, we proceed exactly as in Proposition~\ref{prop:unique}. Smoothing the edges $\partial {S^\A}\times \{\frac{-3}{4}\}$ or  $\partial {S^\A}\times\{ \frac{-1}{4}\}$  creates a right-turning dividing curve; see the final image in Figure~\ref{fig:boundary1}.  However, when smoothing $(\partial {P^\A}\times \{0\}) \cap({S^\A}\times \{0\})$ or $(\partial {P^\A} \times \{1\}) \cap ( {S^\A} \times \{-1\})$, the dividing curve turns left.  See Figure~\ref{fig:boundary2}.

 In the case of the foliated open book $(\{S_i\}, h)$, we choose the two handlebodies to be the sorted 
 ${}_{\varepsilon}H_{1-\varepsilon}\cong \overline{M}$ and a product ${}_{1-\varepsilon}H_{\varepsilon}=S_0\times I$.  Recall that 
 ${}_{\varepsilon}H_{1-\varepsilon}$ was studied in  Propositio n~\ref{prop:sorthandlebody} (as $\overline{M}$), and we can construct the dividing set on 
  $\partial({}_{\varepsilon}H_{1-\varepsilon})$ as in the proof of Proposition~\ref{prop:unique}. Comparing the result to Figure~\ref{fig:boundary2}, one sees that the resulting dividing sets are isotopic
 \begin{figure}[h!]
\begin{center}
\includegraphics[scale=0.9]{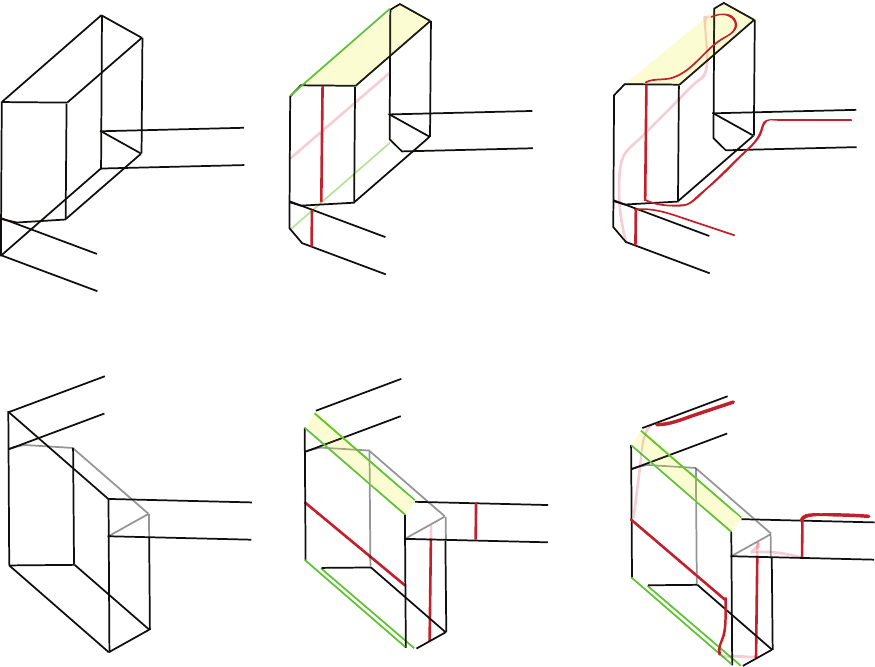}
\caption{ Smoothing the handlebodies in $({S^\A}, {P^\A}, {h^\A})$. }\label{fig:boundary2}
\end{center}
\end{figure}
to the dividing set on 
 ${}_{\varepsilon}H_{1-\varepsilon}$, as desired. 
In the case of each of the product handlebodies, the standard proof applies.   
\end{proof}

\begin{proof}[Proof of Proposition~\ref{prop:fobpobequiv}]
Recall from the  discussion  after   Lemma \ref{lem:pobtofob} that if we start with a sufficiently stabilized partial open book $\mathcal{P}$ and construct the foliated open book ${\A}(\mathcal{P})=(\{S^\mathcal{P}_i\},h^\mathcal{P})$, we can choose the cutting arcs and their order so that the foliated open book is automatically sorted at $S_0^{\mathcal{P}}$.  The construction of Section \ref{sec:pobinfob} then recovers  the partial open book $\mathcal{P}$ without further stabilization, so the result follows from Proposition \ref{prop:pobfobequiv}.
 \end{proof}
 
\subsection{Existence of Foliated Open Books}\label{sec:existfob}
In this subsection we  continue working with the notation established in Section~\ref{sec:sort}.   Beginning with a sufficiently sorted Morse foliated open book, recall  that the subsurfaces $S_t\setminus R_t$ are diffeomorphic for all $t$ and that the associated  partial open book is defined by considering  ${R_0} \subset S_0=S^\A$.}

In fact, all  of the $R_t$ can be read off from $\F_\pit$:
\begin{lemma}\label{lem:Rt} For a sufficiently stabilised {Morse} foliated open book ${\M}=(B, \pi, \F_\pit, {\nabla \pit})$  the diffeomorphism-type of the subsurfaces $R_t$  depends only on $\F_{\pit}$.
\end{lemma}

The proof of this lemma is essentially a generalisation of Lemma~\ref{lem:rplus}:

\begin{proof} 
Fix a regular time $t$. 
Since $\nabla \pit$ is preferred, each leaf $I$ of $\pit^{-1}(t)$ contains disjoint subintervals $I_\pm\subset I$ such that $e_\pm\in I_\pm$ and $I\cap\gamma_k^\pm\subset I_\pm$. Define the subsets $R_t^\pm$ of $R_t$ as follows:
\[R_t^+=N(\bigcup_{e_+} I_+\cup\bigcup_{\pit(h_j)\in[t,1]}\gamma_j^+) \qquad\text{ and }\qquad R_t^-=N(\bigcup_{e_-} I_-\cup\bigcup_{\pit(h_j)\in[0, t]}\gamma_j^-).\]
Thus $R_t\setminus (R_t^+\cup R_t^-)$ is just a union of strips $S_I=N(I\setminus (I_+\cup I_-))$.

As a next step we will push these subsets $R_t^\pm$ to the boundary of $M$. Let $h_j$ be an index 2 critical point with $\pit(h_j)=t_j^*\in[{t}, 1]$. Then  $W^s(h_j)\subset \pi^{-1}[0, t_j^*]$ is a cornered disc with boundary $\gamma_j^+\cup w^s(h_j)'$, where again $w^s(h_j)'$ denotes the stable manifold of $h_j$  with respect to $\nabla\pit$ on $\pit^{-1}[0, t_j^*]$.  (This employs the same slight abuse of notation as in Section~\ref{sec:sort}.)  Then we can isotope $\gamma_j^+$ through $W^s(h_j)$ to $w^s(h_j)'$.  In the neighbourhood of $\partial M$ we can ``turn up'' the half-neighborhood of $I_+$ in $S_i$ into an upper half-neighborhood $N_+(I_+)\subset \pit^{-1}[t,t+\varepsilon]$,
keeping the corners. These two isotopies combine to an isotopy of $R_t^+$ into 

\[Q_t^+=\bigcup_{e_+} N_+ (I_+)\cup\bigcup_{\pit(h_+)\in[t,1]}N(w^s(h_+)').\]
Observe that the orientation of $Q_t^+$ agrees with the orientation of $\partial M$. 
Similarly, we can push $R_t^-$ into 
\[Q_t^-=\bigcup_{e_-} N_- (I_-)\cup\bigcup_{\pit(h_-)\in[0,t] }N(w^u(h_-)'),\]
where 
$N_-(I_-)$ is a lower half-neighborhood of $I_-$ in $\pit^{-1}[t-\varepsilon,t]$. The orientation of $Q_t^-$ is opposite to the orientation of $\partial M$. 

\begin{figure}[h!]
\begin{center}
\includegraphics[scale=0.55]{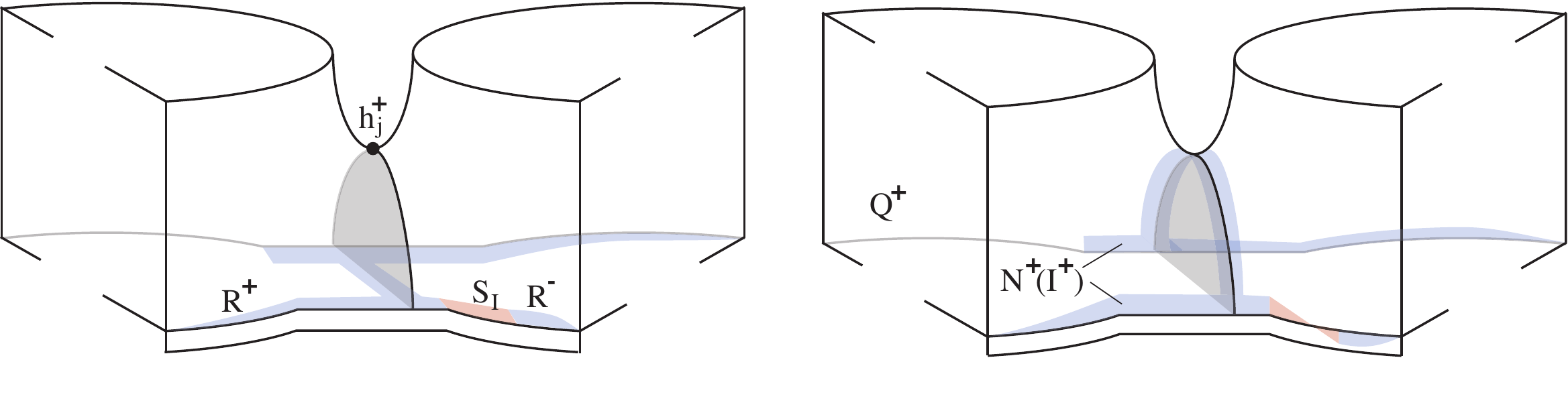}
\caption{$R^\pm_t\subset S_t$ may be isotoped through $M$ to $Q^\pm_t\subset \partial M$ (blue) and connected by twisted bands (red).}\label{fig:rtoq}
\end{center}
\end{figure}

To recover the surface $R_t$, we must isotope the strips $S_I$ to connect $Q_t^\pm$. This can be done via a twist that keeps $\pit^{-1}(t)\cap S_I$ fixed. See Figure~\ref{fig:rtoq}. This construction allows us to reconstruct $R_t$ from the foliation $\F_\pit$,  proving the statement. 

For the critical values,  $R_{t^*}$ is the degenerate surface between  $R_{t^*-\varepsilon}$ and $R_{t^*+\varepsilon}$.
\end{proof}

Recall that the union of the $R_t$ {in $\overline{M}$}, denoted $R$, is a $\nabla \pi$-invariant neighborhood of  the union of $\partial \overline{M}{\setminus (B\times I)}$ and the stable and unstable critical submanifolds. Since each critical submanifold is a cornered disc with half its boundary on $\partial M$, the submanifold $R$ can be compressed inside $\overline{M}$ to a neighborhood of $\partial \overline{M}$ and $R_t\cap\partial \overline{M}$ recovers the foliation $\F_\pit$ cut open at $t=0$. 
As usual, we  choose regular times separating the critical points $0=t_0<t_1^*<\cdots<t_{2k}^*<t_{2k}=1$, and we write $R_{i}=R_{t_i}$. Then up to diffeomorphism, the set $\{R_i\}$ contains all the information about $R$ and $R_t$.   In fact, this is the data needed to  describe a  \textit{partial foliated open book}, which is  an amalgamated foliated- and partial open book; we examine these objects and their applications in a future paper \cite{LV2}.

As motivation  for the following, imagine trying to recover $M(\A)$ from $M(\mathcal{P}(\A))$, where $\mathcal{P}(\A)$ is the partial open book associated to the sufficiently stabilized $\A$.
In the following, we try to mimic the gluing required for this,  but in the case of a partial open book  not necessarily obtained from a foliated open book. More precisely, we will describe  how to glue $R$ 
along $\partial R \setminus \partial \overline{M}$ to the boundary of $S\times [-1,0]\cup P\times [0,1]/\sim$ \footnote{Note that here, and in the rest of the paper, we construct the manifold from the data of a partial open book using the equivalences  $(x,1)\sim (h(x),-1)$ and $(x,t)\sim (x,t')$ for $x\in\partial S$ and $t,t'\in[-1,0]$. This  yields a  manifold diffeomorphic to the one formed by further collapsing $\partial P\setminus \partial S \times [0,1]$ to a collection of intervals.} for some partial open book $(S,P,h)$ while respecting the foliations coming from $R_t$ and $P\times \{t\}$. To move forward,   we examine a further property of the original gluing which will be a key to the general case.

The dividing curve $\Gamma(\mathcal{P}(\A))$ on $\big(M(\mathcal{P}),\xi(\mathcal{P}),\Gamma(\mathcal{P})\big)$ is obtained as $(\partial S^\A\setminus \partial P^\A)\cup (\partial P^\A\setminus \partial S^\A)$. Thus $\partial S^\A\setminus \partial P^\A$ intersects each component $\Gamma_i$ of $\Gamma$ in  $|\Gamma_i\cap \pit^{-1}(0)|$ intervals. The next lemma ensures that we may prescribe these intersections for a partial open book:

\begin{lemma} Let $(M,\xi,\F_\xi)$ be a contact manifold and assume that $\F_\xi$ is strongly topologically conjugate to an open book foliation $\F_{\pit}$ without circle leaves. Then there is a partial open book $(S,P,h)$ supporting $\xi$ so that $\partial S\setminus \partial P$ intersects each component $\Gamma_i$ of $\Gamma$ in $|\Gamma_i\cap \pit^{-1}(0)|$  intervals. 
\end{lemma}\label{lem:pobwithni}

The proof of this lemma is a strengthening of the proof of Lemma~\ref{lem:nocircle}:
\begin{proof}
Recall from Theorem ~1.1 of \cite{HKM} that  the construction of a partial open book requires a polygonal decomposition $K'$ of $\partial M$,   where  $K'$ is a Legendrian graph and the boundary of each 2--cell of $\partial M\setminus K'$ intersects $\Gamma$ in exactly two points. Any such decomposition (with some additional choices explained in \cite{HKM}) gives a partial open book decomposition for $(M,\xi,\Gamma)$ in which $\partial S\setminus \partial P$ intersects each component $\Gamma_i$ of $\Gamma$ in $|\Gamma_i\cap K'|$   intervals.

Define $K'=G_{\scriptscriptstyle{++}}\cup G_{\scriptscriptstyle{--}}\cup \pit^{-1}(0)$. Then we clearly have $|\Gamma_i\cap K'|$=$|\Gamma_i\cap \pit^{-1}(0)|$. So we  need only verify that $K'$ is  the 1--skeleton of a polygonal decomposition. For this we need  each component of $\partial M\setminus K'$ to be a disc whose boundary intersects $\Gamma$ in two points. 

Since the graph $G_{\scriptscriptstyle{\pm\pm}}$ is the spine of $R_\pm$, the components of $R_\pm\setminus G_{\scriptscriptstyle{\pm\pm}}$ are all annuli. And as each component of $\Gamma$  intersects  $\pit^{-1}(0)$, the components of $R_\pm\setminus (G_{\scriptscriptstyle{\pm\pm}}\cup \pit^{-1}(0))$ are rectangles $D_\pm$ with a pair of opposite sides in $\pit^{-1}(0)$ and the other two sides in $\Gamma$ and $G_{\scriptscriptstyle{\pm\pm}}$. The domains of $\partial M\setminus K'$ are obtained by gluing together two such rectangles $D_+$ and $D_-$ along their common boundary in $\Gamma$. This gives a polygonal domain $D$ intersecting $\Gamma$ in two points, as required. 
\end{proof}

Now we are ready to prove the existence result stated as Theorem~\ref{thm:fobexistence}:
\begin{proof}[Proof of Theorem~\ref{thm:fobexistence}] 
 Let $\Gamma$ be the dividing curve for $\F_\xi$ and $\F_\pit$. By Lemma~\ref{lem:pobwithni}, we can find a partial open book $\mathcal{P}=(S,P,h)$ for $(M,\xi,\Gamma)$ so that in the identification $\Gamma=(\partial S\setminus \partial P)\cup (\partial P\setminus \partial S)$, $\Gamma_i$ intersects $\partial P\setminus \partial S$ in exactly $|\pit^{-1}(0)\cap\Gamma_i|$  intervals. 

Construct as above the surfaces $R_i$ for $\F_\pit$. Define \[\widetilde{S_i}^\mathcal{P}=P\cup_{\partial P\setminus \partial S\leftrightarrow \cup \partial R_i\setminus \pit^{-1}(t_i)}R_i,\]
{ where the identification of $\partial P\setminus \partial S$ with $\partial R_i\setminus \pit^{-1}(t_i)$ respects the components of $\Gamma_i$. 
} 

To define the map $\widetilde{h}^\mathcal{P}\colon \widetilde{S_0}^\mathcal{P}\to \widetilde{S_{2k}}^\mathcal{P}$ we have to move between the different identifications of $\widetilde{S_{0}}^\mathcal{P}$ and $\widetilde{S_{2k}}^\mathcal{P}$:
\[\widetilde{S_{0}}^\mathcal{P}=P\cup R_{0}\cong P\cup R_+\cong P\cup (S\setminus P)=S\]
\[\widetilde{S_{2k}}^\mathcal{P}=P\cup R_{2k}\cong P\cup R_-\cong h(P)\cup (S\setminus h(P))=S\]
Under these identifications, $\widetilde{h}^\mathcal{P}\vert_P:=h\colon P\hookrightarrow S$.  Similarly, $\widetilde{h}^\mathcal{P}\vert_{R_0}$ is defined as the embedding into $S$ after we identify $R_0$ with $R_-$.

Since $R_i$ is obtained from $R_{i-1}$ by cutting or gluing, and since the rest of $\widetilde{S_{i}}^\mathcal{P}$, $\widetilde{S_{i}}^\mathcal{P}\setminus R_i=P$ are unchanged, the tuple $\widetilde{\A}(\mathcal{P})=(\{\widetilde{S_{i}}^\mathcal{P}\},\widetilde{h}^\mathcal{P})$ defines a foliated open book.

Notice that $\widetilde{\A}(\mathcal{P})=(\{\widetilde{S_{i}}^\mathcal{P}\},\widetilde{h}^\mathcal{P})$ is sufficiently stabilized and the partial open book associated to it is exactly $\mathcal{P}$. This means that the contact 3-manifold supported by $\widetilde{\A}(\mathcal{P})$ is contactomorphic to $M(\mathcal{\mathcal{P}})$.  Furthermore, the foliation on $\partial M(\widetilde{\A}(\mathcal{P}))$ comes from $R_i$, so $\F_\xi(\widetilde{\A}(\mathcal{P}))$ is strongly topologically conjugate to $\F_\pit$. There is thus a contactomorphism between $M$ and $M(\widetilde{\A}(\mathcal{P}))$
and the image of the abstract foliated open book defines an embedded foliated open book on $(M,\xi,\F_\xi)$, as required.
\end{proof}

\subsection{Proof of Giroux Correspondence}\label{sec:giroux}

The results in this section so far have described a sequence of modifications that may be made to partial and foliated open books via stabilization.  Here, we apply these to prove the Giroux Correspondence for foliated open books which  appeared earlier  as Theorem~\ref{thm:giroux}.

\begin{proof}[Proof of Theorem~\ref{thm:giroux}.]
Suppose that $\A=(\{S_i\},h)$ and $\A'=(\{S_i'\},h'))$ are abstract foliated open books for the same  $(M, \xi, \F_\xi)$. We may assume that each is sufficiently stabilized and that it is possible to define the corresponding partial open books $\mathcal{P}(\A)=({S^\A}, {P^\A}, {h^\A})$ and $\mathcal{P}(\A')=(S^{\A'}, P^{\A'}, h^{\A'})$. As described in the previous section,  we can recover the original abstract foliated open books $\A$ and $\A'$ from the partial open books by  gluing  the $R_i$ determined by the the foliation to the pages. Recall that this gluing is associated to an identification of $\partial P^{\A^*} \setminus \partial S^{\A^*}$ with $\partial R_i\setminus \pit^{-1}(t_i)$. This determines an identification of $\partial P^{\A} \setminus \partial S^{\A}$ with $\partial P^{\A'} \setminus \partial S^{\A'}$, which we will keep in mind as we proceed.

Using \cite{EO},  the partial open books $\mathcal{P}(\A)$ and $\mathcal{P}(\A')$ may be constructed from a Legendrian one-skeleton $K$ and $K'$ as described in \cite{HKM}. These graphs intersect $\partial M$ on $\Gamma$, and the intersections $K^*\cap \Gamma$ are related to $\partial P^{\A^*} \setminus \partial S^{\A^*}$, respectively. Using the identification in the previous section, fix an identification of the endpoints of $K$ and $K'$ and assume that in a small neighborhood $N(\partial M)\cong \partial M\times I$, the graphs $K$ and $K'$ are each of the form $\cup (p_j\times I)$ for some $p_j\in \Gamma$. 
Now apply the  argument of \cite{HKM} to find a common refinement of the Legendrian graph, with the modificaton that $K$ and $K'$ already agree near $N(\partial M)$.  This yields a common stabilization of the two partial open books while preserving the  Legendrian graphs in $N(\partial M)$. This ensures that  
$\partial P \setminus \partial S$ is constant during the process, so at each stage we can  define the a corresponding abstract foliated open book by gluing $R_i$. Since all the added edges were internal,  they may be realized as stabilizations of the foliated open books.  
\end{proof} 

Notice that in the above argument all foliated open books were obtained from partial open books, and thus  were automatically sorted. Consequently, we get a stronger Giroux Correspondence:

\begin{theorem}\label{thm:girouxsorted}[Giroux Correspondence for sorted foliated open books]  Any pair of sorted foliated open books supporting $(M, \xi, \F_\xi)$  are isotopic after a sequence of positive stabilizations. Moreover, each of the intermediate foliated open books in this sequence is sorted. 
\end{theorem}






\bibliographystyle{alpha}

\bibliography{fob}

\end{document}